\documentclass[leqno]{amsart}
\usepackage{paralist}
\usepackage{stmaryrd}
\usepackage{verbatim}
\usepackage{tikz}
\usepackage{afterpage}

\usepackage[linesnumbered,ruled]{algorithm2e}
\usepackage[noend]{algpseudocode}
\let\oldnl\nl
\newcommand{\nonl}{\renewcommand{\nl}{\let\nl\oldnl}}
  
\usepackage{tabularx,ragged2e,booktabs,caption}
\usepackage{graphics} 
\usepackage{epsfig} 
\usepackage{graphicx}  \usepackage{epstopdf}
\usepackage[colorlinks=true]{hyperref}
\usepackage{bm}
\hypersetup{urlcolor=blue, citecolor=red}
\usepackage{hyperref}

\newtheorem{theorem}{Theorem}[section]

\newtheorem{lemma}[theorem]{Lemma}
\newtheorem{proposition}{Proposition}

\theoremstyle{definition}

\newtheorem{remark}{Remark}

\definecolor{light-gray}{gray}{0.95}
\newtheorem{assumption}{Assumption}

\title[Nonlinear ultrasound focusing] 
{Isogeometric shape optimization for nonlinear ultrasound focusing}
\author[Markus Muhr, Vanja Nikoli\' c, Barbara Wohlmuth and Linus Wunderlich]{}


\email{\href{mailto:muhr@ma.tum.de}{muhr@ma.tum.de}}
\email{\href{mailto:vanja.nikolic@ma.tum.de}{vanja.nikolic@ma.tum.de}}
\email{\href{mailto:wohlmuth@ma.tum.de}{wohlmuth@ma.tum.de}}
\email{\href{mailto:linus.wunderlich@ma.tum.de}{linus.wunderlich@ma.tum.de}}

\thanks{$^*$ Corresponding author: Vanja Nikoli\' c}

\begin{document}
\maketitle

\centerline{\scshape Markus Muhr, Vanja Nikoli\' c$^*$, Barbara Wohlmuth and Linus Wunderlich}

\medskip
{\footnotesize
   \centerline{Technical University of Munich, Department of Mathematics, Chair of Numerical Mathematics}
   \centerline{Boltzmannstra\ss e 3, 85748 Garching, Germany}
}

\bigskip

\begin{abstract}
The goal of this work is to improve focusing of high-intensity ultrasound by modifying the geometry of acoustic lenses through shape optimization. The shape optimization problem is formulated by introducing a tracking-type cost functional to match a desired pressure distribution in the focal region. Westervelt's equation, a nonlinear acoustic wave equation, is used to model the pressure field. We apply the optimize first, then discretize approach, where we first rigorously compute the shape derivative of our cost functional. A gradient-based optimization algorithm is then developed within the concept of isogeometric analysis, where the geometry is exactly represented by splines at every gradient step and the same basis is used to approximate the equations. Numerical experiments in a $2$D setting illustrate our findings.  
\end{abstract}
\section{Introduction}
\begin{figure}[h] 
	\centering
	\begin{tikzpicture}[scale=0.65,>=stealth,label distance=.2cm]
	
	\colorlet{Blau}{blue!50!cyan!75!white};
	\colorlet{Gruen}{green!80!black!75!white};
	\colorlet{Oranje}{orange!75!white};
	\colorlet{ltblue}{blue!20!cyan!75!white};
	\colorlet{ltgray}{black!30!white};
	\colorlet{lltgray}{black!10!white};
	\colorlet{llltgray}{black!2!white};
	
	\coordinate (L1) at (0,0);
	\coordinate (L2) at (4.6,0);
	\coordinate (B1) at (0,-4);
	\coordinate (B2) at (4.6,-4);
	\coordinate (D1) at (0.5,-4);
	\coordinate (D2) at (4.1,-4);
	\coordinate (D3) at (0.5,-4.1);
	\coordinate (D4) at (4.1,-4.1);
	\coordinate (T1) at (0,4);
	\coordinate (T2) at (4.6,4);
	\coordinate (F)  at (2.3,2.5);
	
	\coordinate (w1) at (0,-4);
	\coordinate (w2) at (0,-0.8);
	\coordinate (w3) at (0,0);
	
	\coordinate (W1) at (1.15,-4);
	\coordinate (W2) at (1.15,-0.8);
	\coordinate (W3) at (1.3, 0);
	
	\coordinate (WW1) at (3.45,-4);
	\coordinate (WW2) at (3.45,-0.8);
	\coordinate (WW3) at (3.3, 0);
	
	\coordinate (WWW1) at (2.3,-4);
	\coordinate (WWW2) at (2.3,-0.8);
	\coordinate (WWW3) at (2.3, 0);
	
	\coordinate (ww1) at (4.6,-4);
	\coordinate (ww2) at (4.6,-0.8);
	\coordinate (ww3) at (4.6,0);

	\filldraw[llltgray] (-0.5,-4) -- (-0.5,4) -- (5.1,4) -- (5.1,-4) -- cycle;
	\draw[black] (-0.5,-4) -- (-0.5,4) -- (5.1,4) -- (5.1,-4) -- cycle;
	
	\filldraw[lltgray] (L1) arc (-140:-40:3) -- (L2) -- cycle;
	\draw[black,thick] (L1) arc (-140:-40:3);
	\draw[black,thick] (L1) -++ (4.6,0);
	
	\draw[draw=gray] (w1) -- (w2) -- (w3) --(F);
	\draw[draw=gray] (ww1) -- (ww2) -- (ww3) --(F);
	
	\draw[draw=gray] (W1) -- (W2) -- (W3) --(F);
	\draw[draw=gray] (WW1) -- (WW2) -- (WW3) --(F);
	\draw[draw=gray] (WWW1) -- (WWW2) -- (WWW3) --(F);
	
	\filldraw[gray,ultra thick] (D1) --  (D2) -- (D4) -- (D3) -- cycle;
	\draw[black,ultra thick] (B1) --  (B2);
	\filldraw[black] (F) circle (2pt);
	
	\draw[black, thick, ->] (0.5,-3.8) -- (0.5,-3);
	
	\node[shape=circle,inner sep=2pt] at (2.8,-0.5) {lens};
	\node[shape=circle,inner sep=2pt] at (3.75,2.5) {focal point};
	\node[shape=circle,inner sep=2pt] at (0.1,1) {fluid};
	\node[shape=circle,inner sep=2pt] at (2.3,-4.4) {excitator};

	\coordinate (TR1) at (8,-2.5);
	\coordinate (TR2) at (14,-2.5);
	\coordinate (TB1) at (8,-4);
	\coordinate (TB2) at (14,-4);
	\coordinate (TT1) at (8,4);
	\coordinate (TT2) at (14,4);
	\coordinate (FF)  at (11,2.5);
	
	\filldraw[llltgray] (TR1) -- (TT1) -- (TT2) -- (TR2) arc (-60:-120:6) -- cycle;
	\draw[black] (TR1) -- (TT1) -- (TT2) -- (TR2) arc (-60:-120:6) -- cycle;
	
	\filldraw[lltgray] (TR1) -++ (0,0.2) arc (-120:-60:6) -++ (0,-0.2) arc (-60:-120:6);
	\draw[black] (TR1) -++ (0,0.2) arc (-120:-60:6) -++ (0,-0.2) arc (-60:-120:6);
	\draw[black, ultra thick] (8,-2.3) arc (-120:-60:6);
	\draw[black, thick] (8.5,-2.75) -++ (0.08,0.17);
	\draw[black, thick] (9.25,-3.05) -++ (0.065,0.17);
	\draw[black, thick] (10.1,-3.25) -++ (0.038,0.19);
	\draw[black, thick] (11,-3.28) -++ (0.00,0.17);
	\draw[black, thick] (22-8.5,-2.75) -++ (-0.08,0.17);
	\draw[black, thick] (22-9.25,-3.05) -++ (-0.065,0.17);
	\draw[black, thick] (22-10.1,-3.25) -++ (-0.038,0.19);

	\draw[gray] (8.26,-2.4) --(FF);
	\draw[gray] (8.9,-2.7) --(FF);
	\draw[gray] (9.75,-2.95) --(FF);
	\draw[gray] (10.5,-3.05) --(FF);
	\draw[gray] (22-8.26,-2.4) --(FF);
	\draw[gray] (22-8.9,-2.7) --(FF);
	\draw[gray] (22-9.75,-2.95) --(FF);
	\draw[gray] (22-10.5,-3.05) --(FF);
	
	\filldraw[black] (FF) circle (2pt);
	
	\node[shape=circle,inner sep=2pt] at (9+3.6,2.5) {focal point};
	\node[shape=circle,inner sep=2pt] at (9+0.1,1) {fluid};
	\node[shape=circle,inner sep=2pt] at (11,-3.7) {piezoelectric transducers};
	
	\draw[black, thick, ->] (8.6,-2.4) -++ (0.5,0.95);
	\end{tikzpicture}\\[-15mm]
	
	\caption{Sketches of different ultrasound focusing approaches; \textbf{left:} Focusing by a lens, \textbf{right:} Focusing by an array of transducers placed on a curved surface. \label{fig:Sketches}}
\end{figure}
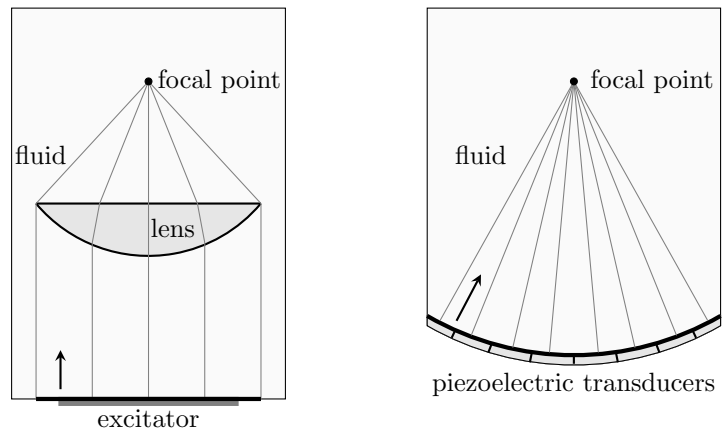

High-intensity focused ultrasound (HIFU) has become a powerful medical tool in recent years and a research area under active development. As a non-invasive therapeutic technology, it is used for breaking up kidney stones by lithotripsy \cite{Lee2, Paterson, Yoshizawa} and it is currently involved in various stages of clinical trials for treatments of cancer in different organs, including the kidney, liver, prostate and brain  \cite{Blana, Maloney, Kennedy1, Kennedy2, Wu}. \\ 
\indent Focusing in  ultrasound applications is achieved geometrically, with lenses or spherical focusing  arrays  of  piezoelectric  transducers (see~\cite{manfred} and Figure~\ref{fig:Sketches} for geometry description). Such setups can naturally benefit from shape optimization: we can find the shape of the lens or the transducer surface which results in a desired pressure distribution around the focal point. The lens design has been quite often determined by heuristic laboratory-based approaches  \cite{Mancini, Neisius, Zhong}. The goal of the present work is to improve the focusing of high-intensity ultrasound by employing a shape optimization approach to modify geometry of a focusing lens.  \\
\indent The use of high-power ultrasound sources enhances nonlinear effects in the wave propagation, inducing the formation of a sawtooth shape of the wave and generation of higher harmonics. When the focused acoustic wave propagates toward the focal point, its amplitude increases in a nonlinear fashion. Shape optimization for problems depending on linear partial differential equations is a well-explored topic; see, for instance, \cite{Henrot, Zolesio} for results concerned with analysis in the continuous setting and works \cite{Haslinger, Laporte, Paganini}, which tackle numerical aspects. Investigations of nonlinear problems are not as common, we refer, for example, to~\cite{Gangl} for the treatment of an optimal design problem in nonlinear magnetostatics and~\cite{Brandenburg} for optimization of instationary Navier-Stokes flow. \\
\indent Results on gradient-based shape optimization for problems depending on second order hyperbolic equations are mainly restricted to the analysis with some numerical results available for linear problems. In \cite{Maestre}, an optimal design problem governed by a linear damped one-dimensional wave equation was investigated. Results on optimal design for the support of the control for the $2$D linear wave equation can be found in \cite{Muench}. In~\cite{NikolicKaltenbacher}, shape sensitivity analysis was provided for the problem of optimizing the shape of the acoustic lens; the case of ultrasound focusing by a concave array of transducers was treated analytically in~\cite{KaltenbacherPeichl}. First steps toward numerical shape optimization for linear ultrasound were made in~\cite{Veljovic}.  Numerical treatment of the shape optimization problems for \emph{nonlinear} ultrasound has  to the authors' best knowledge been missing up to now. \\
 \indent Another goal of the present work is to investigate shape optimization subject to nonlinear waves within the concept of isogeometric analysis. Isogeometric analysis (IGA) was introduced by Hughes et al. in~\cite{hughes:09, hughes:05}. The   idea behind IGA is that the non-uniform rational B-spline (NURBS) basis used to model geometry is also used as a basis for the desired solution. The idea of using spline functions for the approximation of partial differential equations can even be found in earlier works, see, e.g.,~\cite{Hoellig} and references therein. However with IGA, geometry and analysis always go hand in hand, which is convenient especially for shape optimization. Numerical isogeometric shape optimization has been by now applied to a number of linear problems arising in fluid mechanics~\cite{Nortoft}, electromagnetic scattering~\cite{Nguyen} and structural mechanics \cite{Cho,Kiendl,Wall}, while foundations of a mathematical framework for shape optimization problems within isogeometric analysis were laid out in \cite{Simeon,Simeon2}. In comparison to classical shape optimization methods, shape optimization based on the isogeometric approach offers several   advantages. Even complex geometries can be represented precisely by multi-patch NURBS geometries and the output of the optimization can then directly be exported to a computer aided design (CAD) system. In addition, NURBS-bases can easily achieve higher order of global regularity than classical finite element bases.

\indent The rest of the paper is organized as follows. In Section~\ref{Section_ProblemSetting}, we formulate the shape optimization problem together with its underlying nonlinear wave equation. Section~\ref{Section_ShapeDer} contains basics of shape calculus and provides a result on continuity of the state with respect to domain perturbations that will be needed to formally compute the derivative. Section~\ref{Section_ComputationShapeDer} is then dedicated to the calculation of the shape derivative of the cost function. In Section~\ref{Section_IGA}, we give a brief overview of the isogeometric analysis and then introduce the numerical schemes we use to solve our state and adjoint problems. The shape optimization algorithm is presented in Section~\ref{Section_ShapeOptAlgorithm}. Finally, in Section~\ref{Section_NumericalExperiments}, we illustrate the efficiency of our algorithm through several numerical experiments. 
\section{Problem setting} \label{Section_ProblemSetting}
\begin{figure}[h]
	\includegraphics[trim = 3cm 9cm 3cm 10cm, clip, scale = 0.5]{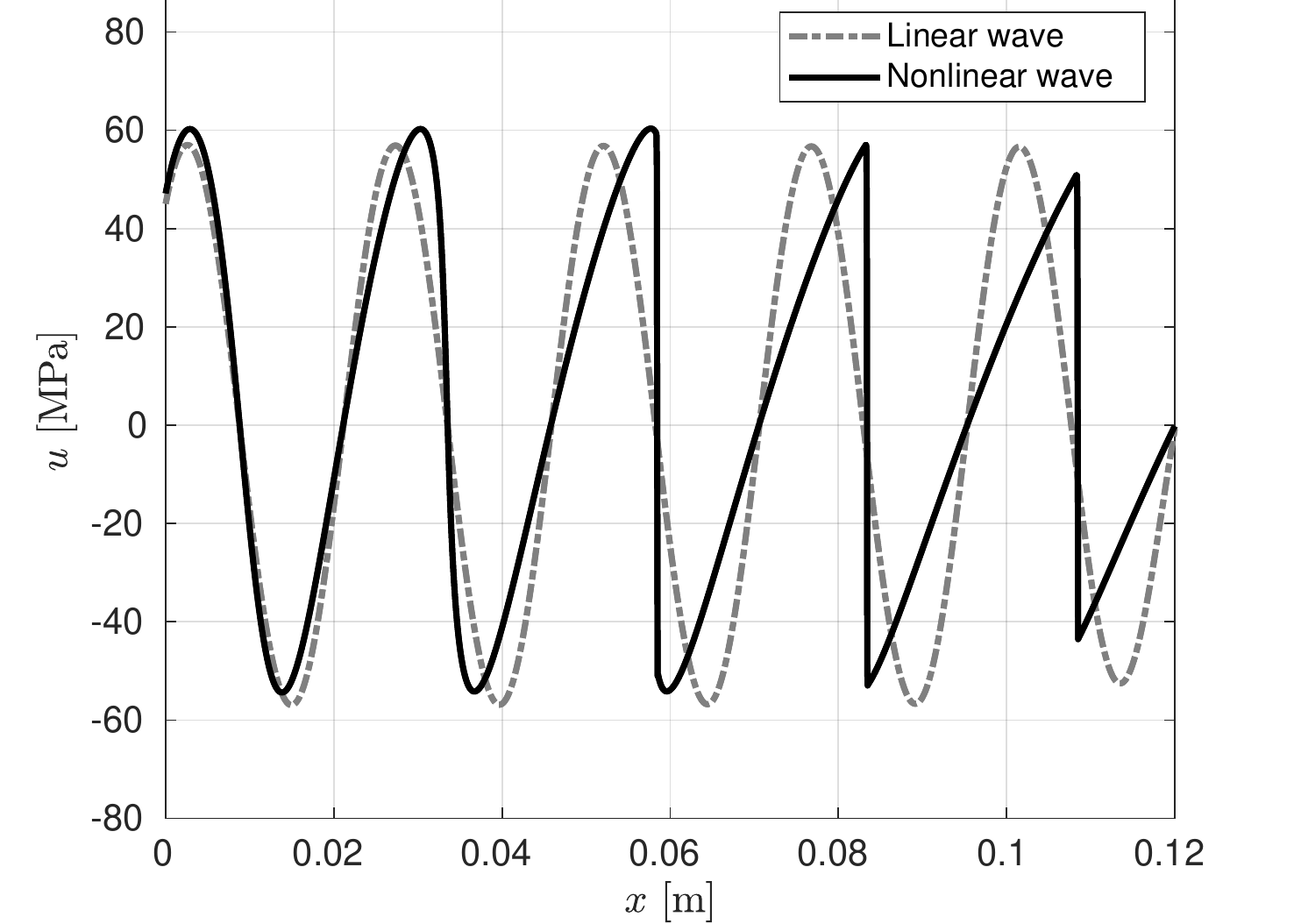}
	\caption{Formation of a saw-tooth pattern of a nonlinear wave in a channel setting with a sinusoidal excitation signal, compared to a linear wave propagation. \label{fig:Sawtooth}}
\end{figure}

\indent To properly capture the nonlinear propagation behavior (see Figure~\ref{fig:Sawtooth}) we have to employ a model of \emph{nonlinear} acoustic for computing the pressure field. There is a number of so-called weakly nonlinear models in nonlinear acoustics obtained as a one-equation approximation of the compressible Navier-Stokes system. In the present work, we will employ a classical model, the Westervelt equation, to compute the acoustic pressure $u$:
\begin{align*} 
\ddot{u}-c^2 \Delta u -b\Delta \dot{u}=2k(\dot{u}^2+u \ddot{u}).
\end{align*}
In the equation, $c$ denotes the speed of sound in the medium and $b$ the diffusivity. The strong $b$-damping in the model accounts for the energy losses that occur due to viscosity and thermal conductivity of the medium. The coefficient $k$ is given by $$\displaystyle k=\frac{1}{\varrho c^2}(1+\frac{1}{2}B/A),$$ where $\varrho$ is the mass density. Here the parameter of nonlinearity $B/A$ accounts for the nonlinear pressure-density relation up to second order and indicates the strength of the nonlinearity of a medium. Values of $B/A$ for different nonlinear fluids and gases are given, for example, in~\cite{Rossing}. A detailed derivation of the model from the governing equations and an overview of well-established acoustic models can be found, for instance, in \cite{Crighton, manfred}. 

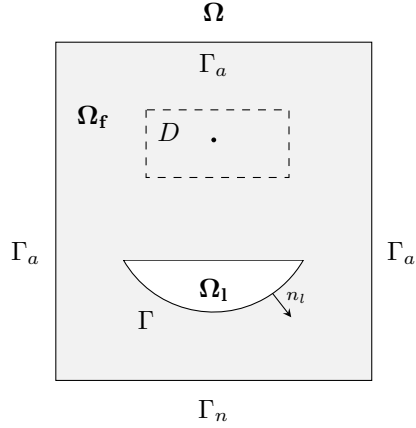
\begin{figure}[h] 
\centering
\begin{tikzpicture}[>=stealth] 
    \draw [fill=light-gray]  (-2.1,-2.1) rectangle (2.1, 2.4);
    \draw [dashed] (-0.9, 1.5) rectangle (1, 0.6);
    \draw[fill] (0, 1.1) circle [radius=0.025];
    \begin{scope}
        \clip (-1.2,-0.5) rectangle (1.2,-1.7);
        \draw [fill=white] (-1.2,-0.5) arc (-150:-30:1.38) -- cycle;
    \end{scope}
    \draw [black,->] (0.78,-0.94) -- (1.03, -1.26);
    \node at (1.1, -1) {\scriptsize $n_l$};
    \node at (0, -0.9) {$\mathbf{\Omega_l}$};
    \node at (-0.6, 1.2) {$D$};
    \node at (-0.9, -1.3) {$\large \Gamma$};
    \node at (-1.6, 1.4) {$\mathbf{\Omega_f}$};
    \node at (0, -2.5) {$\Gamma_n$};
    \node at (0, 2.1) {$\Gamma_a$};
    \node at (0, 2.8) {$\mathbf{\Omega}$};
    \node at (2.5, -0.4) {$\Gamma_a$};
    \node at (-2.5, -0.4) {$\Gamma_a$};
\end{tikzpicture} 
\caption{Domain $\Omega$, consisting of lens $\Omega_l$ and fluid $\Omega_f$.  \label{figure_domain}}
\end{figure}

\indent Typically in ultrasound applications the acoustic lens is immersed in an acoustic fluid. Let $\Omega \subset \mathbb{R}^d$, $d \in \{1,2,3\}$, be a fixed bounded domain with Lipschitz  boundary $\partial \Omega$ and $\Omega_{l}$ its subdomain, representing the lens, such that $\Omega_{l} \subset \Omega$ and $\Omega_{l}$  has Lipschitz boundary $\partial \Omega_l$ (see Figure~\ref{figure_domain}). We denote by $\Omega_{f}=\Omega \setminus \overline{\Omega}_{l}$ the part of the domain representing the fluid region. $\Gamma=\partial \Omega_l \cap \partial \Omega_f$ then represents the interface between $\Omega_l$ and $\Omega_f$. \\
\indent The magneto-mechanical excitation will be modeled by Neumann boundary conditions on the part of the boundary $\partial \Omega$ denoted $\Gamma_n$:
\begin{align*}
 \displaystyle \frac{\partial u}{\partial n}=g \, \text{  on  }  \Gamma_n,
\end{align*}
where $n$ represents the outer unit normal to $\partial \Omega$. \\
 \indent An important issue in acoustics is how to truncate an open domain problem to a bounded domain $\Omega \subset \mathbb{R}^d$, $d=1,2,3$, without causing spurious reflections of the waves. In this work we will employ the zero order Enquist-Majda absorbing boundary conditions \cite{Engquist}:
\begin{align*} 
 \displaystyle \dot{u}+c\frac{\partial u}{\partial n}=0 \, \text{  on  }  \Gamma_a,
\end{align*} 
where $\Gamma_a =\partial \Omega \setminus \Gamma_n$ denotes the part of the boundary with absorbing conditions prescribed. \\
\indent In what follows, we will denote by $n_l$ the unit outer normal to the lens domain $\Omega_l$. Restriction of a function $v$ to $\Omega_{i}$ will be denoted by $v_{i}$, where $i \in \{l, f\}$, and $\llbracket v \rrbracket =v_l-v_f$ will stand for the jump over the interface $\Gamma$. The acoustic pressure can be calculated via the following initial boundary value problem:
\begin{align} \label{state_constraint_strong}
\begin{cases}
\ddot{u}-\text{div}(c^2 \nabla u)
-\text{div}(b \nabla \dot{u}) =2k(\dot{u}^2+u \ddot{u}) \, \text{  in  } \Omega_{l} \cup \Omega_{f}, \vspace{2mm} \\
\llbracket u \rrbracket=0 \, \text{  on  }  \Gamma, \vspace{2mm} \\
 \displaystyle \Bigl \llbracket  c^2 \frac{\partial u}{\partial n_l}+b\frac{\partial \dot{u}}{\partial n_l} \Bigr \rrbracket=0 \, \text{  on  }  \Gamma, \vspace{2mm} \\
 \displaystyle \frac{\partial u}{\partial n}=g \, \text{  on  }  \Gamma_n, \vspace{2mm}\\
 \displaystyle \dot{u}+c\frac{\partial u}{\partial n}=0 \, \text{  on  }  \Gamma_a, \vspace{2mm}\\
(u,\dot{u})\vert_{t=0}=(0, 0).
\end{cases}
\end{align}
The coefficients are assumed to be piecewise constant functions, i.e.
\begin{equation} \label{coeff_assumptions}
\begin{aligned}
& c(x)=\chi_{\Omega_{l}}(x)c_l+(1-\chi_{\Omega_{l}}(x))c_f, \\
& b(x)=\chi_{\Omega_{l}}(x)b_l+(1-\chi_{\Omega_{l}}(x))b_f, \\
& k(x)=\chi_{\Omega_{l}}(x)k_l+(1-\chi_{\Omega_{l}}(x))k_f, \\
& c_{l, f} > 0, \ b_{l,f} >0, \ k_{l,f} \in \mathbb{R}, \ k_l^2+k_f^2 \neq 0,
\end{aligned}
\end{equation}
where $\chi_{\Omega_{l}}$ is the indicator function of $\Omega_l$. Throughout the paper, we assume that $t \in [0,T]$, with a finite time horizon $T$. In order to simplify the presentation and in accordance with ultrasound applications, we have set the initial conditions to zero. \\
\\
\indent We will study the variational form of~\eqref{state_constraint_strong}, given by
\begin{align} \label{state_constraint_weak}
\begin{cases} 
\text{Find }u \in \mathcal{X} =H^{2}(0, T; L^2(\Omega)) \cap H^1(0,T;H^1(\Omega)),\text{ such that} \vspace{1.5mm} \\
\displaystyle \int_0^T \int_{\Omega} \{\ddot{u} \phi+c^2\nabla u \cdot \nabla \phi + b\nabla \dot{u} \cdot \nabla \phi 
-2k(\dot{u}^2+u\ddot{u}) \phi\} \, dx \, ds \vspace{1.5mm} \\
+\displaystyle \int_0^T \int_{\Gamma_a} (c \dot{u}+\frac{b}{c}\ddot{u})\phi \, dx \, ds- \int_0^T  \int_{\Gamma_n} (c^2g+b \dot{g})\phi  \, dx \, ds=0, \vspace{1.5mm}\\
\text{for all test functions} \ \phi \in X=L^2(0,T;H^{1}(\Omega)) \ \text{and} \vspace{1.5mm}\\
 (u,\dot{u})\vert_{t=0}=(0, 0).
\end{cases}
\end{align} 
\subsection{On the well-posedness of \eqref{state_constraint_weak}} Westervelt's equation is a strongly damped quasilinear wave equation. An interesting feature of the model becomes more evident when rewritten as
 \begin{align*}
(1-2ku)\ddot{u}-c^2 \Delta u -b\Delta \dot{u}=2k\dot{u}^2.
\end{align*}
If  $\displaystyle u=1/(2k)$, the equation degenerates, which means that special attention in the analysis has to be paid to obtaining a bound in the $L^\infty$ norm for the pressure to keep it away from $1/(2k)$. Typically, this is resolved by a higher-regularity result and an embedding, e.g. $H^2(\Omega) \hookrightarrow L^\infty(\Omega)$, together with appropriate energy estimates. For details we refer to various papers \cite{KL11, KLV10, KaltenbacherPeichl, MeyerWilke} dealing with well-posedness of this model with constant coefficients under different boundary conditions. \\
\indent Going forward, we assume that problem~\eqref{state_constraint_weak} is well-posed locally in time and space:
\begin{assumption} \label{assumption_wellposedness} We assume that there exists $\bar{m}>0$ and a final time $T>0$ such that under sufficiently regular boundary data $g$ problem~\eqref{state_constraint_weak} has a unique solution $u \in \mathcal{U}$, where
\begin{equation} \label{space_U}
\begin{aligned} 
 \mathcal{U} =\Bigl \{u \in \mathcal{X}  \ \vert& \  u \in  C([0,T]; L^\infty(\Omega)) \cap W^{1, \infty}(0,T; H^1(\Omega)), \vspace{0.3mm}\\
 & \ \displaystyle \text{tr}_{\Gamma_a} \dot{u} \in H^1(0,T; L^2(\Gamma_a)), \ (u,\dot{u})\vert_{t=0}=(0, 0), \vspace{0.3mm}\\
& \ \|\ddot{u}\|_{L^2(0,T;L^2(\Omega))} \leq \bar{m}, \ \|\dot{u}\|_{L^\infty(0,T;L^2(\Omega))} \leq \bar{m}, \vspace{0.3mm}\\
& \  \|\nabla \dot{u}\|_{L^2(0,T;L^2(\Omega))} \leq \bar{m}, \ \|u\|_{L^\infty(0,T; L^\infty(\Omega))} \leq \bar{a} < \frac{1}{2|k|_{L^\infty(\Omega)}}\Bigr \}.
\end{aligned}
\end{equation}
\end{assumption}
\indent Well-posedness of the problem \eqref{state_constraint_weak} with constant coefficients $c, b, k>0$ was shown to hold in~\cite[Theorem 2]{KaltenbacherPeichl}, with a solution $u \in \mathcal{U} \cap L^\infty(0,T; H^{3/2+\tilde{s}}(\Omega))$, $\tilde{s} \in (0, 1/2)$, however well-posedness in the case of piecewise-constant coefficients is to the authors' best knowledge still an open problem. The main difficulty is  that the $H^{3/2+\tilde{s}}$ regularity in space, which is employed to avoid degeneracy, is too high to achieve when the normal derivative of the pressure is allowed to jump, as is the case in our problem. \\
\indent Note that the last condition in \eqref{space_U} is imposed to ensure that the degeneracy of the model is avoided. In practice, this condition is easily fulfilled since the coefficient $k$ is quite small. For instance, in water  $k \approx 1.5 \cdot 10^{-9} \, \displaystyle \text{s}^2\text{m}/\text{kg}$, whereas the maximal pressure values in HIFU applications typically do not exceed $80$ MPa~\cite{Canney, manfred, Neisius}.  This will be reflected in our numerical experiments as well.
\subsection{The shape optimization problem}\label{Section_ShapeOptProblem}
  \indent We will consider a shape optimization problem of the following form \index{cost functional}
 \begin{align} \label{shape_opt_problem}
 \begin{cases}
 \displaystyle \min_{\substack{(u, \Omega_l)}} \quad &J(u, \Omega_l) \\
\text{s.t.} \ \quad  &E(u, \Omega_l)=0, \ \ u \in \mathcal{U}, \ \Omega_l \in \mathcal{O}_{\mathrm{ad}},
 \end{cases}
 \end{align}
with $\mathcal{U}$ defined as in~\eqref{space_U}. For the objective function, we choose an $L^2$-tracking type functional
\begin{align} \label{objective_functional} 
\displaystyle  J(u, \Omega_l)=\,  \int_0 ^T \int_{\Omega} (u-u_{d})^2 \, \chi_{D}\, dx \, ds,
\end{align}
with $D$ being a fixed bounded Lipschitz domain, $\overline{D} \subset \Omega_f$,  where the pressure distributions should match (see Figure~\ref{figure_domain}). The constraint is given by
\begin{align*} 
&\langle E(u,\Omega_l), \phi \rangle_{X^{\star}, X} \vspace{1.5mm}\\
 =&\displaystyle \int_0^T \int_{\Omega} \{\ddot{u} \phi+c^2\nabla u \cdot \nabla \phi + b\nabla \dot{u} \cdot \nabla \phi 
-2k(\dot{u}^2+u\ddot{u}) \phi\} \, dx \, ds \vspace{1.5mm} \\
&+\displaystyle \int_0^T \int_{\Gamma_a} (c \dot{u}+\frac{b}{c}\ddot{u})\phi \, dx \, ds- \int_0^T  \int_{\Gamma_n} (c^2g+b \dot{g})\phi  \, dx \, ds.
\end{align*}
We assume that the desired pressure $u_{d} \in L^2(0,T; L^2(\Omega)) \cap H^1(0,T; (H^1(\Omega))^\star)$, where $(H^1(\Omega))^\star$ denotes the dual space of $H^1(\Omega)$. Together with the assumed regularity of the state, we have $u-u_{d} \in L^2(0,T; L^2(\Omega))$. The set of admissible domains $\mathcal{O}_{\text{ad}}$ can be defined as
\begin{align*} 
\mathcal{O}_{\text{ad}}=\{\ \Omega_l \ \vert \ & \Omega_l \subset \Omega, \ \Omega_l \ \text{is open and Lipschitz} \\
&\text{with uniform Lipschitz constant} \ L_{\mathcal{O}} \ \}.
\end{align*} 
It can be shown that this set is compact with respect to the Hausdorff metric (cf.~\cite[Theorem 2.4.10]{Henrot}).
\begin{remark}
From the point of view of applications it is crucial that the target pressure is reached around the focal point, which is why we demand that the acoustic pressure achieves prescribed values in a fixed subregion $D$. This means that our cost functional is compactly supported. The class of shape optimization problems with an $L^2(D)$-tracking cost functional subject to Poisson equation with either Dirichlet or Neumann data are well-investigated and known to be ill-posed~\cite{Eppler1, Eppler2, Harbrecht1}. The problems were proven to be unstable by considering the shape Hessian and showing that it is not strictly $H^{1/2}$-coercive at the critical domain. In our setting, however, the fixed domain $D$ does not intersect the domain $\Omega_l$ subject to optimization. The numerical experiments we conducted were shown to be effective without any regularization, which is why we do not introduce it here.
\end{remark}
\subsection{The adjoint problem} \label{Section_Adjoint}
\noindent Since we will use the adjoint approach to compute the shape derivative, here we investigate the adjoint problem, given by the very-weak-in-time form 
\begin{align*} 
\langle E_u(u, \Omega_l)\psi, p \rangle_{X^\star, X} =\displaystyle \int_0^T \int_{\Omega}j^{\prime}(u)\, \psi \, dx ds,
\end{align*}
with $\psi \in \mathcal{X}$, $j(u)=(u-u_d)^2 \chi_D$. The weak form of the adjoint problem is then given by
\begin{align} \label{adjoint_problem_weak}
\begin{cases}
\displaystyle \int_0^T \int_{\Omega} \{(1-2ku)\ddot{p} \zeta+c^2\nabla p \cdot \nabla \zeta-b\nabla \dot{p} \cdot \nabla \zeta \} \, dx \, ds \vspace{1.5mm} \\
+\displaystyle \int_0^T \int_{\Gamma_a} (-c \dot{p}+\frac{b}{c}\ddot{p})\zeta\, dx ds = 2\int_0^T \int_{\Omega}  (u-u_d) \chi_D \zeta \, dx \, ds,  \vspace{1.5mm} \\
\text{for all} \ \zeta \in X=L^2(0,T;H^1(\Omega)), \, \text{with} \vspace{1.5mm} \\
(p,\dot{p})\vert_{t=T}=(0, 0).
\end{cases}
\end{align}

\begin{proposition}\label{well_posedness_adjoint} Let Assumption \ref{assumption_wellposedness} on well-posedness of the state problem hold. Then for sufficiently small final time $T>0$ and $\bar{m}>0$, with $\bar{m}$ as in \eqref{space_U},  problem \eqref{adjoint_problem_weak} has a unique solution 
\begin{align*}
p \in \mathcal{P}=\{ \phi \ \vert& \ \phi \in H^2(0,T;L^2(\Omega)) \cap W^{1, \infty}(0,T;H^1(\Omega)), \\
& \ \textup{tr}_{\Gamma_a} \dot{\phi} \in H^1(0,T; L^2(\Gamma_a)) \}.
\end{align*}
\end{proposition}
\begin{proof}
By time reversal $\tilde{p}(t)=p(T-t)$, $\tilde{u}(t)=u(T-t)$, $\tilde{u}_d(t)=u_d(T-t)$, we can obtain a forward problem with initial conditions $(\tilde{p},\dot{\tilde{p}})\vert_{t=0}=(0, 0)$.\\
\indent Well-posedness of a forward problem obtained in that way, but with coefficients assumed to be constant $c, b, k>0$ was proven in  \cite[Corollary 7]{KaltenbacherPeichl}. The proof is based on a fixed-point argument, together with Galerkin-Faedo approximations in space and lower and higher-order energy estimates. By closely inspecting the proof, we see that the well-posedness result still holds with our piecewise constant coefficients assumed to be as in \eqref{coeff_assumptions}.
\end{proof}

\section{Shape calculus} \label{Section_ShapeDer}
\indent In order to define the shape derivative, we take the well-established approach of Murat and Simon~\cite{MuratSimon}, where changing domains are described by transformations of a reference domain. To this end, we introduce a fixed vector field $\Theta \in C^{1,1}(\overline{\Omega},\mathbb{R}^d)$ with $\text{supp}~\Theta \cap (\overline{D} \cup \partial \Omega)=\emptyset$, and a family of transformations $F_\tau:  \overline{\Omega} \mapsto \mathbb{R}^d$ given by perturbations of the identity mapping
\begin{align} \label{perturbation of identity mapping} 
F_\tau=\textup{id}+\tau \Theta, \ \ \ \tau \in \mathbb{R}.
\end{align}
There exists $\tau_0>0$ such that for $|\tau|<\tau_0$, $F_\tau(\Omega)=\Omega$ and  $F_\tau$ is a $C^{1,1}$ diffeomorphism  (cf.~\cite{IKP}). We define the perturbed lens as $$\Omega_{l, \tau}=F_\tau(\Omega_l)$$
and the perturbed lens boundary as 
$$\Gamma_\tau=F_\tau(\Gamma).$$ 
Then $\Gamma_\tau$ is Lipschitz continuous (cf.~\cite[Theorem 4.1]{Hofmann}) and $\Omega_{l,\tau} \subset \Omega$ for $|\tau|<\tau_0$. Due to assumptions on $\Theta$, the boundary of $\Omega$ and subdomain $D$ remain fixed as $\tau$ changes. \\
\indent For $\tau \in (-\tau_0, \tau_0)$, we consider also the transported problem on the domain with the perturbed lens $\Omega_{l, \tau}$
\begin{equation*} 
\langle E(u_{\tau}, \Omega_{l, \tau}), \phi \rangle_{X^{\star}, X}=0,
\end{equation*}
with $(u_\tau,\dot{u}_\tau)\vert_{t=0}=(0, 0)$. Under Assumption~\ref{assumption_wellposedness} on well-posedness of the problem on the reference domain, the transported problem has a unique solution.  \\
\indent The \index{Eulerian derivative} Eulerian derivative of $J$ at $\Omega_l$ in the direction of the vector field $\Theta$ is defined as
$$dJ(u,\Omega_l)\Theta=\displaystyle \lim_{\tau \rightarrow 0} \frac{1}{\tau}(J(u_\tau,\Omega_{l,\tau})-J(u,\Omega_l)).$$
\indent  It is said that the functional $J$ is \emph{shape differentiable} at $\Omega_l$ if $dJ(u,\Omega_l)\Theta$ exists for all $\Theta \in C^{1,1}(\overline{\Omega},\mathbb{R}^d)$ and defines a continuous linear functional on $C^{1,1}(\overline{\Omega},\mathbb{R}^d)$.\\
\indent Since the shape derivative defines a linear and bounded functional, according to the Riesz representation theorem it has a unique representative in a Hilbert space. This representative is commonly referred to as the shape gradient, i.e. the gradient is defined as the solution of
\begin{equation*} 
 (\nabla J, \phi)_H=dJ(u, \Omega_l)\phi, \ \forall \phi \in H.
\end{equation*}
This allows a certain freedom in choosing the Hilbert space $H$ and in turn the descent direction. A comparison of $L^2$ and weighted $H^1$  scalar products defined on the boundary can be found in~\cite{Dogan}. In the present work we will consider the common case $H=L^2(\partial \Omega_l)$. \\
\indent For sufficiently small $|\tau|$, we then bring back the transported solution to the domain with the fixed lens by introducing 
\begin{align*}
u^\tau(x,t) =u_\tau(F_\tau(x),t),
\end{align*}
which will satisfy
\begin{equation} \label{weak_form_transported_solution} 
\begin{aligned} 
& \int_0^T \int_{\Omega} \Bigl\{(1-2ku^\tau)\ddot{u}^\tau  \phi+c^2A_\tau\nabla u^\tau \cdot A_\tau\nabla \phi \\
&+ b A_\tau\nabla \dot{u}^\tau \cdot A_\tau\nabla \phi 
-2k(\dot{u}^\tau)^2\phi \Bigr\} \, I_\tau \, dx \, ds \\
&+\displaystyle \int_0^T \int_{\Gamma_a} (c \dot{u}^\tau+\frac{b}{c}\ddot{u}^\tau)\phi \, dx \, ds- \int_0^T  \int_{\Gamma_n} (c^2g+b \dot{g})\phi  \, dx \, ds =0,
\end{aligned}
\end{equation}
for all $\phi \in X$, where we have employed the following notation
\begin{equation*}\label{defIAw}
I_\tau=\text{det}(DF_\tau)\,, \quad A_\tau=(DF_\tau)^{-T},
\end{equation*}
with $DF_\tau$ being the Jacobian of the transformation $F_\tau$. Note that we have made use of the fact that $F_\tau$ vanishes on $\Gamma_n \cup \Gamma_a$ to end up with the boundary terms in \eqref{weak_form_transported_solution}. \\
\indent In what is to follow we will rely on some useful properties of $F_\tau$, which we for the convenience of the reader recall here. 
\begin{lemma} \label{lemma_mappingF}\cite{IKP08, Zolesio} There exists $\tau_0>0$ such that the mapping $F_\tau$, defined by \eqref{perturbation of identity mapping}, has the following properties for $\tau \in (-\tau_0, \tau_0)$:
\begin{equation*}
\begin{aligned}
&\tau \mapsto F_\tau\in C^1(-\tau_0,\tau_0;C^1(\overline{\Omega},\mathbb{R}^d)),
&&\tau\mapsto A_\tau\in C(-\tau_0,\tau_0;C(\overline{\Omega},\mathbb{R}^{d\times d})),\\	
&\tau\mapsto F_\tau^{-1}\in C(-\tau_0,\tau_0;C^1(\overline{\Omega},\mathbb{R}^d)),
&&\tau\mapsto I_\tau\in C^1(-\tau_0,\tau_0;C(\overline{\Omega})).\\
\end{aligned}
\end{equation*}
Furthermore, it holds
\begin{equation*}
\begin{aligned}
&\frac{d}{d\tau} F_\tau\vert_{\tau=0}=\Theta,		&&\frac{d}{d\tau} F_\tau^{-1}\vert_{\tau=0}=-\Theta,\\
&\frac{d}{d\tau} DF_\tau\vert_{\tau=0}=D\Theta,	&&\frac{d}{d\tau} DF_\tau^{-1}\vert_{\tau=0}=\frac{d}{d\tau} (A_\tau)^T\vert_{\tau=0}=-D\Theta,\\
&\frac{d}{d\tau} I_\tau\vert_{\tau=0}=\textnormal{div}\, \Theta, && \frac{d}{d\tau} A_\tau\vert_{\tau=0} =-(D\Theta)^T.
\end{aligned}
\end{equation*}
\end{lemma}
\noindent From Lemma~\ref{lemma_mappingF}, it follows that there exist $\alpha_0, \alpha_1, \alpha_2 >0$, such that 
\begin{equation}  \label{uniform_bound_I}
\begin{aligned}
& 0<\alpha_0 \leq I_\tau(x) \leq \alpha_1,  \quad  \text{for} \  x \in \overline{\Omega}, \ \tau \in
 [-\tau_0,\tau_0], \\
&  |A_\tau|_{L^\infty(\Omega)} \leq \alpha_2, \quad \ \tau \in
 [-\tau_0,\tau_0].
\end{aligned}
\end{equation}
Note that we have used the fact that $I_\tau=\text{det}(DF_\tau)>0$ on $\Omega$ for sufficiently small $|\tau|$ when applying the transformation of integrals formula to obtain~\eqref{weak_form_transported_solution}.
\subsection{An auxiliary result on H\" older continuity} We will compute the shape derivative of our pressure-tracking shape functional by the so-called variational approach, introduced in~\cite{IKP08}. This method avoids the condition of shape differentiability of the state variable by replacing it with H\" older continuity with respect to domain perturbations.\\
\indent Although the derivation of the volume representation of the shape derivative is analogous to the one presented in~\cite{NikolicKaltenbacher}, we include it here for the sake of completeness. We remark however, that the computation of the boundary representation of the shape derivative differs from~\cite{NikolicKaltenbacher} due to different boundary conditions (influencing piece-wise regularity of the solution) and a lack of nonlinear damping in the present model, which allows for simplifications of the final derivative expression.\\ 
\indent By employing the transformation of integrals formula it can be shown that $u^\tau$ is uniformly bounded with respect to $\tau$, i.e.
\begin{equation} \label{bound_transported_solution}
\begin{aligned}
&\|\nabla \dot{u}^\tau\|_{L^{2}(0,T;L^{2}(\Omega))} \leq \, \frac{\tau_0\,|D\Theta|_{L^\infty(\Omega)}}{\sqrt{\alpha_0}}\,\|\nabla \dot{u}_\tau\|_{L^{2}(0,T; L^{2}(\Omega))} \leq  \frac{\tau_0|D\Theta|_{L^\infty(\Omega)}}{\sqrt{\alpha_0}}\,\bar{m}, \\
&\|\dot{u}^\tau\|_{L^{\infty}(0,T;L^{2}(\Omega))} \leq \, \frac{1}{\sqrt{\alpha_0}}\,\|\dot{u}_\tau\|_{L^{\infty}(0,T;L^{2}(\Omega))} \leq \frac{1}{\sqrt{\alpha_0}} \,\bar{m}, \\
&\|\ddot{u}^\tau\|_{L^{2}(0,T;L^{2}(\Omega))} \leq \, \frac{1}{\sqrt{\alpha_0}}\,\|\ddot{u}_\tau\|_{L^{2}(0,T;L^{2}(\Omega))}\leq \frac{1}{\sqrt{\alpha_0}}\, \bar{m},
\end{aligned}
\end{equation} 
for all $\tau \in (-\tau_0, \tau_0)$, with $\alpha_0$ being the lower bound of $I_\tau$ as given in \eqref{uniform_bound_I}. Similarly we have
\begin{align*}
\|u^{\tau}\|_{L^\infty(0,T; L^\infty(\Omega))} \leq \|u_\tau\|_{L^\infty(0,T; L^\infty(\Omega))} \leq \bar{a} < \frac{1}{2|k|_{L^\infty(\Omega)}}.
\end{align*}
\indent Let us define 
\begin{align*}
\bar{\mathcal{X}}=\{& \,z \in W^{1, \infty}(0,T; L^2(\Omega))\, \cap \, H^1(0,T; H^1(\Omega)): \\
 & \,  z_{\vert t=0}=\dot{z}_{\vert t=0}=0, \ \textup{tr}_{\Gamma_a} \dot{z} \in  L^\infty(0,T; L^2(\Gamma_a)) \, \}.
\end{align*}
\begin{theorem} \label{theorem_continuity} Under Assumption \ref{assumption_wellposedness}, the solution of problem~\eqref{state_constraint_weak} is H\" older continuous with exponent $1/2$ with respect to domain perturbations; i.e. it holds
\begin{align}
\displaystyle \lim_{\tau \rightarrow 0} \frac{\|u^\tau-u\|^2_{\bar{\mathcal{X}}}}{\tau}=0.
\end{align}
\end{theorem}
\begin{proof}
We begin by noting that the difference $z^\tau=u^\tau-u$ satisfies
\begin{equation} \label{weak_form_difference}
\begin{aligned}
& \int_0^T \int_\Omega \{(1-2ku^\tau)\,\ddot{z}^\tau \phi +c^2 \nabla z^\tau \cdot \nabla \phi+b  \nabla \dot{z}^\tau \cdot \nabla \phi \\ & \hspace{5mm}-2k(\dot{u}+\dot{u}^\tau)\dot{z}^\tau \phi-2k z^\tau \ddot{u} \phi\} \, dx ds +\displaystyle \int_0^T \int_{\Gamma_a} (c \dot{z}^\tau+\frac{b}{c}\ddot{z}^\tau)\phi \, dx \, ds\\
 =& \, \langle \, f_\tau(u^\tau, \Omega_l), \phi \, \rangle_{X^\star, X},
\end{aligned}
\end{equation}
for all $\phi \in X=L^2(0,T; H^1(\Omega))$, where the right hand side is given by
\begin{equation*} 
\begin{aligned}
&\langle \, f_\tau(u^\tau, \Omega_l), \phi \, \rangle_{X^\star, X} \\
=& \int_0^T \int_\Omega \Bigl \{-(I_\tau-1)(1-2ku^\tau)\,\ddot{u}^\tau\phi +2k\,(I_\tau-1)\,(\dot{u}^\tau)^2 \phi \vspace{0.2mm}\\
&-c^2 \,((I_\tau-1)A_\tau \nabla u^\tau \cdot A_\tau \nabla \phi +(A_\tau-I)\nabla u^\tau \cdot A_\tau \nabla \phi + \nabla u^\tau \cdot (A_\tau-I)\nabla \phi) \vspace{1mm} \\
&-b\,((I_\tau-1)A_\tau \nabla \dot{u}^\tau \cdot A_\tau \nabla \phi +(A_\tau-I)\nabla \dot{u}^\tau \cdot A_\tau \nabla \phi + \nabla \dot{u}^\tau \cdot (A_\tau-I)\nabla \phi) \Bigr \}\, dx ds.
\end{aligned}
\end{equation*}
Testing~\eqref{weak_form_difference} by $\phi=\dot{z}^\tau\chi_{[0,t]} \in X$ then leads to
\begin{align*}
& \frac{1}{2}\int_\Omega (1-2ku^\tau(t))\,(\dot{z}^\tau(t))^2 \, dx+\frac{1}{2}\int_\Omega c^2 \,|\nabla z^\tau(t)|^2\, dx+\frac{1}{2} \int_0^t\int_\Omega b \,|\nabla \dot{z}^\tau|^2\, dx ds \\
&+ \int_0^T \int_{\Gamma_a} c \, |\dot{z}^\tau|^2 \, dx ds+\frac{1}{2}\int_{\Gamma_a} \frac{b}{c} \, |\dot{z}^\tau|^2 \, dx  \\
=& \, \displaystyle \int_0^t \int_\Omega 2k\,(z^\tau \dot{z}^\tau\ddot{u}+(\dot{u}+\dot{u}^\tau)\,(\dot{z}^\tau)^2) \, dx ds+\langle \,f_\tau(u^\tau, \Omega_l), \dot{z}^\tau\chi_{[0,t]} \,\rangle_{\bar{\mathcal{X}}^\star, \bar{\mathcal{X}}} \\
\leq& \,  2|k|_{L^\infty(\Omega)}\,(C^{\Omega}_{H^1,L^4})^2\,\Bigl[\,
\|\ddot{u}\|_{L^{2}(0,T;L^2(\Omega))}\,\|z^\tau\|_{L^\infty(0,T;H^1(\Omega))}\,\|\dot{z}^\tau\|_{L^2(0,T;H^1(\Omega))}\\
&+(\|\dot{u}^\tau\|_{L^{\infty}(0,T;L^{2}(\Omega))}+\|\dot{u}\|_{L^{\infty}(0,T;L^{2}(\Omega))})\,\|\dot{z}^\tau\|^2_{L^{2}(0,T;H^1(\Omega))}\,\Bigr] \\
&+\langle \, f_\tau(u^\tau, \Omega_l), \dot{z}^\tau\chi_{[0,t]} \, \rangle_{\bar{\mathcal{X}}^\star, \bar{\mathcal{X}}}.
\end{align*}
In the last estimate, we have made use of the continuous embedding $H^1(\Omega) \hookrightarrow  L^4(\Omega)$, cf.~\cite[Corollary 4.53]{Demengel}, with the embedding constant $C^\Omega_{H^1,L^4}$. For further estimating the terms in the last line we employ
\begin{align*}
\|z^\tau\|_{L^{\infty}(0,T;H^1(\Omega))} \leq& \, \sqrt{T}\,\|\dot{z}^\tau\|_{L^{2}(0,T;H^1(\Omega))}, 
\end{align*} (since $z^\tau \vert_{t=0}=\dot{z}^\tau \vert_{t=0}=0$), to conclude that for sufficiently small $\bar{m}$ and final time $T$ it holds
\begin{equation} \label{convergence1_}
\begin{aligned}
&\|z^\tau\|^2_{\bar{\mathcal{X}}} \leq C\,|\,\langle \, f_\tau(u^\tau, \Omega_l), \dot{z}^\tau \chi_{[0,t]}\, \rangle_{\bar{\mathcal{X}}^\star, \bar{\mathcal{X}}}\,|,
\end{aligned}
\end{equation}
for some sufficiently large $C>0$ that does not depend on $\tau$. We can estimate the terms on the right hand side in~\eqref{convergence1_} as follows
\begin{equation} \label{estimate_rhs}
\begin{aligned} 
&|\,\langle\,  f_\tau(u^\tau, \Omega_l), \dot{z}^\tau \chi_{[0,t]}\, \rangle_{\bar{\mathcal{X}}^\star, \bar{\mathcal{X}}}\,|   \\
 \leq& \, |I_\tau-1|_{L^\infty(\Omega)}\,\Bigl ((1+2|k|_{L^\infty(\Omega)}\,\bar{a}\,)\|\ddot{u}^\tau\|_{L^2(0,T;L^2(\Omega))}\|\dot{z}^\tau\|_{L^2(0,T;L^2(\Omega))} \\
 &+2|k|_{L^\infty(\Omega)}(C^{\Omega}_{H^1,L^4})^2\|\dot{u}^\tau\|^2_{L^2(0,T;H^1(\Omega))}\|\dot{z}^\tau\|_{L^{\infty}(0,T;L^{2}(\Omega))}\Bigr ) \\
&+(|I_\tau-1|_{L^\infty(\Omega)}\,\alpha_2^2 +|A_\tau-I|_{L^\infty(\Omega)}(1+\alpha_2)) 
\Bigl(|c|_{L^\infty(\Omega)}^2\|\nabla u^\tau\|_{L^2(0,T;L^2(\Omega))} \\
&+|b|_{L^\infty(\Omega)}\|\nabla \dot{u}^\tau\|_{L^2(0,T;L^2(\Omega))}\Bigr)\,\|\nabla \dot{z}^\tau\|_{L^2(0,T;L^2(\Omega))},
\end{aligned}
\end{equation}
where $\bar{a}$ is defined as in~\eqref{space_U} and $\alpha_2$ as in~\eqref{uniform_bound_I}. Recalling that $u_0=0$, we can further derive
\begin{align*}
& \|\dot{z}^\tau\|_{L^2(0,T;L^2(\Omega))} \leq T\,\|\dot{z}^\tau\|_{L^\infty(0,T;L^2(\Omega))}, \\
& \|\nabla u^\tau\|_{L^2(0,T;L^2(\Omega))} \leq T\sqrt{T}\,\|\nabla \dot{u}^\tau\|_{L^2(0,T;L^2(\Omega))},
\end{align*}
which, together with uniform boundedness of $u^\tau$ in the sense of~\eqref{bound_transported_solution}, leads to
\begin{equation} \label{estimate_rhs_improved}
\begin{aligned} 
&|\,\langle \, f_\tau(u^\tau, \Omega_l), \dot{z}^\tau \chi_{[0,t]} \, \rangle_{\bar{\mathcal{X}}^\star, \bar{\mathcal{X}}}\,|   \\
 \leq& \, C(\bar{a}, \bar{m}, T)\,(|I_\tau-1|_{L^\infty(\Omega)} +|A_\tau-I|_{L^\infty(\Omega)}) \,\|z^\tau\|_{\bar{\mathcal{X}}}.
\end{aligned}
\end{equation}
 By combining~\eqref{estimate_rhs_improved} with estimate~\eqref{convergence1_} and then employing the  properties of the mapping $F_\tau$ from Lemma~\ref{lemma_mappingF}, we first conclude that
\begin{equation} \label{convergence2}
\begin{aligned}
&\displaystyle \lim_{\tau \rightarrow 0}\|z^\tau\|_{\bar{\mathcal{X}}}=0.
\end{aligned}
\end{equation} Secondly, we divide~\eqref{convergence1_} by $\tau \neq 0$, and then it remains to show that $$\displaystyle \lim_{\tau \rightarrow 0} \frac{|\,\langle \, f_\tau(u^\tau, \Omega_l), \dot{z}^\tau\chi_{[0,t]}\, \rangle_{\bar{\mathcal{X}}^\star, \bar{\mathcal{X}}}\,|}{\tau}=0.$$ This now follows from estimate~\eqref{estimate_rhs_improved}, Lemma~\ref{lemma_mappingF} and~\eqref{convergence2}.
\end{proof}
 
\section{Shape derivative of the cost function} \label{Section_ComputationShapeDer}
Next, we will make use of the tools of shape calculus and the continuity result of the previous section to compute the shape derivative. We will begin by calculating the derivative of $J$ in terms of the volume integral over the whole domain $\Omega$, which is sometimes regarded as the weak shape derivative.
\begin{proposition} \label{thm_shape_derivative_volume}
Let Assumption~\ref{assumption_wellposedness} on the well-posedness of the state problem hold and let $u_d \in L^2(0, T; L^2(\Omega))\, \cap \, H^1(0,T; (H^1(\Omega))^\star)$. The volume representation of the shape derivative of $J$ at $\Omega_l$ in the direction of $\Theta \in C^{1,1}(\overline{\Omega}, \mathbb{R}^d)$, where $\textup{supp}\,\Theta \cap (\overline{D} \cup \partial \Omega)=\emptyset$, is given by 
\begin{equation} \label{formula_weak_shape_derivative}
\begin{aligned}
& dJ(u,\Omega_l)\Theta \\ 
=& \, \int_0^T \int_{\Omega}\Bigl\{ 
(c^2 \nabla u^T+b \nabla \dot{u}^T)(D\Theta^T \nabla p+D\Theta \nabla p) \Bigr\} \, dx \, ds  \\
& - \int_0^T \int_{\Omega} \Bigl\{(1-2ku )\ddot{u} p  + c^2\nabla u \cdot \nabla p +b\nabla \dot{u} \cdot \nabla p 
-2k(\dot{u})^2p \Bigr\} \,\textnormal{div} \, \Theta \, dx \, ds. \\
\end{aligned}
\end{equation}
\end{proposition}
\begin{proof}
By recalling the definition of the shape derivative and the fact that $\textup{supp}~\Theta \cap \overline{D}=\emptyset$, together with the continuity result of Theorem \ref{theorem_continuity}, we arrive at
\begin{align*}
d J(u, \Omega_l)\Theta 
=& \,  \lim_{\tau \rightarrow 0} \frac{1}{\tau} \int_0^T \int_\Omega (u+u^\tau-2u_d)(u^{\tau}-u)\, \chi_D \, dx ds \\
=& \, \lim_{\tau \rightarrow 0} \frac{1}{\tau} \int_0^T \int_\Omega 2(u-u_d)(u^{\tau}-u)\, \chi_D \, dx ds. 
\end{align*}
Note that the last integral appears as the right hand side in the adjoint problem~\eqref{adjoint_problem_weak} when the test function is $u^\tau-u$. We employ again the notation $z^\tau=u^\tau-u$, recalling that $z^\tau_{\vert t=0}=\dot{z}^\tau_{\vert t=0}=0$ and $p_{\vert t=T}=\dot{p}_{\vert t=T}=0$. Testing the adjoint problem~\eqref{adjoint_problem_weak} with $z^\tau \in L^2(0,T; H^1(\Omega))$ and integrating by parts with respect to time leads to
\begin{align*}
& \int_0^T \int_\Omega 2(u-u_d)\chi_D \,z^{\tau}\, dx ds \\
=& \, \displaystyle \int_0^T \int_{\Omega} \{(1-2ku)\,\ddot{p}\, z^{\tau}+c^2\nabla p \cdot \nabla z^{\tau}-b\nabla \dot{p} \cdot \nabla z^{\tau} \} \, dx \, ds\\
&+\displaystyle \int_0^T \int_{\Gamma_a} (-c \dot{p}+\frac{b}{c}\ddot{p})z^\tau\, dx ds \\
=& \,\displaystyle \int_0^T \int_{\Omega} \{(1-2ku)\,\ddot{z}^{\tau}p+c^2 \nabla z^{\tau} \cdot \nabla p +b \nabla \dot{z}^{\tau} \cdot \nabla p -4k\dot{u}\dot{z}^{\tau}p-2k\ddot{u}z^{\tau}p\} \, dx \, ds \\
&+\displaystyle \int_0^T \int_{\Gamma_a} (c \dot{z}^\tau+\frac{b}{c}\ddot{z}^\tau)p\, dx ds \\
=&\, \langle \, f_\tau(u^\tau, \Omega_l), p \, \rangle_{X^\star, X}+\int_0^T \int_\Omega 2k\,((\dot{z}^{\tau})^2+z^{\tau}\ddot{z}^\tau p) \, dx ds,
\end{align*}
where in the last line we have then also used the weak form~\eqref{weak_form_difference} satisfied by $z^\tau$. Due to the continuity result of Theorem~\ref{theorem_continuity}
\begin{align*}
\displaystyle \lim_{\tau \rightarrow 0} \frac{1}{\tau} \int_0^T \int_\Omega 2k\,((\dot{z}^{\tau})^2+z^{\tau}\ddot{z}^\tau p) \, dx ds 
= \, \displaystyle \lim_{\tau \rightarrow 0} \frac{1}{\tau} \int_0^T \int_\Omega -2k \,z^{\tau} \dot{z}^{\tau} \dot{p} \, dx ds 
=\, 0,
\end{align*} 
and our computations reduce to
\begin{align*}
&d J(u,\Omega_l)\Theta= \, \displaystyle \lim_{\tau \rightarrow 0} \frac{1}{\tau}\langle \,f_\tau(u^\tau, \Omega_l), p \,\rangle_{X^\star, X}. 
\end{align*}
Thanks to Theorem~\ref{theorem_continuity} and Lemma~\ref{lemma_mappingF}, it is straightforward to verify that
\begin{align}\label{rhs_auxillary_rearranged}
\displaystyle \lim_{\tau \rightarrow 0} \frac{1}{\tau} \langle \,f_\tau(u^\tau, \Omega_l), p \,\rangle_{X^\star, X}=\displaystyle \lim_{\tau \rightarrow 0} \frac{1}{\tau}\langle \,f_\tau(u, \Omega_l), p \,\rangle_{X^\star, X}.
\end{align}
By employing Lemma~\ref{lemma_mappingF} again, we can directly compute the limit on the right hand side of~\eqref{rhs_auxillary_rearranged}, which leads to the expression \eqref{thm_shape_derivative_volume}.
\end{proof}
\subsubsection{Imposing stronger regularity assumptions} According to the Hadamard-Zol\' esio Structure Theorem~\cite[Theorem 3.5]{Zolesio}, when $J(u, \Omega_l)$ and the domain are smooth enough, the shape derivative can be represented in the form
\begin{align*}
dJ(u, \Omega_l)\Theta=\displaystyle \int_ {\Gamma} \mathcal{G}(x)\, \Theta \cdot n_l \, dx.
 \end{align*}
To be able to express the derivative in terms of an integral over the boundary of the lens we have to impose stronger assumptions with respect to regularity of $u$, $p$, $u_d$ and domains $\Omega_l$, $\Omega_f$. From this point on let $\Omega_l$, $\Omega_f$ be of class $C^{1,1}$. Moreover, we introduce the following assumption:
\begin{assumption}\label{stronger_reg_assumptions} 
 We assume that the state and the adjoint variable are sufficiently regular in space on each of the subdomains, i.e.
\begin{equation} 
\begin{aligned}
& u \in \{ u \in \mathcal{U} \ \vert \ u_{\vert \Omega_{l, f}} \in H^1(0,T; \mathcal{M}(\Omega_{l,f}))\},\\
& p \in \{ p \in \mathcal{P} \ \vert \ p_{\vert \Omega_{l,f}} \in H^1(0,T; \mathcal{M}(\Omega_{l,f}))\},
\end{aligned}
\end{equation}
where $\mathcal{M}(\Omega)$ is the linear space induced by the norm
\begin{align*}
\displaystyle \|w\|_\mathcal{M}=\|w\|_{H^1(\Omega)}+\|\Delta w\|_{L^2(\Omega)}+\displaystyle \Bigl \|\frac{\partial w}{\partial n} \Bigr \|_{H^{\tilde{s}}(\partial \Omega)},
\end{align*}
for some $\tilde{s} \in (0, \frac{1}{2})$. 
\end{assumption}
\noindent Let us recall a useful identity that generalizes integration by parts in $\mathcal{M}(\Omega)$, which we will employ to transform expression \eqref{formula_weak_shape_derivative}.
\begin{lemma}\cite[Lemma 6]{KaltenbacherPeichl}
Let $\Omega \subset \mathbb{R}^3$ be a domain of class $C^{1,1}$. For $u, v \in \mathcal{M}(\Omega)$ and $\Theta \in C^{1,1}(\overline{\Omega}, \mathbb{R}^d)$, the following identity holds
\begin{equation} \label{identity}
\begin{aligned}
&\int_ \Omega (I \textup{div}\, \Theta-D\Theta-D\Theta^T)\nabla u\cdot \nabla v \, dx \\
=& \, \int_{\Omega} (\Delta u (\Theta \cdot \nabla v)) + \Delta v (\Theta \cdot \nabla u)) \, dx-\int_{\partial \Omega} (\frac{\partial u}{\partial n}(\Theta \cdot v)+\frac{\partial v}{\partial n}(\Theta \cdot u))\, dx \\
&+\int_{\partial \Omega} \nabla u \cdot \nabla v (\Theta \cdot n) \, dx,
\end{aligned}
\end{equation}
where $n$ denotes the outer unit normal to $\partial \Omega$.
\end{lemma}
\noindent We are now ready to state the main result.
\begin{theorem}
Let $\partial \Omega_l$ and $\partial \Omega_f$ be $C^{1,1}$ regular. Let Assumption~\ref{assumption_wellposedness} on the well-posedness of the state problem hold as well as Assumption~\ref{stronger_reg_assumptions} on regularity of solutions to the state and the adjoint problems. The shape derivative of $J$ at $\Omega_l$ in the direction of the deformation field $\Theta \in C^{1,1}(\overline{\Omega},\mathbb{R}^d)$, where $\textup{supp}\,\Theta \cap (\overline{D} \cup \partial \Omega_f\setminus \partial \Omega_l)=\emptyset$, is then given by
\begin{equation} \label{shape_derivative_boundary}
\begin{aligned}
& dJ(u,\Omega_l)\Theta \\
=& \, -\displaystyle \int_0^T \int_{\Gamma} \Bigl ( 2\llbracket \, k \,\rrbracket \, u\dot{u} \dot{p} 
+ \llbracket \,c^2 \, \rrbracket  \, \nabla u_l \cdot \nabla p_f +\llbracket \, b \,\rrbracket \, \nabla \dot{u}_l \cdot \nabla p_f\Bigr )\, \Theta^Tn_l \, dx \, ds,
\end{aligned} 
\end{equation}
where $\Gamma$ denotes the interface between $\Omega_l$ and $\Omega_f$.
\end{theorem}
\begin{proof}
 Under Assumption \ref{stronger_reg_assumptions}, we can rewrite the expression~\eqref{formula_weak_shape_derivative} for the weak shape derivative as the sum of integrals over $\Omega_l$ and $\Omega_f$ and directly apply formula~\eqref{identity} to the terms containing $\nabla u \cdot \nabla p$ and $\nabla \dot{u} \cdot \nabla p$. The remaining terms containing $\textnormal{div} \, \Theta$ can be dealt with by employing integration by parts with respect to time and space
\begin{equation} \label{simplify_divTheta}
\begin{aligned}
&- \int_0^T \int_{\Omega} ((1-2ku )\ddot{u} p-2k(\dot{u})^2p ) \,\textnormal{div} \, \Theta \, dx \, ds \\
=& \,\int_0^T \int_{\Omega} (1-2ku )\dot{u} \dot{p} \,\textnormal{div} \, \Theta \, dx \, ds \\
=&\, \int_0^T \int_{\Gamma} \llbracket\, (1-2ku )\dot{u} \dot{p} \,\rrbracket \, \Theta \cdot n_l \, dx \, ds-\int_0^T \int_{\Omega} \nabla((1-2ku )\dot{u} \dot{p} ) \cdot \Theta \, dx \, ds.
\end{aligned}
\end{equation}
Here we have also made use of the fact that $(u, \dot{u})_{\vert t=0}=(p, \dot{p})_{\vert t=T}=(0,0)$. The last term in \eqref{simplify_divTheta} can be further rewritten as 
\begin{align*}
&-\int_0^T \int_{\Omega} \nabla((1-2ku)\dot{u} \dot{p} ) \cdot \Theta \, dx \, ds\\
=&\,\int_0^T \int_{\Omega}\Bigl\{((1-2ku )\ddot{u}-2k\dot{u}^2)\,\nabla p \cdot \Theta+(1-2ku)\ddot{p}\,  \nabla u \cdot \Theta \Bigr \}\, dx \, ds.
\end{align*}
In this way we transform~\eqref{formula_weak_shape_derivative} into
\begin{equation*}
\begin{aligned}
& dJ(u,\Omega_l)\Theta  \\
=&\, \displaystyle \sum_{i \in \{l,f\}}  \int_0^T \int_{ \Omega_i} ((1-2k_iu_i)\ddot{u}_i-c_i^2\Delta u_i 
-b_i\Delta \dot{u}_i-2k_i\dot{u}_i^2 )\, \nabla p_i \cdot \Theta \, dx \, ds  \\
& +\displaystyle \sum_{i \in \{l,f\}}  \int_0^T \int_{ \Omega_i} ((1-2k_iu_i)\ddot{p}_i-c_i^2\Delta p_i 
+b_i\Delta \dot{p}_i)\, \nabla u_i \cdot \Theta \, dx \, ds  \\
& +\int_0^T \int_{\Gamma} \Bigl \llbracket\, (1-2ku)\dot{u}\dot{p}-c^2\nabla u \cdot \nabla p-b\nabla \dot{u} \cdot \nabla p \,\Bigr \rrbracket \, \Theta \cdot n_l \, dx \, ds \\
&+\int_0^T \int_{\Gamma} \Bigl \llbracket \,(c^2\nabla u+b\nabla \dot{u})\cdot n_l (\nabla p \cdot \Theta) \,\Bigr \rrbracket \, dx \, ds \\
&+\int_0^T \int_{\Gamma} \Bigl \llbracket \,\Bigl(c^2\nabla p-b\nabla \dot{p} \Bigr)\cdot n_l (\nabla u \cdot \Theta)\,\Bigr \rrbracket \, dx \, ds.\\
\end{aligned}
\end{equation*}
The first two sums on the right-hand side vanish and we are left with
\begin{equation} \label{strong_form_der_1}
\begin{aligned}
dJ(u,\Omega_l)\Theta
=& \,\int_0^T \int_{\Gamma} \Bigl \llbracket\, (1-2ku)\dot{u}\dot{p}-c^2\nabla u \cdot \nabla p-b\nabla \dot{u} \cdot \nabla p \,\Bigr \rrbracket \, \Theta \cdot n_l \, dx \, ds \\
&+\int_0^T \int_{\Gamma} \Bigl \llbracket\, (c^2\nabla u+b\nabla \dot{u})\cdot n_l (\nabla p \cdot \Theta) \,\Bigr \rrbracket \, dx \, ds \\
&+\int_0^T \int_{\Gamma} \Bigl \llbracket\, \Bigl(c^2\nabla p-b\nabla \dot{p} \Bigr)\cdot n_l (\nabla u \cdot \Theta)\,\Bigr \rrbracket \, dx \, ds. \\
\end{aligned}
\end{equation}
It is possible to slightly simplify expression~\eqref{strong_form_der_1}, mimicking the steps taken in \cite[Section 3.2]{IKP}. If we denote by $\displaystyle \frac{\partial}{\partial \nu}$ the tangential derivative, then due to the interface conditions on $\Gamma$ it holds
\begin{equation} \label{interface_conditions}
\begin{aligned}
&\displaystyle \Bigl \llbracket\, \frac{\partial u}{\partial \nu}\,\Bigr  \rrbracket=\Bigl \llbracket\, \frac{\partial p}{\partial \nu} \,\Bigr \rrbracket=0,\\
&\displaystyle \Bigl \llbracket\, c^2 \frac{\partial u}{\partial n_l}+b \frac{\partial \dot{u}}{\partial n_l} \,\Bigr \rrbracket=\Bigl \llbracket \,c^2 \frac{\partial p}{\partial n_l}-b \frac{\partial \dot{p}}{\partial n_l} \,\Bigr \rrbracket=0.
\end{aligned}
\end{equation}
Hence we have the following identity 
 \begin{align*}
\Bigl \llbracket \,(c^2\nabla u+b\nabla \dot{u})\cdot n_l (\nabla p \cdot \Theta)\,\Bigr \rrbracket = \,\Bigl \llbracket \,c^2\frac{\partial  u}{\partial n_l}\frac{\partial  p}{\partial n_l}+b\frac{\partial  \dot{u}}{\partial n_l}\frac{\partial  p}{\partial n_l}\,\Bigr \rrbracket\,(\Theta \cdot n_l),
\end{align*}
as well as
\begin{align*}
&\Bigl \llbracket \, (c^2\nabla p-b\nabla \dot{p})\cdot n_l (\nabla u \cdot \Theta)\,\Bigr \rrbracket
=\, \Bigl \llbracket \,c^2\frac{\partial  u}{\partial n_l}\frac{\partial p}{\partial n_l}-b\frac{\partial  u}{\partial n_l}\frac{\partial \dot{p}}{\partial n_l}\,\Bigr \rrbracket \,(\Theta \cdot n_l). 
\end{align*}
By plugging these identities into~\eqref{strong_form_der_1} and making use of $\llbracket \, u \dot{u} \dot{p} \, \rrbracket =0$, the derivative reduces to
\begin{equation*} 
\begin{aligned}
& dJ(u,\Omega_l)\Theta  \\
=& \, -\int_0^T \int_{\Gamma} \Bigl( 2 \, \llbracket \, k \, \rrbracket u\dot{u}\dot{p}+\llbracket \, c^2\frac{\partial  u}{\partial \nu}\frac{\partial  p}{\partial \nu}+b\frac{\partial  \dot{u}}{\partial \nu}\frac{\partial  p}{\partial \nu}\\
&\hspace{12mm} -c^2\frac{\partial  u}{\partial n_l}\frac{\partial p}{\partial n_l}-b\frac{\partial  \dot{u}}{\partial n_l}\frac{\partial p}{\partial n_l}\,\Bigr \rrbracket\, \Bigr ) \Theta^Tn_l \, dx \, ds.
\end{aligned} 
\end{equation*}
Furthermore, due to conditions~\eqref{interface_conditions} on the interface and again the fact that $(u, \dot{u})_{\vert t=0}=(0,0)$ and $(p, \dot{p})_{\vert t=T}=(0,0)$, it holds
\begin{align*}
&\int_0^T \int_{\Gamma}\Bigl(\llbracket\, c^2 \,\rrbracket \frac{\partial u}{\partial \nu} \frac{\partial p}{\partial \nu}+\llbracket \,b \,\rrbracket \frac{\partial \dot{u}}{\partial \nu} \frac{\partial p}{\partial \nu}-\Bigl \llbracket \,c^2\frac{\partial  u}{\partial n_l}\frac{\partial p}{\partial n_l} 
+b\frac{\partial  \dot{u}}{\partial n_l}\frac{\partial p}{\partial n_l} \,\Bigr \rrbracket \Bigr )\, \Theta^Tn_l \, dx ds\\
=&\, \frac{1}{2} \int_0^T \int_{\Gamma} \llbracket\,c^2 \,\rrbracket (\nabla u_l \cdot \nabla p_f+\nabla u_f \cdot \nabla p_l)+\llbracket\, b \,\rrbracket (\nabla \dot{u}_l \cdot \nabla p_f+\nabla \dot{u}_f \cdot \nabla p_l) \, \Theta^Tn_l \, dx ds,
\end{align*}
which helps us further to obtain
\begin{equation*} 
\begin{aligned}
& dJ(u,\Omega_l)\Theta  \\
=& \,- \int_0^T \int_{\Gamma} \Bigl ( 2\llbracket\, k \,\rrbracket u\dot{u} \dot{p} 
+\frac{1}{2}\llbracket\, c^2\, \rrbracket (\nabla u_l \cdot \nabla p_f+\nabla u_f \cdot \nabla p_l)) \\
&+\frac{1}{2}\llbracket\, b \,\rrbracket (\nabla \dot{u}_l \cdot \nabla p_f+\nabla \dot{u}_l \cdot \nabla p_f) \Bigr )\, \Theta^Tn_l \, dx \, ds.
\end{aligned} 
\end{equation*}
Lastly, we can make use of the fact that again due to interface conditions~\eqref{strong_form_der_1} it holds
\begin{align*}
&\int_0^T \int_{\Gamma}(\llbracket \, c^2 \,\rrbracket  \nabla u_l\cdot \nabla p_f +\llbracket \,b \,\rrbracket \nabla \dot{u}_l \cdot \nabla p_f )\, \Theta^Tn_l \,  dx ds \\
 =& \,  \int_0^T \int_{\Gamma}(\llbracket \, c^2 \,\rrbracket  \nabla u_f \cdot \nabla p_l +\llbracket\, b \,\rrbracket \nabla \dot{u}_f \cdot \nabla p_l )\, \Theta^Tn_l \,  dx ds,
\end{align*}
which finally leads to \eqref{shape_derivative_boundary}.
\end{proof}
\section{Isogeometric analysis} \label{Section_IGA} In this section, we give a brief overview of basic concepts of isogeometric analysis which will be used to solve our state and adjoint problems and set the notation for future use. For an in-depth study of the isogeometric approach, we refer the reader to classical literature \cite{hughes:09, Piegl, Shumaker}. 
\subsection{Basics of B-splines and NURBS} \subsubsection{Univariate B-splines} We begin by introducing a knot vector in one dimension
\begin{align*}
\Xi=\{0=\xi_1, \ldots, \xi_{n+l+1}=1\}
\end{align*}
where $l$ is the polynomial order and $n$ is the number of basis functions used to construct the B-spline curve. As knot values are allowed to repeat, let us also introduce $Z=\{\zeta_1, \ldots, \zeta_m\}$ as a knot vector without any repetitions. $Z$ partitions the parameter space $\hat{\Omega}=[0,1]$ into elements. \\
\indent We will denote by $B_{i, q}$, $i=1,\ldots,d$, B-Splines defined on $\Xi$ recursively
\begin{align*}
& \hat{B}_{i, 0}(\hat{x})=\begin{cases}\, 1 \quad \text{ if } \xi_i \leq \hat{x} < \xi_{i+1}, \\ \, 0 \quad \text{ otherwise},\end{cases} \vspace{2mm} \\
&\hat{B}_{i, q}(\hat{x})=\frac{\hat{x}-\xi_i}{\xi_{i+q}-\hat{x}}\hat{B}_{i, q-1}(\hat{x})+\frac{\xi_{i+q+1}- \hat{x}}{\xi_{i+q+1}-\xi_{i+1}}\hat{B}_{i+1, q-1}(\hat{x}), \quad l=1,2,3,\ldots 
\end{align*}
 where in case one of the denominators is zero the contribution is assumed to be zero and omitted. The spline space can then be introduced as
\begin{align*}
\mathcal{S}^q(\Xi)=\text{span}\{\, \hat{B}_{i, q}(\hat{x}), \ i=1, \ldots,d \, \}.
\end{align*}
\subsubsection{Multivariate B-splines} Multivariate splines are defined as a tensor product of univariate B-splines. For every parametric direction $s=1, \ldots, d$ we introduce the number $n_s$ of univariate $B$-spline functions, the open knot vector $\Xi_s$ and the knot vector without repetitions $Z_s$. To simplify exposition we assume that the degree $l$ is the same in all directions. Then we can define multivariate knot vectors as follows:
\begin{align*}
&\mathbf{\Xi}=\Xi_1 \times \Xi_2 \ldots \times \Xi_d, \\
&\mathbf{Z}=Z_1 \times Z_2 \ldots \times Z_d.
\end{align*}
$\mathbf{Z}$ forms a mesh of the parametric domain $\hat{\Omega}=[0,1]^d$. \\
\indent We define multivariate $B$-spline functions for each multi-index $i \in \mathcal{I}=\{i=(i_1,\ldots,i_d)\}$ as follows
\begin{align*}
\hat{B}_{i, q}(\bm{\hat{x}})=\hat{B}_{i_1, q}(\hat{x}_1)\cdot\ldots \cdot \hat{B}_{i_d, q}(\hat{x}_d),
\end{align*}
and the corresponding spline space is given by
\begin{align*}
\mathcal{S}^l(\mathbf{\Xi})=\text{span}\{\, \hat{B}_{i, q}, \ i \in \mathcal{I} \, \}.
\end{align*}
\subsubsection{NURBS} Going one step further, NURBS are then introduced as rational functions of B-spline functions. This creates a possibility to exactly represent various geometries including conic sections, which makes them especially popular in computer aided design (CAD). For a set of positive weights $\{w_i, \ i \in \mathcal{I}\}$, we define the NURBS functions as
\begin{align*}
\hat{N}_{i, q}(\bm{\hat{x}})=\displaystyle \frac{w_i \hat{B}_{i, q}(\bm{\hat{x}})}{\displaystyle \sum_{r \in \mathcal{I}} w_r \hat{B}_{r, q}(\bm{\hat{x}})}.
\end{align*}
Note that $B$-splines are obtained as special cases of NURBS when all weights are equal.
\subsubsection{NURBS geometries}
For a given set of control points $\mathcal{C}_i \in \mathbb{R}^d$, $i \in \mathcal{I}$ we can define a NURBS surface ($d=2$) or a NURBS solid ($d=3$) by taking a linear combination of NURBS basis functions:
\begin{align*}
G(\bm{\hat{x}})=\sum_{i \in \mathcal{I}} \, \mathcal{C}_i \hat{N}_{i, q}(\bm{\hat{x}}).
\end{align*}
The NURBS geometry is then defined as $\Omega=G(\hat{\Omega})$. We note that in isogeometric analysis \emph{one map}  $G: \hat{\Omega} \rightarrow \Omega$ takes the patch from the parameter space to the physical space. 
However, from an implementational point of view, it still makes sense to speak about elements in which the domain gets divided through $\mathbf{Z}$, since, for instance, quadrature nodes can then be defined in a uniform way for all the elements.
\subsection{Description of the computational lens/fluid domain in the multipatch setting}
\begin{figure}[h] 
\centering
\begin{tikzpicture}[scale=1.7,>=stealth,label distance=.2cm]
		
		\colorlet{Blau}{blue!50!cyan!75!white};
		\colorlet{Gruen}{green!80!black!75!white};
		\colorlet{Oranje}{orange!75!white};
		\colorlet{ltblue}{blue!20!cyan!75!white};
		\colorlet{ltgray}{black!30!white};
		
		\coordinate (A1) at (0,-0.33);
		\coordinate (B1) at (1.37,-0.33);
		\coordinate (C1) at (2.25,-0.33);
		\coordinate (A2) at (0,1);
		\coordinate (B2) at (1.37,1);
		\coordinate (C2) at (2.25,1);
		\coordinate (A3) at (0,2);
		\coordinate (B3) at (1.37,2);
		\coordinate (C3) at (2.25,2);
		\coordinate (A4) at (0,3);
		\coordinate (B4) at (1.37,3);
		\coordinate (C4) at (2.25,3);
		\coordinate (L) at (0,0.37);


		\filldraw[Blau] (B2)--(B1)--(A1)--(L) arc (-90:-30:1.75)--cycle;
		\draw[black,thick] (B2)--(B1)--(A1)--(L) arc (-90:-30:1.75)--cycle;
		
		\filldraw[Oranje] (L) arc (-90:-30:1.75) -- (B2) -- (A2) -- cycle;
		\draw[black,very thick] (L) arc (-90:-30:1.75) -- (B2) -- (A2) -- cycle;
		
		\filldraw[Blau] (B1)--(C1)--(C2)--(B2)--cycle;
		\draw[black] (B1)--(C1)--(C2)--(B2)--cycle;
		
		\filldraw[Blau] (A2)--(B2)--(B3)--(A3)--cycle;
		\draw[black] (A2)--(B2)--(B3)--(A3)--cycle;
		
		\filldraw[Blau] (B2)--(C2)--(C3)--(B3)--cycle;
		\draw[black] (B2)--(C2)--(C3)--(B3)--cycle;
		
		\filldraw[ltblue] (A3)--(B3)--(B4)--(A4)--cycle;
		\draw[black] (A3)--(B3)--(B4)--(A4)--cycle;
		
		\filldraw[Blau] (B3)--(C3)--(C4)--(B4)--cycle;
		\draw[black] (B3)--(C3)--(C4)--(B4)--cycle;
		
		\filldraw[light-gray] (B2)--(B1)--(A1)--(L) arc (-90:-30:1.75)--cycle;
		\draw[black,thick] (B2)--(B1)--(A1)--(L) arc (-90:-30:1.75)--cycle;
		
		\filldraw[white] (L) arc (-90:-30:1.75) -- (B2) -- (A2) -- cycle;
		\draw[black,very thick] (L) arc (-90:-30:1.75) -- (B2) -- (A2) -- cycle;
		
		\filldraw[light-gray] (B1)--(C1)--(C2)--(B2)--cycle;
		\draw[black] (B1)--(C1)--(C2)--(B2)--cycle;
		
		\filldraw[light-gray] (A2)--(B2)--(B3)--(A3)--cycle;
		\draw[black] (A2)--(B2)--(B3)--(A3)--cycle;
		
		\filldraw[light-gray] (B2)--(C2)--(C3)--(B3)--cycle;
		\draw[black] (B2)--(C2)--(C3)--(B3)--cycle;
		
		\filldraw[light-gray] (A3)--(B3)--(B4)--(A4)--cycle;
		\draw[black] (A3)--(B3)--(B4)--(A4)--cycle;
		
		\filldraw[light-gray] (B3)--(C3)--(C4)--(B4)--cycle;
		\draw[black] (B3)--(C3)--(C4)--(B4)--cycle;
		
		\node[shape=circle,draw,inner sep=2pt] at (0.625,0.1) {\tiny 1};
		\node[shape=circle,draw,inner sep=2pt] at (1.8,0.1) {\tiny  2};
		\node[shape=circle,draw,inner sep=2pt] at (0.625,0.75) {\tiny 3};
		\node[shape=circle,draw,inner sep=2pt] at (0.625,1.5) {\tiny 4};
		\node[shape=circle,draw,inner sep=2pt] at (1.8,1.5) {\tiny 5};
		\node[shape=circle,draw,inner sep=2pt] at (0.625,2.5) {\tiny 6};
		\node[shape=circle,draw,inner sep=2pt] at (1.8,2.5) {\tiny 7};

		\draw[<->] (0,-0.53)--(1.37,-0.53);
		\node[shape=circle,inner sep=2pt] at (0.69,-0.67) {$W$};
		\draw[dashed,ltgray] (0,-0.33)--(0,-0.78);
		\draw[dashed,ltgray] (1.37,-0.33)--(1.37,-0.53);
		
		\draw[<->] (0,-0.78)--(2.25,-0.78);
		\node[shape=circle,inner sep=2pt] at (1.125,-0.93) {$B$};
		\draw[dashed,ltgray] (2.25,-0.33)--(2.25,-0.78);

		\draw[<->] (2.9,-0.33)--(2.9,2);
		\node[shape=circle,inner sep=2pt] at (3.1,0.83) {$S$};
		\draw[dashed,ltgray] (2.25,-0.33)--(2.9,-0.33);
		\draw[dashed,ltgray] (2.25,2)--(2.9,2);
		
		\draw[<->] (2.65,-0.33)--(2.65,1);
		\node[shape=circle,inner sep=2pt] at (2.775,0.33) {$K$};
		\draw[dashed,ltgray] (2.25,1)--(2.65,1);
		
		\draw[<->] (-0.2,0.37)--(-0.2,1);
		\node[shape=circle,inner sep=2pt] at (-0.35,0.67) {$P$};
		\draw[dashed,ltgray] (-0.2,0.37)--(0,0.37);
		\draw[dashed,ltgray] (-0.2,1)--(0,1);

		\draw[<->] (-0.2,-0.33)--(-0.2,0.37);
		\node[shape=circle,inner sep=2pt] at (-0.35,0.05) {$R$};
		
		\draw[<->] (-0.45,-0.33)--(-0.45,3);
		\node[shape=circle,inner sep=2pt] at (-0.6,1.39) {$L$};
		\draw[dashed,ltgray] (-0.45,-0.33)--(0,-0.33);
		\draw[dashed,ltgray] (-0.45,3)--(0,3);

		\draw[->,very thick] (0,-0.33)--(2.5,-0.33);
		\node[shape=circle,inner sep=2pt] at (2.55,-0.43) {$\normalsize x$};
		\draw[->,very thick] (0,-0.33)--(0,3.25);
		\node[shape=circle,inner sep=2pt] at (-0.1,3.25) {$\normalsize y$};
		\filldraw[black] (0,-0.33) circle (0.015cm);
\end{tikzpicture}

\caption{Patches used for the discretization. The lens domain is given by patch 3. \label{fig:Patches}}
\end{figure}
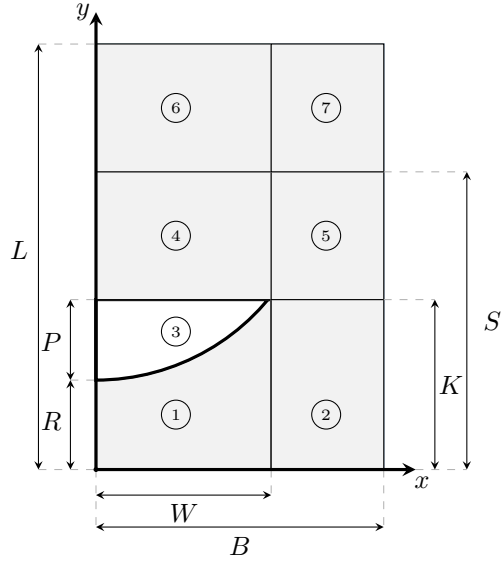

In our numerical experiments, we will assume that the problem is symmetric with respect to the $y$-axis through the center of the lens. Hence it suffices to only discretize half of the domain $\Omega$ and use symmetry boundary conditions at the $y$ axis. \\
\indent In practical examples computational domains are often described with several patches. In our case the full (half-)domain is decomposed into $N_{\mathrm{ptc}}=7$ patches, $\Omega^{(j)}$, $j=1, \ldots, N_{\mathrm{ptc}}$, six of which represent water and one the lens region (see Figure~\ref{fig:Patches}). It is assumed that $\overline{\Omega}=\bigcup_{j=1}^{N_{\mathrm{ptc}}} \, \overline{\Omega^{(j)}}$ and $\Omega^{(j)}$ are assumed to be disjoint.\\
\indent Every patch has one map $G^{(j)}$, $j \in \{1, \ldots, N_{\mathrm{ptc}}\}$, which maps the parametric domain to the physical one.  Since the pressure $u$ is continuous across the interface, but its normal derivative is not, we assume that the meshes and the control points coincide on the interface and impose strong $C^0$- but not $C^1$-continuity between patches. Furthermore, this assumption holds in every optimization step, i.e. for every updated geometry. 
\subsection{Numerical treatment of the state problem}
We will employ a semi-discrete method, where we discretize the space using a Galerkin finite element scheme and formulate the problem for continuous time. The main principle of isogeometric analysis is that the NURBS basis which we used to exactly model our domain also serves as the basis for the solution space of our numerical method.\\
\indent  Following \cite{BuffaVazquez}, in each patch $\Omega^{(j)}$ we define a discrete space
\begin{align*}
V_h^{(j)}=\text{span} \{\, N^{(j)}_{i, q} =\hat{N}^{(j)}_{i,q} \circ G^{(j)^{-1}}, \ i \in \mathcal{I}^{(j)} \, \}, \ \ j=1, \ldots, N_{\mathrm{ptc}}.
\end{align*}
The discrete domain on the whole space will consist of continuous functions which, restricted on each domain, belong to $V_h^{(j)}$. In other words, we seek an approximate solution $u^h$ such that
\begin{align*}
u^h(\cdot, t) \in V_h=\{\, v \in C(\Omega) \ | \ v_{\vert \Omega^{(j)}} \in V_h^{(j)}, \, j=1, \ldots, N_{\mathrm{ptc}} \, \}, \ t \in [0,T].
\end{align*}
\begin{remark}
One has to distinguish between a patch-wise and a global set of degrees of freedom. Within an individual patch the numerical solution has the form $$u^{h, (j)}(x,t)=\sum_{i\in \mathcal{I}^{(j)}}N_{i,q}^{(j)}(x)\,u^i(t).$$ If two patches are ``glued together'', due to our assumptions on coinciding control points and meshes, it is possible to identify degrees of freedom with each other and also to add their respective ansatz functions. \\
\indent For example, when using two rectangular patches of $3\times 3$ degrees of freedom each and aligning them on one common boundary we would have to identify the three degrees of freedom of that boundary of patch 2 with the respective degrees of freedom of patch 1, giving a total of $15$ instead of $9+9=18$ global degrees of freedom. This approach ensures global continuity and hanging nodes are avoided. In this global set of ansatz functions one introduces a new (global) numbering of degrees of freedom and calls the set of all global ansatz functions (after identification on all interfaces) $\mathcal{I}$ such that the global numerical solution can be written in the following way:
\begin{equation*}
	u^h(x,t)=\sum_{i\in \mathcal{I}}N_{i,q}(x)\, u^i(t),
\end{equation*}
where the $N_{i,q}$ are the global ansatz functions (after identification).
\end{remark}

Let us denote coefficients of a function $u$ in the global NURBS basis by  $$\underline{u}=[{u}^1, \ {u}^2, \ \ldots ,\ {u}^{\textup{ndof}} ]^T,$$ 
where time dependency is omitted for notational convenience. After discretizing both state and domain with NURBS, we arrive at the following semidiscrete formulation of our state problem: 
\begin{align} \label{state_semidiscrete}
 \begin{cases}
\ \mathbf{M} \ddot{\underline{u}}+\mathbf{C} \dot{\underline{u}}+\mathbf{K} \underline{u}-\bm{\mathcal{T}}(\dot{\underline{u}})\dot{\underline{u}}-\bm{\mathcal{T}}(\underline{u})\ddot{\underline{u}}=\mathbf{F}-\mathbf{{A}}_1\dot{\underline{u}}-\mathbf{{A}}_2\ddot{\underline{u}}, \vspace{1mm} \\
\ \underline{u}(0)=0, \vspace{1mm} \\
\ \dot{\underline{u}}(0)=0,
\end{cases}
 \end{align}
 where mass, damping, and stiffness matrices are defined by 
 
 \begin{enumerate}
 	\item[1.)] Computing
				 	\begin{equation*}
				 	\mathbf{M}' =\left(\begin{array}{lcccccr}
				 	\mathbf{M}^{(1)} &  &  &  &  &  &  \\
				 	& \mathbf{M}^{(2)} &  &  &  &  &  \\
				 	& & \mathbf{M}^{(3)} &  &  & \scalebox{1.5}{$0$} &  \\
				 	&  &  & \mathbf{M}^{(4)} &  &  &  \\
				 	&  \scalebox{1.5}{$0$} &  &  & \mathbf{M}^{(5)} &  &  \\
				 	 &  &  &  &  & \mathbf{M}^{(6)} &  \\
				 	 &   &  &  &  && \mathbf{M}^{(7)} \\
				 	\end{array}\right),
				 	\end{equation*}
			 	with $\mathbf{M}^{(j)}$, $j=1,...,N_{\mathrm{ptc}}$ being the patch-wise mass matrices, defined by
			 	\begin{equation*}
			 		\mathbf{M}^{(j)} = \left[m^{(j)}_{i_1i_2}\right], \ m^{(j)}_{i_1i_2}=\int_{\hat{\Omega}}\hat{N}_{i_1}\hat{N}_{i_2} \, |\text{det}DG^{(j)}| \,  d\hat{x}, 
			 	\end{equation*}
			 	where $DG^{(j)}$ denotes the Jacobian, $(DG^{(j)})^{-T}$ the inverse transposed of the patch-wise mappings $G^{(j)}$ and we have omitted polynomial degree $q$ of basis functions to simplify notation. \\
			 	
	\item[2.)] Secondly, identifying degrees of freedom at interfaces by adding the respective rows and columns in the global matrix $\mathbf{M}'$. For instance, at the interface between patch 1 and 2 we add the rows corresponding to interface degrees of freedom of patch 2 to the respective ones of patch 1 and then remove them from the global matrix. Proceeding analogously with the columns and repeating the process for all interfaces, we obtain $\mathbf{M}$ from $\mathbf{M'}$. This corresponds to transforming $\mathbf{M}'$ via $\mathbf{M} = \mathbf{E^TM}'\mathbf{E}$, where $\mathbf{E}$ is a rectangular matrix that only consists of ones and zeros.
 \end{enumerate}
 
We proceed in the same manner to obtain the damping matrix $\mathbf{C}$ and stiffness matrix $\mathbf{K}$, with the only difference being their patch-wise definitions compared to the mass matrix, i.e.:
	 	\begin{align*}
		 	\mathbf{C}^{(j)} &= \left[c^{(j)}_{i_1i_2}\right], \ c^{(j)}_{i_1i_2}=\int_{\hat{\Omega}} b(DG^{(j)})^{-T}\nabla \hat{N}_{i_1} \cdot (DG^{(j)})^{-T}\nabla\hat{N}_{i_2} \, |\text{det}DG^{(j)}| \, d\hat{x}, \\
		 	\mathbf{K}^{(j)} &= \left[k^{(j)}_{i_1i_2}\right], \ k^{(j)}_{i_1i_2}=\int_{\hat{\Omega}} c^2(DG^{(j)})^{-T}\nabla \hat{N}_{i_1} \cdot (DG^{(j)})^{-T}\nabla \hat{N}_{i_2} \, |\text{det}DG^{(j)}| \, d\hat{x}. 
	 	\end{align*}

\begin{remark}
 	\noindent Even though the isogeometric approach does not rely on individual elements and hence mass and stiffness matrices can analytically be calculated as above via integrals over the whole (parametric) domain $\hat{\Omega}$, it is still more convenient to adopt the notion of elements from classical FEA to the support boundaries of the individual NURBS basis functions, that are given by $\mathbf{Z}$. This speeds up the computation as the integrals are only evaluated over parts of the domain $\hat{\Omega}$, where the respective basis functions do not vanish. In this sense one can also reformulate the definition of $\mathbf{M}^{(j)},\mathbf{C}^{(j)}$ and $\mathbf{K}^{(j)}$ in terms of the assembly operator known from classical finite elements with $n^{(j)}_e$ being the number of such elements per patch:
 
 \begin{align*}
& \mathbf{M}^{(j)}= \bigwedge\limits_{e=1}^{n^{(j)}_e} [m^{(j)}_{i_1i_2}], \ m^{(j)}_{i_1i_2}= \int_{\hat{\Omega}^e} N_{i_1}N_{i_2} \, |\text{det}DG^{(j)}| \,  d\hat{x}, \\
& \mathbf{C}^{(j)}= \bigwedge\limits_{e=1}^{n^{(j)}_e} [c^{(j)}_{i_1i_2}], \ c^{(j)}_{i_1i_2}=\int_{\hat{\Omega}^e} b(DG^{(j)})^{-T}\nabla \hat{N}_{i_1} \cdot (DG^{(j)})^{-T}\nabla \hat{N}_{i_2} \, |\text{det}DG^{(j)}| \, d\hat{x}, \\
& \mathbf{K}^{(j)}= \bigwedge\limits_{e=1}^{n^{(j)}_e} [k^{(j)}_{i_1i_2}], \ \displaystyle k^{(j)}_{i_1i_2}=\int_{\hat{\Omega}^e} c^2(DG^{(j)})^{-T}\nabla \hat{N}_{i_1} \cdot (DG^{(j)})^{-T}\nabla \hat{N}_{i_2} \, |\text{det}DG^{(j)}| \, d\hat{x}. 
 \end{align*}

In the same way also the tensor $\bm{\mathcal{T}}$ is defined as the global structure after interface identification of degrees of freedom of the patch-wise tensors
\begin{align*}
 & \bm{\mathcal{T}}^{(j)}(\underline{u})= \bigwedge\limits_{e=1}^{n^{(j)}_e} [\mathcal{T}^{(j)}_{i_1i_2}(\underline{u})], \quad \mathcal{T}^{(j)}_{i_1i_2}(\underline{u})=\int_{\hat{\Omega}^e} 2k \hat{N}_{i_1}\hat{N}_{i_2} \left(\displaystyle \sum_{i_3=1}^{\textup{ndof}^{(j)}}\hat{N}_{i_3} u^{i_3}\right) \, |\text{det}DG^{(j)}| \, d \hat{x}, 
 \end{align*}
and the right-hand side vector and the matrices corresponding to the absorbing conditions are patch-wise given as
\begin{align*}
 & \mathbf{F}^{(j)}=\bigwedge\limits_{e=1}^{n^{(j)}_e} [f^{(j)}_{i}], \quad f^{(j)}_{i}=\int_{\hat{\Gamma}_n^e} \hat{N}_i (c_f^2g \circ G^{(j)}+b_f\dot{g}\circ G^{(j)})\, |(DG^{(j)})^{-T} \hat{n}|\,|\text{det}DG^{(j)}| \,d \hat{x}, \\
 & \mathbf{A}^{(j)}_{1}=\bigwedge\limits_{e=1}^{n^{(j)}_e} [a^{(j)}_{1, i_1i_2}], \quad a^{(j)}_{1, i_1i_2}=\int_{\hat{\Gamma}_a^e} c_f \hat{N}_{i_1} \hat{N}_{i_2}\,|(DG^{(j)})^{-T} \hat{n}|\,|\text{det}DG^{(j)}| \, d \hat{x}, \\
 & \mathbf{A}^{(j)}_{2}=\bigwedge\limits_{e=1}^{n^{(j)}_e} [a^{(j)}_{2, i_1i_2}], \quad a^{(j)}_{2, i_1i_2}=\int_{\hat{\Gamma}_a^e} \frac{c_f}{b_f} \hat{N}_{i_1} \hat{N}_{i_2}\, |(DG^{(j)})^{-T} \hat{n}|\,|\text{det}DG^{(j)}| \, d\hat{x},
 \end{align*}
 which are then - as shown for $\mathbf{M}$ before - put together into global structures $\mathbf{F},\mathbf{A_1}$ and $\mathbf{A_2}$ by interface identification of degrees of freedom.
 \end{remark}
\subsection{Time stepping scheme} Due to the nonlinear nature of our problem, the higher harmonics of the fundamental frequency of the initial wave will grow with the distance from the source. In numerical algorithms cutting off higher harmonics is known to cause spurious oscillations at the shock peaks; this is often regarded as the Gibbs phenomenon.   \\
\indent We tackle this issue by employing the generalized $\alpha$-scheme for time stepping, introduced first by Chung and Hulbert in \cite{ChungHulbert}. The generalized $\alpha$-method allows to control the level of numerical dissipation of high frequencies with minimal unwanted damping in the low frequency range. In this way the high frequency modes which are insufficiently resolved by the spatial or the temporal discretization are eliminated. \\
\indent Let $0=t_0<t_1<\ldots t_m=T$ be a given uniform partition of the time domain and let us define the time step as $\Delta t=t_{n+1}-t_n$. The $\alpha$-method applied to the semi-discrete equation \eqref{state_semidiscrete} yields the modified equation
\begin{equation} \label{GenAlpha}
\begin{aligned}
&\mathbf{M} \ddot{\underline{u}}_{n+1-\alpha_m}+\mathbf{C} \dot{\underline{u}}_{n+1-\alpha_f}+\mathbf{K} \underline{u}_{n+1-\alpha_f} \\
&-\bm{\mathcal{T}}(\dot{\underline{u}}_{n+1-\alpha_f})\dot{\underline{u}}_{n+1-\alpha_f}-\bm{\mathcal{T}}(\underline{u}_{n+1-\alpha_f})\ddot{\underline{u}}_{n+1-\alpha_f}\\
=& \,\mathbf{F}_{n+1-\alpha_f}-\mathbf{{A}}_1\dot{\underline{u}}_{n+1-\alpha_f}-\mathbf{{A}}_2\ddot{\underline{u}}_{n+1-\alpha_m},
\end{aligned}
\end{equation}
where
\begin{align*}
& \ddot{\underline{u}}_{n+1-\alpha_m}=(1-\alpha_m)\ddot{\underline{u}}_{n+1}+\alpha_m \ddot{u}_n, \\
& \dot{\underline{u}}_{n+1-\alpha_f}=(1-\alpha_f)\dot{\underline{u}}_{n+1}+\alpha_f \dot{u}_n,\\
& \underline{u}_{n+1-\alpha_f}=(1-\alpha_f)\underline{u}_{n+1}+\alpha_f u_n,
\end{align*}
and
\begin{align*}
\mathbf{F}_{n+1-\alpha_f}=\mathbf{F}((1-\alpha_f)t_{n+1}+\alpha_ft_n).
\end{align*}
Note that we have adopted a mid-point approach for the nonlinearities in \eqref{GenAlpha}. Morover, the classical Newmark relations \cite{Newmark} are used
\begin{align*}
\underline{u}_{n+1}&=\underline{u}_{n}+\Delta t \, \underline{\dot{u}}_{n+1}+\frac{1}{2}\Delta t^2((1-2\beta)\underline{\ddot{u}}_{n}+2\beta\underline{\ddot{u}}_{n+1}), \\
\underline{\dot{u}}_{n+1}&=\underline{\dot{u}}_{n}+\Delta t((1-\gamma)\underline{\ddot{u}}_{n}+\gamma\underline{\ddot{u}}_{n+1}),
\end{align*}
which will be numerically realized as a predictor-corrector algorithm. When $\alpha_m=\alpha_f=0$, the time stepping reduces to the standard Newmark scheme. Relations among parameters $\alpha_m$, $\alpha_f$, $\beta$, and $\gamma$ that guarantee stability and desired dissipation levels are given in \cite{Erlicher} for systems with nonlinear internal forces involving $u$, but not its time derivatives. Adopting the strategy from~\cite{Hoffelner, manfred}, we define predictors for $\underline{u}_{n+1}$ and $\underline{\dot{u}}_{n+1}$ as follows
\begin{align*}
&\underline{u}_{\textup{pred}}=\underline{u}_{n}+\Delta t \dot{\underline{u}}_n+\frac{\Delta t^2}{2}(1-2\beta)\ddot{\underline{u}}_n \\
&\dot{\underline{u}}_{\textup{pred}}=\dot{\underline{u}}_n+(1-\gamma)\Delta t \ddot{\underline{u}}_n,
\end{align*} 
and the effective mass matrix $$\mathbf{\overline{M}}=(1-\alpha_m)\mathbf{M}+\gamma(1-\alpha_f) \Delta t \mathbf{C}+ \beta (1-\alpha_f)\Delta t^2 \mathbf{K},$$
which leads to the following system
\begin{equation} \label{GenAlpha_system}
\begin{aligned}
&\hspace{-12mm}\mathbf{\overline{M}} \, \ddot{\underline{u}}_{n+1} 
+ ((1-\alpha_f)\underline{u}_{\textup{pred}}+\alpha_f\underline{u}_n)+\mathbf{C} ((1-\alpha_f)\dot{\underline{u}}_{\textup{pred}}+\alpha_f\dot{\underline{u}}_n)\\
=&\,\mathbf{F}_{n+1-\alpha_f}+\bm{\mathcal{T}}(\underline{\dot{u}}_{n+1-\alpha_f})\underline{\dot{u}}_{n+1-\alpha_f}+\bm{\mathcal{T}}(\underline{u}_{n+1-\alpha_f})\underline{\ddot{u}}_{n+1-\alpha_m} \\
& -\mathbf{{A}}_1\dot{\underline{u}}_{n+1-\alpha_f} -\mathbf{{A}}_2\ddot{\underline{u}}_{n+1-\alpha_m}.
\end{aligned}
\end{equation} 
\noindent The system~\eqref{GenAlpha_system} is nonlinear, so we have to solve it by an iterative scheme. Following~\cite{manfred}, we will employ a fixed-point iteration with respect to the second time derivative. Note that the effective mass matrix remains unchanged during the time stepping. An LU decomposition can therefore be performed once in the first time step, while in subsequent steps matrix inversion is avoided. \\
\vspace*{0cm}

\begin{algorithm}[H]
\caption{Solving the state problem}\label{algorithm_state} \vspace{2mm}
$t=0$, $n=0$\vspace{2mm}\\
$\underline{u}_n=0$, $\dot{\underline{u}}_n=0$, $\mathbf{M} \, \ddot{\underline{u}}_n=\mathbf{F}_n-\mathbf{C} \dot{\underline{u}}_n-\mathbf{K} \underline{u}_n$\vspace{2mm}\\
\text{compute the effective mass matrix:} \vspace{2mm} \\
$\mathbf{\overline{M}}=(1-\alpha_m)\mathbf{M}+\gamma(1-\alpha_f) \Delta t \mathbf{C}+ \beta (1-\alpha_f)\Delta t^2 \mathbf{K}$  \vspace{2mm} \\
\SetKwBlock{whileloop}{while $t\leq T$ do}{}
\whileloop{
\text{perform predictor step:} \vspace{2mm} \\
$\underline{u}_{\textup{pred}}=\underline{u}_{n}+\Delta t \dot{\underline{u}}_n+\frac{\Delta t^2}{2}(1-2\beta)\ddot{\underline{u}}_n$\vspace{2mm}\\
$\dot{\underline{u}}_{\textup{pred}}=\dot{\underline{u}}_n+(1-\gamma)\Delta t \ddot{\underline{u}}_n$\vspace{1mm}\\
$\kappa=0$\vspace{2mm}\\
$\underline{u}_{n+1}^{\kappa}=\underline{u}_{\textup{pred}}$\vspace{2mm}\\
$\dot{\underline{u}}_{n+1}^{\kappa}=\underline{\dot{u}}_{\textup{pred}}$\vspace{2mm}\\
\SetKwBlock{whileloopi}{while $\displaystyle \frac{\|\ddot{\underline{u}}^{\kappa+1}_{n+1}-\ddot{\underline{u}}^\kappa_{n+1}\|}{\|\ddot{\underline{u}}^{\kappa+1}_{n+1}\|} > \textup{TOL}_{u}$ do}{}
\whileloopi{
 \text{solve the algebraic system of equations:} \vspace{2mm} 
 \begin{align*}
 	&\hspace{-12mm}\mathbf{\overline{M}} \, \ddot{\underline{u}}^{\kappa+1}_{n+1} 
 	= \mathbf{F}_{n+1-\alpha_f}-\mathbf{K} ((1-\alpha_f)\underline{u}_{\textup{pred}}+\alpha_f\underline{u}_n)-\mathbf{C} ((1-\alpha_f)\dot{\underline{u}}_{\textup{pred}}+\alpha_f\dot{\underline{u}}_n)\\
 	&\hspace*{0.5cm}+\bm{\mathcal{T}}(\underline{\dot{u}}_{n+1-\alpha_f}^\kappa)\underline{\dot{u}}_{n+1-\alpha_f}^\kappa+\bm{\mathcal{T}}(\underline{u}_{n+1-\alpha_f}^\kappa)\underline{\ddot{u}}_{n+1-\alpha_m}^\kappa \\
 	&\hspace*{0.5cm} -\mathbf{{A}}_1\dot{\underline{u}}_{n+1-\alpha_f}^\kappa -\mathbf{{A}}_2\ddot{\underline{u}}_{n+1-\alpha_m}^\kappa 
 \end{align*} \\}}
\end{algorithm}
\newpage

\begin{algorithm}[H]
\nonl \setcounter{AlgoLine}{13}
\SetKwBlock{whileloop}{}{end}
\nonl \whileloop{
	\SetKwBlock{whileloopi}{}{end}
	\nonl \whileloopi{
\hspace{-0mm}\text{perform corrector step:} \vspace{2mm}  \\
$\underline{u}_{n+1}^{\kappa+1}=\underline{u}_{\textup{pred}}+\beta \Delta t^2 \ddot{\underline{u}}^{\kappa+1}_{n+1}$\vspace{2mm}\\
$\dot{\underline{u}}^{\kappa+1}_{n+1}=\dot{\underline{u}}_{\textup{pred}}+\gamma \Delta t \ddot{\underline{u}}^{\kappa+1}_{n+1}$\vspace{2mm}\\
$\kappa \mapsto \kappa+1$\\
}  
$n \mapsto n+1$, $t \mapsto t+\Delta t$\\
}
\end{algorithm} 

\vspace{5mm}

\subsection{Numerical treatment of the adjoint problem} For solving the adjoint problem, we will use the same NURBS space and the same time discretization as for the state solver. Let us denote by $$\underline{u}_d=[u_{d}^1 \ u_{d}^2 \ \ldots u_{d}^{\textup{ndof}}\ ]^T$$ the coefficients of the target pressure $u_d$ in the NURBS base and by $\underline{u}_{d,n}$ the coefficients at time step $n$. After discretizing also the adjoint state with NURBS, we obtain the semidiscrete formulation of our adjoint problem
\begin{align} \label{adjoint_semidiscretized}
\begin{cases}
\ (\mathbf{M}-\bm{\mathcal{T}}(u)) \ddot{\underline{p}}-\mathbf{C} \dot{\underline{p}}+\mathbf{K} \underline{p}=2\mathbf{M}^{D}(\underline{u}-\underline{u}_d)+\mathbf{{A}}_1\dot{\underline{p}}-\mathbf{{A}}_2\ddot{\underline{p}}, \vspace{1mm}\\
\ \underline{p}(T)=0, \vspace{1mm} \\
\ \dot{\underline{p}}(T)=0.
\end{cases}
\end{align}
\indent The matrices $\mathbf{M}$, $\mathbf{C}$, $\mathbf{K}$ and tensor $\bm{\mathcal{T}}$ in \eqref{adjoint_semidiscretized} are defined in the same way as for the state problem. Note that the integral appearing on the right hand side in the weak formulation \eqref{adjoint_problem_weak} of the adjoint problem only needs to be numerically evaluated on the focal region $D$. Therefore, in order to define the matrix $\mathbf{M}^{D}$, we first introduce
\begin{equation*}
\mathbf{M}^{D,(j)} = \left[m^{D,(j)}_{i_1i_2}\right], \quad \quad m^{D,(j)}_{i_1i_2}=\int_{\hat{\Omega}}\hat{N}_{i_1}\hat{N}_{i_2} \, |\text{det}DG^{(j)}| \,(\chi_D\circ G^{(j)})\ d\hat{x},
\end{equation*}
 where $j=1,...,N_{\textup{ptc}}$. Compared to $\mathbf{M}^{(j)}$, now we only integrate over the part of $D$ where the pressure is tracked. This is expressed by the use of the characteristic function $\chi_D$ taking the value $1$ within $D$ and $0$ elsewhere. We define
\begin{equation*}
{\mathbf{M}^D}' =\left(\begin{array}{lccccc}
~~ \hspace*{-0mm}\mathbf{M}^{D,(1)} &  &  &  &  & \\
& \hspace*{-0mm}\mathbf{M}^{D,(2)} &  &  &  & \\
&  & \hspace*{-0mm}\mathbf{M}^{D,(3)} &  &  & \scalebox{1.5}{$0$} \\
&  &  & \hspace*{-0mm}\mathbf{M}^{D,(4)} &  & \\
& \hspace*{-5mm}\scalebox{1.5}{$0$} &  &  & \hspace*{-0mm}\mathbf{M}^{D,(5)} &  \\
 &  &  &  &  & \hspace*{-14mm}\mathbf{M}^{D,(6)} \\
 &   &  &  &  &\hspace*{12mm}\mathbf{M}^{D,(7)} \\
\end{array}\right).
\end{equation*}	
The matrix $\mathbf{M}^D$ is then obtained from ${\mathbf{M}^D}'$ by again performing interface identification of degrees of freedom as it has been already done for the global mass, damping, and stiffness matrices.\\
\indent Although the adjoint problem is linear, the modified mass matrix $\mathbf{M}-\bm{\mathcal{T}}(u)$ is time-dependent and thus would have to be inverted in every time step of our numerical solver.  To avoid this computationally expensive step, the adjoint problem is also solved iteratively. The time stepping is performed \emph{backward in time} by a Newmark predictor-corrector scheme to obtain $\underline{p}$ and thus $p^h$ as the solution. The Newmark parameters $\beta_p$ and $\gamma_p$ for the adjoint problem are chosen in such a way that second order of accuracy is achieved. 

\vspace*{0.5cm}
\begin{algorithm}[H]
\caption{Solving the adjoint problem}\label{algorithm_adjoint} \vspace{2mm}
$t=T$, \ $n=\text{number of time steps}$\vspace{2mm}\\
$\underline{p}_n=0$, $\dot{\underline{p}}_n=0$, $\ddot{\underline{p}}_n=0$\vspace{2mm}\\
\text{compute the effective mass matrix:} \vspace{2mm} \\
$\mathbf{\overline{M}} =\mathbf{M}+\gamma_p(-\Delta t) \mathbf{C}+ \beta_p \Delta t^2 \mathbf{K}$ \vspace{2mm}\\
\While{$t \geq 0$ \vspace{2mm}}{
\text{perform predictor step:} \vspace{2mm} \\
$\underline{p}_{\textup{pred}}=\underline{p}_n+(-\Delta t) \dot{\underline{p}}_n+\frac{\Delta t^2}{2}(1-2\beta_p)\ddot{\underline{p}}_n$\vspace{2mm}\\
$\dot{\underline{p}}_{\textup{pred}}=\dot{\underline{p}}_n+(1-\gamma_p)(-\Delta t) \ddot{\underline{p}}_n$ \vspace{2mm}\\
$ \mathbf{\bar{F}}_{n-1}=2\mathbf{M}^{D}(\underline{u}_{\, n-1}-\underline{u}_{\, d, n-1})$ \vspace{2mm} \\
$\kappa=0$ \vspace{2mm}  \\
$p^\kappa_{n-1}=\underline{p}_{\textup{pred}}$ \vspace{2mm}\\
$\dot{p}^\kappa_{n-1}=\dot{\underline{p}}_{\textup{pred}}$\vspace{2mm}\\
$\textup{Res}\,
=~\mathbf{M} \ddot{\underline{p}}^{\kappa}_{n-1}-\mathbf{C} \dot{\underline{p}}^{\kappa}+\mathbf{K} \underline{p}^{\kappa}-\mathbf{\bar{F}}_{n-1}-\bm{\mathcal{T}}(u_{d, n-1})\ddot{\underline{p}}^{\kappa}_{n-1} - \mathbf{{A}}_1\dot{\underline{p}}^{\kappa}_{n-1} +\mathbf{{A}}_2\ddot{\underline{p}}^\kappa_{n-1}$\vspace{2mm} \\
\While{$\textup{Res}> \textup{TOL}_p$ \vspace{1mm}}{
\text{solve the algebraic system of equations:} \vspace{2mm} 
\begin{align*}
&\hspace{-15mm}\mathbf{\overline{M}}  \, \ddot{\underline{p}}^{\kappa+1}_{n-1}= \mathbf{\bar{F}}_{n-1}+\mathbf{C} \dot{\underline{p}}_{\textup{pred}}-\mathbf{K} \underline{p}_{\textup{pred}}+\bm{\mathcal{T}}(u_{d, n-1})\ddot{\underline{p}}^{\kappa}_{n-1} +\mathbf{{A}}_1\dot{\underline{p}}^\kappa_{n-1} -\mathbf{{A}}_2\ddot{\underline{p}}^\kappa_{n-1}
\end{align*}
\hspace{-1.7mm}\text{perform corrector step:} \vspace{2mm}  \\
$\underline{p}^{\kappa+1}_{n-1}=\underline{p}_{\textup{pred}}+\beta_p \Delta t^2 \ddot{\underline{p}}^{\kappa+1}_{n-1}$\vspace{2mm}\\
$\dot{\underline{p}}^{\kappa+1}_{n-1}=\dot{\underline{p}}_{\textup{pred}}+\gamma_p (-\Delta t) \ddot{\underline{p}}^{\kappa+1}_{n-1}$\vspace{2mm}\\
$\textup{Res}\,
=~\mathbf{M} \ddot{\underline{p}}^{\kappa+1}_{n-1}-\mathbf{C} \dot{\underline{p}}^{\kappa+1}+\mathbf{K} \underline{p}^{\kappa+1}-\mathbf{\bar{F}}_{n-1}-\bm{\mathcal{T}}(u_{d, n-1})\ddot{\underline{p}}^{\kappa+1}_{n-1} - \hspace*{1cm} \mathbf{{A}}_1\dot{\underline{p}}^{\kappa+1}_{n-1} +\mathbf{{A}}_2\ddot{\underline{p}}^\kappa_{n-1}$ \vspace{2mm}\\
$\kappa \mapsto \kappa+1$ \\}
$n \mapsto n-1$, $t \mapsto t-\Delta t$\\
}
\end{algorithm} 
\section{Shape optimization algorithm} \label{Section_ShapeOptAlgorithm}
Having seen how to numerically solve the state and the adjoint problem, we can now move to the general form of our optimization algorithm with step-size control. \\

\begin{algorithm}[H] \label{shapeopt_algorithm}
\caption{Shape optimization algorithm}\label{algorithm} \vspace{2mm}
\textbf{input:} initial lens geometry $\Omega_l^0$ \\
 $s=0$\\
\While{$\Omega_l^s$ is not optimal  \textbf{\textup{and}} $s<s_{\textup{max}}$\vspace{1mm} \label{line_optimal}}{
\text{solve the state problem on the domain with lens $\Omega_l^s$ \vspace{1mm}}\\
\text{solve the adjoint problem on the domain with lens $\Omega_l^s$ }\vspace{1mm}\\
 \text{compute the shape gradient}\label{step1} \vspace{1mm}\\
\text{update geometry $\Omega_l^s \mapsto \Omega_l^{s+1, \text{new}}$ with a descent step $\alpha$}\label{step2}\vspace{1mm} \\ 
\uIf{$J(\Omega_l^{s+1, \mathrm{new}}) \leq J(\Omega_l^{s})$ \vspace{1mm} \label{geom_update}}
{accept the new geometry, $\Omega_l^{s+1}=\Omega_l^{s+1, \text{new}}$ \vspace{1mm} \\
\If{previous step did not need to be repeated \vspace{1mm}}{make the step size larger}
proceed to the next step, $s \mapsto s+1$ \vspace{1mm} \\}
\Else {decrease the step size \vspace{1mm}\\
\uIf{$\mathrm{step size} < \textup{TOL}_{\mathrm{step}}$ \vspace{1mm}}{stop, accept the current geometry as the final one}
\Else{\textbf{goto} \ref{step2}, repeat with the current descent step \vspace{1mm}} 
}\label{end}}
\end{algorithm} 
\vspace{4mm}
\indent The optimality condition in line \ref{line_optimal} will in practice be realized by testing if the norm of the shape gradient at the current step is sufficiently small, i.e. below a certain tolerance level $\textup{TOL}_{\textup{grad}}$, relative to the gradient norm at the first step. We also introduce $s_{\textup{max}}$ as the maximal number of gradient steps. \\
\indent Lines \ref{geom_update} to \ref{end} in Algorithm \ref{shapeopt_algorithm} describe the step size control, where the new geometry is accepted if it is better than the previous one in the sense of the value of the cost function. If not, the update is repeated with a smaller descent step. We introduce a lower step size bound $\textup{TOL}_{\text{step}}$ so that the search for a smaller step can be stopped after a certain number of attempts when the cost and gradient norm don't decrease significantly anymore. However, if the newly obtained geometry is better than the previous one and there was no repetition of the previous step, we allow the current descent step to become larger. 
\subsection{Updating geometry}
\indent Let us reflect on step \ref{step2} of Algorithm \ref{algorithm} where the geometry is updated. Within the isogeometric approach shapes are easily changed by moving control points. In our solver, we restrict ourselves to a 2D setting, where control points corresponding to the lens boundary are treated as design variables. We denote by $\mathcal{I}_{\partial\Omega_l}$ the set of indices of degrees of freedom contributing to the boundary of the lens, i.e. $$\mathcal{I}_{\partial\Omega_l}=\lbrace \, i\in \mathcal{I}~|~ \left.N_{i,q}\right|_{\partial\Omega_l}\neq \boldsymbol{0} \, \rbrace.$$
\indent In every optimization step, control points $\mathcal{C}_i=(x_i, y_i)$, $i \in \mathcal{I}_{\partial \Omega_l}$, are moved in the $y$-direction (while they remain fixed in the $x$-direction), according to
\begin{align} \label{moving_control_points}
y^{\mathrm{new}}_{i}=y_i-\alpha \, dJ(u^h, p^h, \Omega_l)(\left.N_{i, q}\right|_{\partial\Omega_l}e_y), 
\end{align}
where $\alpha$ is the current descent step and $e_y=(0,1)^T$. Since we identify degrees of freedom on the interfaces, the fluid boundaries common with the lens ones are moved in the same way, whereas the others remain fixed. \\
\indent In \eqref{moving_control_points}, we emphasized the dependence of the shape gradient on the solution of the adjoint problem. When $u^h$ is the exact solution of the state and $p^h$ of the adjoint problem, then $dJ(u^h, p^h, \Omega_l)\Theta=dJ(u, \Omega_l)\Theta$. 

\subsubsection{Extending the boundary changes to the interior}
After updating the lens boundary, the internal control points of the lens and other patches are recomputed using the Coons patch approach~\cite{Qian}. \\
\indent Let us denote the updated control points by $\mathcal{C}^{\mathrm{new}}_i=(x_i,y^{\mathrm{new}}_i)$. The updated boundary curves of, for instance, the lens patch are given by 
\begin{align*}
	l_{\textup{top}}(\hat{x}_1) =\sum_{i\in\mathcal{I}_{\partial\Omega_l}}\mathcal{C}_i^{\mathrm{new}}\hat{N}_{i,q}(\hat{x}_1,1), & \hspace{1cm}	l_{\textup{bottom}}(\hat{x}_1) =\sum_{i\in\mathcal{I}_{\partial\Omega_l}}\mathcal{C}_i^{\mathrm{new}}\hat{N}_{i,q}(\hat{x}_1,0), \\
	l_{\textup{left}}(\hat{x}_2) =\sum_{i\in\mathcal{I}_{\partial\Omega_l}}\mathcal{C}_i^{\mathrm{new}}\hat{N}_{i,q}(0,\hat{x}_2), & \hspace{1cm}	l_{\textup{right}}(\hat{x}_2) =\sum_{i\in\mathcal{I}_{\partial\Omega_l}}\mathcal{C}_i^{\mathrm{new}}\hat{N}_{i,q}(1,\hat{x}_2),	
\end{align*} 
each mapping the parametric interval $[0,1]$ to the respective boundary in the physical space. A Coons surface is then computed as a combination of the following three surfaces:\\
\begin{itemize}
\item Two ruled NURBS surfaces by linear interpolation of opposite boundary curves:
\begin{align*}
		L_1(\hat{x}_1,\hat{x}_2) &=(1-\hat{x}_2)\, l_{\textup{bottom}}(\hat{x}_1)+\hat{x}_2\, l_{\textup{top}}(\hat{x}_1),\\
		L_2(\hat{x}_1,\hat{x}_2) &=(1-\hat{x}_1)\, l_{\textup{left}}(\hat{x}_2)+\hat{x}_1\, l_{\textup{right}}(\hat{x}_2).
		\end{align*}
	\item The bilinear NURBS surface interpolating between the corner points:
	\begin{align*}
			B(\hat{x}_1,\hat{x}_2) =& \,(1-\hat{x}_1)(1-\hat{x}_2)\ G(0,0) + (1-\hat{x}_1)\hat{x}_2\,G(0,1)  \\
		 & +\,\hat{x}_1(1-\hat{x}_2)\,G(1,0)+ \hat{x}_1\hat{x}_2\,G(1,1).
	\end{align*}
\end{itemize}
 The new patch mapping $G^{\mathrm{new}}$ is then defined as
	\begin{equation*}
		G^{\mathrm{new}}(\hat{x}_1,\hat{x}_2) = L_1(\hat{x}_1,\hat{x}_2) + L_2(\hat{x}_1,\hat{x}_2) - B(\hat{x}_1,\hat{x}_2).
	\end{equation*}
For simplicity, we have omitted the index $(j)$ of the current patch in notation above. For instance, $G$ is used instead of $G^{(j)}$ to denote the current NURBS mapping. We note that the construction of the Coons patch can be motivated purely geometrically and is hence independent of the chosen coordinate frame.

\remark{In the present work we do not employ tools to avoid mesh tangling, since in our numerical experiments this unwanted effect was never present. One possibility would be to use relative positioning when control points are moved only in the $y$-direction as suggested in \cite{Simeon}}.
 \section{Numerical experiments} \label{Section_NumericalExperiments}
 
 We will now test our algorithm in two different scenarios. In the first one, we validate our code by creating artificial pressure data with a known goal shape for the lens and then trying to get back to that goal shape by starting the optimization process with a different lens geometry. In the second scenario, the goal shape will not be known and the target pressure distribution will be given analytically.\\
 \indent In both cases, we will perform experiments with linear B-splines (where IGA and classical FEM coincide) as well as quadratic splines. In all the experiments we will keep the lens corner point fixed, i.e. the shape will be pinned at the right corner. \\
 \indent The numerical results presented here have been obtained with the help of the GeoPDEs package in MATLAB~\cite{Falco}.

 \subsection{Coefficient values} For all the performed tests, the fluid is taken to be water. The values of physical parameters in typical media can be found, for instance, in \cite{Folds, Rossing}. The coefficients used in the experiments are provided in Table \ref{table_coeff}.

 \begin{center}
 	\renewcommand{\arraystretch}{1.6}
 	\captionof{table}{Physical parameter values} \label{table_coeff} 
 	\begin{tabular}{||c c ||} 
 		\hline
 		fluid & lens   \\ 
 		\hline\hline
 		$c_f=1500 \, \frac{m}{s}$ & $c_l=1100 \, \frac{m}{s}$    \\
 		\hline
 		$b_f=6\cdot 10^{-9} \, \frac{m^2}{s}$ & $ b_l=4\cdot 10^{-9} \, \frac{m^2}{s}$    \\ 
 		\hline
 		$\varrho_f=1000 \, \frac{kg}{m^3}$ & $\varrho_l=1250 \, \frac{kg}{m^3}$    \\ 
 		\hline
 		$B/A=5$  & $B/A=4$    \\ 
 		\hline
 	\end{tabular}\vspace{2mm} \\
 \end{center}
 
\noindent The sound diffusivity in fluid $b_f$ and in lens $b_l$ are quite small relative to the amplitudes in the acoustic pressure field. Their influence on the damping of the wave was proven to be neglectable in all our numerical experiments. \\
  \subsection{Excitation} 
\begin{figure}[h]
 	\includegraphics[trim = 9cm 0cm 7cm 0cm, clip, scale = 0.16]{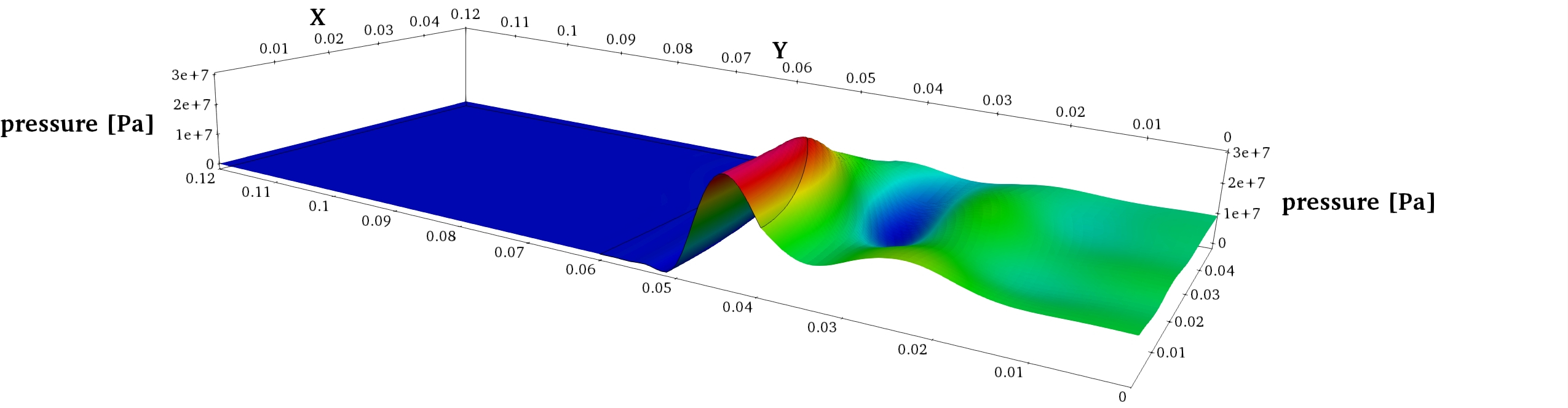}\parbox{2cm}{~~}\\
 	\includegraphics[trim = 9cm 0cm 15cm 0cm, clip, scale = 0.16]{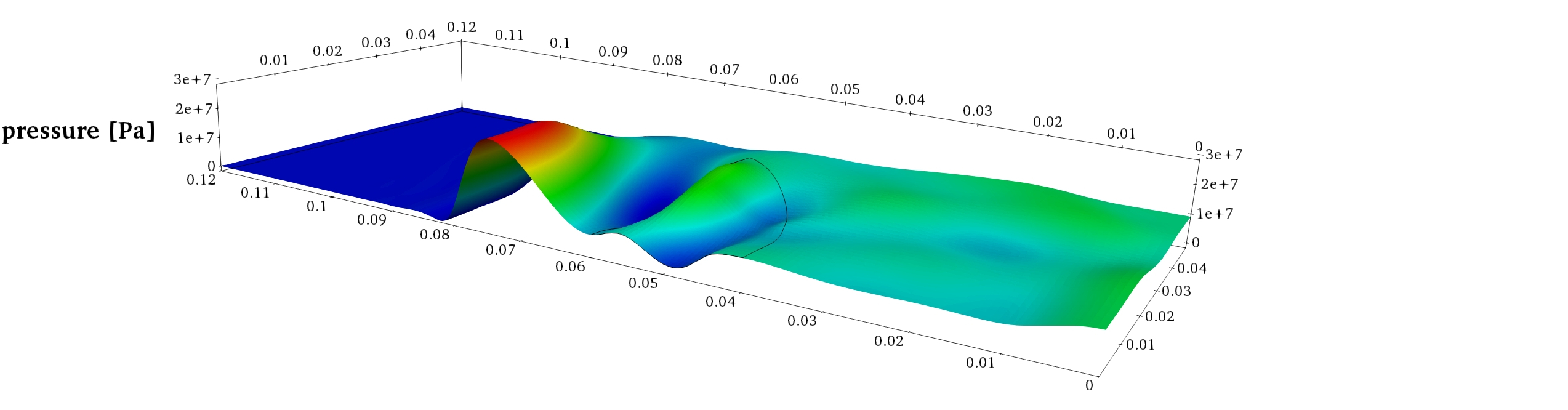}\hspace*{-5mm}\includegraphics[trim = 0cm 0cm 0cm 0cm, clip, scale = 0.16]{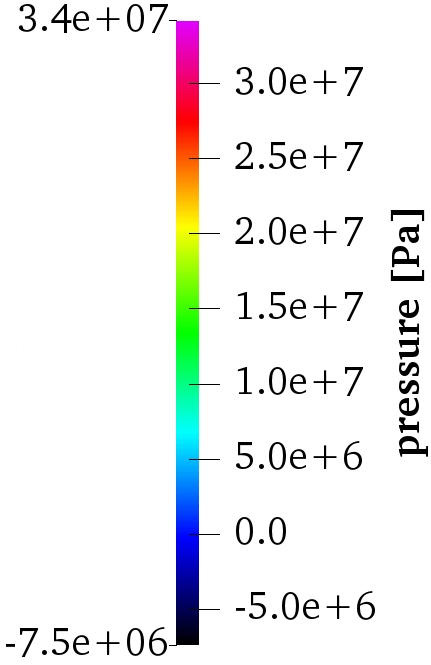}\\
 	\includegraphics[trim = 9cm 0cm 7cm 0cm, clip, scale = 0.16]{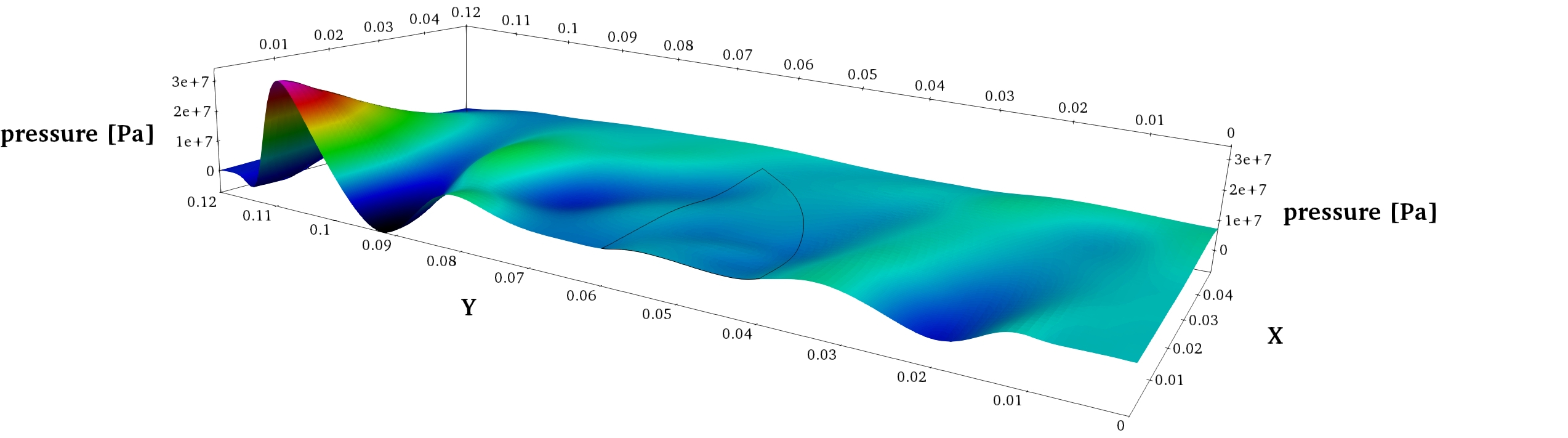}\parbox{2cm}{~~}
 	\caption{Nonlinear wave propagation in the lens setting with steepening toward the end of the geometry. The black lines in the picture show the position of the lens. \label{fig:3DPicture}}
 \end{figure}

For the source condition, we consider a modulated sinusoidal wave (cf. \cite[Chapter 2]{Lee})
 \begin{align*}
 	g=g_0\,\displaystyle e^{-\left(\omega\, t/8\right)^2}\, \sin(\omega t).
 \end{align*}
In this way our waveform resembles the short tone bursts used in medical ultrasound practices, while still being continuous in time and space. We set the source amplitude value to be $g_0=4\cdot 10^9\,\textup{Pa}$ and frequency $f=70$kHz. The angular frequency is then given by $\omega=2\pi f$. \\
 \indent Figure~\ref{fig:3DPicture} shows the propagation of the signal in an exemplary lens setting. Toward the end of the geometry the steepening due to the nonlinear effects can be observed, which is also where the focal point is assumed to be and where high pressure values are required.
 \clearpage
 
\subsection{Test 1: Exact solution}
Since the analytical solution of our problem is not known, we will generate artificial data for testing purposes. Using the same solver to generate data and then to invert the problem is known as committing the "inverse crime" \cite{ColtonKress}, since it can cause unrealistically good solutions. To avoid the inverse crime, we produce the synthetic pressure data on a domain with lens shape $\Omega_l^\star$ by solving the forward problem on a staggered grid (compared to the grid that is used for the optimization process) and pollute it with random Gaussian noise of 2\,\% of the maximal pressure amplitude. $\Omega_l^\star$ is then taken to be the exact lens shape and the artificially produced data the target pressure distribution $u_d$.

\subsubsection{Domain measures}
In the forthcoming numerical experiments we will use several lens shapes as initial $\Omega_l^0$ and goal domains $\Omega_l^{\star}$. Table \ref{table_domain_size} lists all the specific measurements.
\clearpage
\begin{center}
	\renewcommand{\arraystretch}{1.6}
	\captionof{table}{Lens measurements} \label{table_domain_size} 
	\begin{tabular}{|| c | c  | c  | c | c ||} 
		\hline
		$\Omega_{l,\textup{upper straight}}$& $\Omega_{l,\textup{upper curved}}$& $\Omega_{l,\textup{both perturbed}}$& $\Omega_{l,\textup{both down}}$& $\Omega_{l,\textup{Gauss}}$\\
		\hline\hline
		$L=0.12$\,m & $L=0.12$\,m & $L=0.12$\,m & $L=0.12$\,m & $L=0.12$\,m \\
		\hline
		$B=0.05$\,m & $B=0.05$\,m & $B=0.05$\,m & $B=0.05$\,m & $B=0.05$\,m \\
		\hline
		$K=0.06$\,m & $K=0.06$\,m & $K=0.06$\,m & $K=0.06$\,m & $K=0.06$\,m\\
		\hline
		$W=0.04$\,m & $W=0.04$\,m & $W=0.04$\,m & $W=0.04$\,m  & $W=0.04$\,m  \\
		\hline
		$P=0.02$\,m & $P=0.015$\,m & $P = 0.016$\,m & $P=0.021$\,m & $P=0.025$\,m  \\
		\hline
		$S=0.09$\,m & $S=0.09$\,m & $S=0.09 $\,m & $S=0.09$\,m & $S=0.09$\,m   \\
		\hline
		$R=0.04$\,m & $R=0.04$\,m & $R= 0.042$\,m & $R=0.037$\,m & $R=0.035$\,m \\
		\hline
	\end{tabular}\vspace{2mm}
\end{center}
The geometry parameters are diagrammed in Figure~\ref{fig:Patches}. The curved lens boundaries are always chosen to be arc segments with horizontal tangents at the point belonging to the axis of symmetry. Unless stated otherwise, the subregion $D$, where the pressure should match $u_d$, is taken to be patch number $6$ of our computational domain.

\subsubsection{Spatial discretization}
Table~\ref{spatial_discretization} specifies the number of degrees of freedom used in the numerical experiments.
\begin{center}
	\renewcommand{\arraystretch}{1.6}
	\captionof{table}{Spatial grid sizes} \label{spatial_discretization} 
	\begin{tabular}{|| c | c||} 
		\hline
		\multicolumn{2}{||c||}{Spatial degrees of freedom}\\
		Linear NURBS & Quadratic NURBS \\
		\hline\hline
		$\textup{ndof}_x = 46$&$\textup{ndof}_x = 48$ \\
		\hline
		$\textup{ndof}_y = 181$ & $\textup{ndof}_y = 185$ \\
		\hline
		$\textup{ndof} = 7976$ & $\textup{ndof} = 8484$\\
		\hline
	\end{tabular}\vspace{2mm}
\end{center}
Here $\textup{ndof}_x$ denotes the number of degrees of freedom in $x$-direction and $\textup{ndof}_y$ in $y$-direction of the computational domain. The overall number of degrees of freedom is denoted by $\textup{ndof}$. Such a fine mesh is necessary to cope with the nonlinear effects. In the focal region, it amounts to about 35 elements in $y$-direction per fundamental wavelength. With this refinement, we have $36$ design control points for the upper as well as for the lower lens boundary.

\subsubsection{Time stepping} The time resolution as well as the generalized $\alpha$-parameters used when solving the state problem are given in Table \ref{table_num_param}. For the adjoint solver only the parameters $\gamma_p$ and $\beta_p$ are relevant. 
\clearpage

\begin{center}
	\renewcommand{\arraystretch}{1.6}
	\captionof{table}{Time integration and numerical parameter values} \label{table_num_param} 
	\begin{tabular}{|| c| c | c ||} 
		\hline
		Time discretization & Method parameter  & Tolerances \\ 
		\hline\hline
		Final time $T=90 \, \mu$s & $\gamma = 0.75, \beta = 0.45$  & $\textup{TOL}_u=10^{-6}$\\
		\hline
		$\textup{ndof}_t=3801$   & $\gamma_p=0.5, \beta_p = 0.25$ & $\textup{TOL}_p=10^{-8}$\\
		\hline
		$\Delta t=23.684\,$ns& $\alpha_m =1/2, \alpha_f = 1/3$ &$\textup{TOL}_{\textup{grad}}=10^{-4}$ \\
		\hline
	\end{tabular}\vspace{2mm}
\end{center}

\noindent Such a fine time step is again needed to resolve the nonlinear behavior. This refinement amounts to about 600 time samples per fundamental time period. Around 7 iterations per time step were necessary to achieve convergence in the forward problem with the tolerance $\textup{TOL}_u$. In the adjoint solver, typically 8 iterations were conducted to converge to a solution with the residual tolerance $\textup{TOL}_p$. \vspace{-1mm}

\subsubsection{Only upper lens boundary moving} 
In our first test, we fix the lower lens boundary to be already at its target position and only perturb the upper boundary to see if the algorithm can retrieve the goal shape in this simple setting. For the initial shape, we choose the one with a straight upper boundary $$\Omega_l^0 = \Omega_{l,\textup{upper straight}}.$$ The goal shape is given by $$\Omega_l^{\star} = \Omega_{l,\textup{upper curved}}.$$ This is a thinner shape than the initial one, so the upper boundary should move downward in order to reach its goal. The exact domain measurements are given in Table \ref{table_domain_size}. \\
\indent Let us introduce the discrete shape gradient $\nabla J(\Omega_l^s)$, where $s$ denotes the current gradient step,  as
\begin{equation*}
	\nabla J(\Omega_l^s) =\left(\begin{array}{c}
		\vdots\\
		dJ(u^h,\Omega_l^s)(\left.N_{i,q}\right|_{\partial\Omega_l}e_y)\\
		\vdots
	\end{array}\right)_{i\in\mathcal{I}_{\partial\Omega_l}},
\end{equation*}
containing the discrete sensitivities of all boundary control points. We are interested in the relative change of the cost functional $J(\Omega_l^s)/J(\Omega_l^0)$ and the relative change of the norm of the shape gradient $\|\nabla J(\Omega_l^s)\|/\|\nabla J(\Omega_l^0)\|$ over the course of optimization. \\
\indent A scaling factor, depending on the computational domain and the norm of the shape gradient, has to be introduced to the descent step size to ensure that geometry changes lead to a feasible lens shape, i.e. a shape which is still contained within the domain $\Omega$ and does not intersect itself. Its order of magnitude is taken to be $10^{-3}/\|\nabla J(\Omega_l^s)\|$. We take the lower bound for the descent step to be $\textup{TOL}_{\textup{step}}=10^{-8}$.\\
\begin{figure}[h]
	\includegraphics[trim = 0cm 7cm 0cm 7cm, clip, scale=0.35]{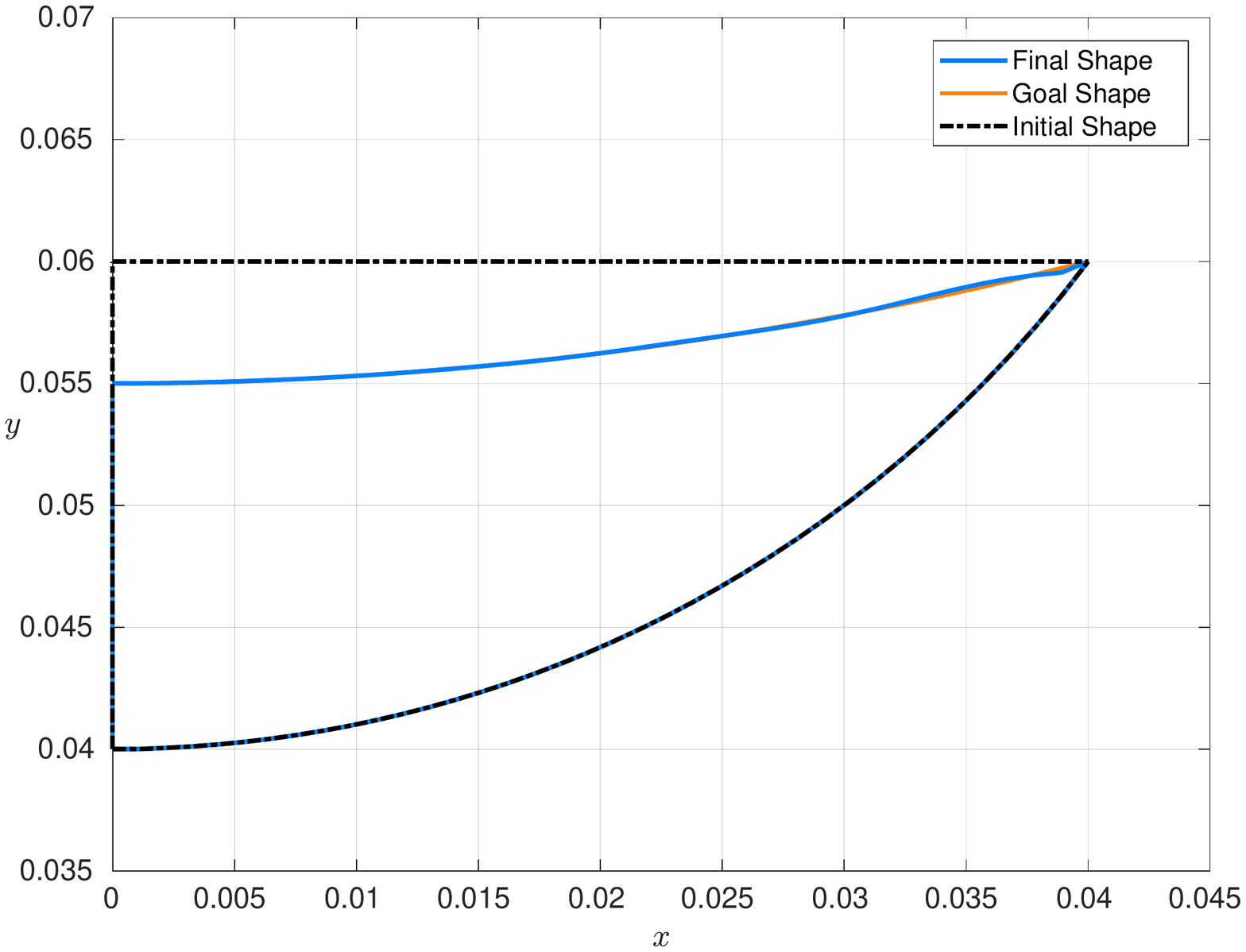}
	\caption{Final lens shape together with initial and goal shapes. \label{fig:Lens_overlay}}
\end{figure}
\begin{figure}[h]
	\raisebox{2.5mm}{\includegraphics[trim = 4cm 9cm 4.2cm 9cm, clip, scale=0.465]{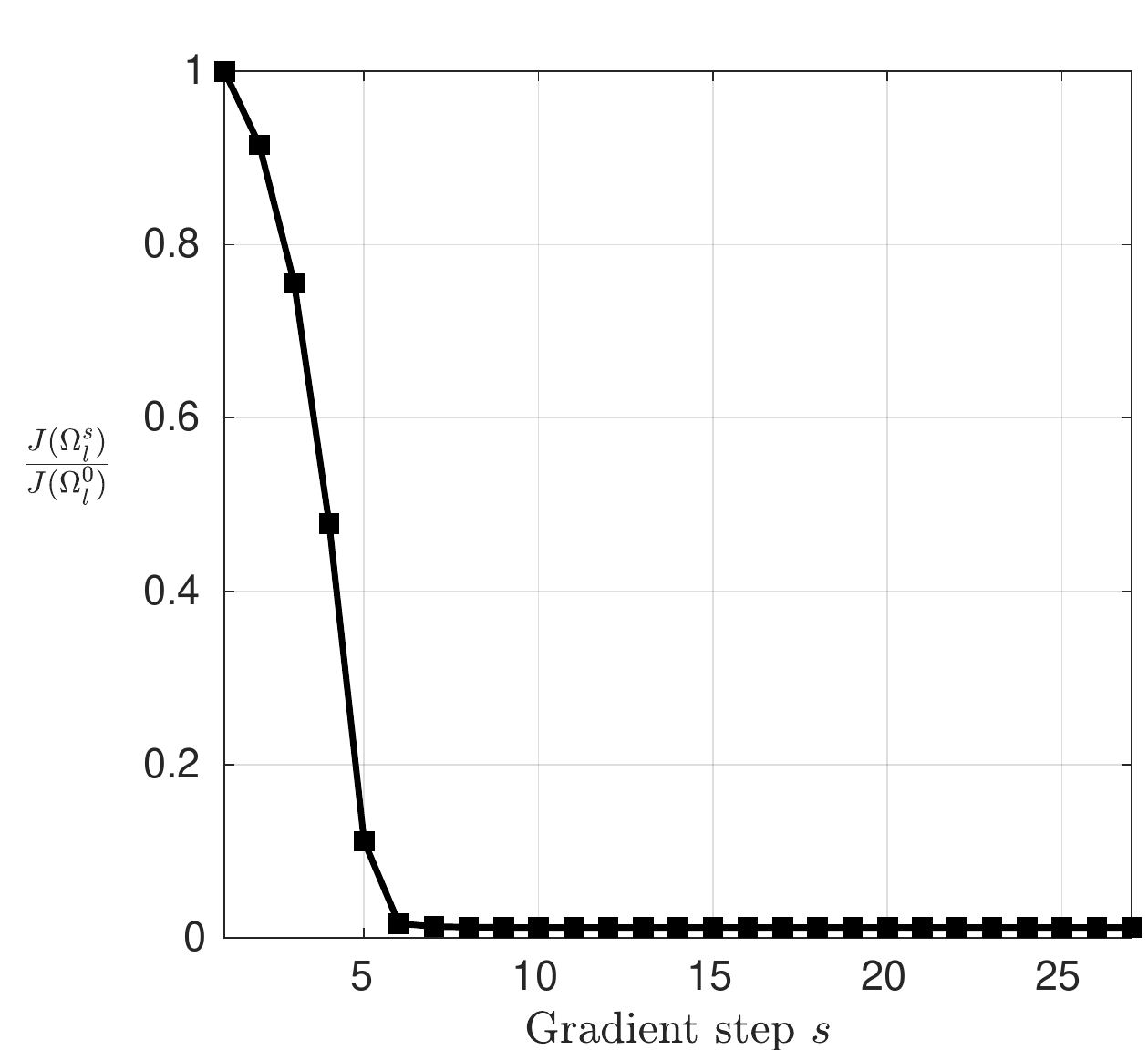}}
	\includegraphics[trim = 4cm 9cm 4.2cm 9cm, clip, scale=0.5]{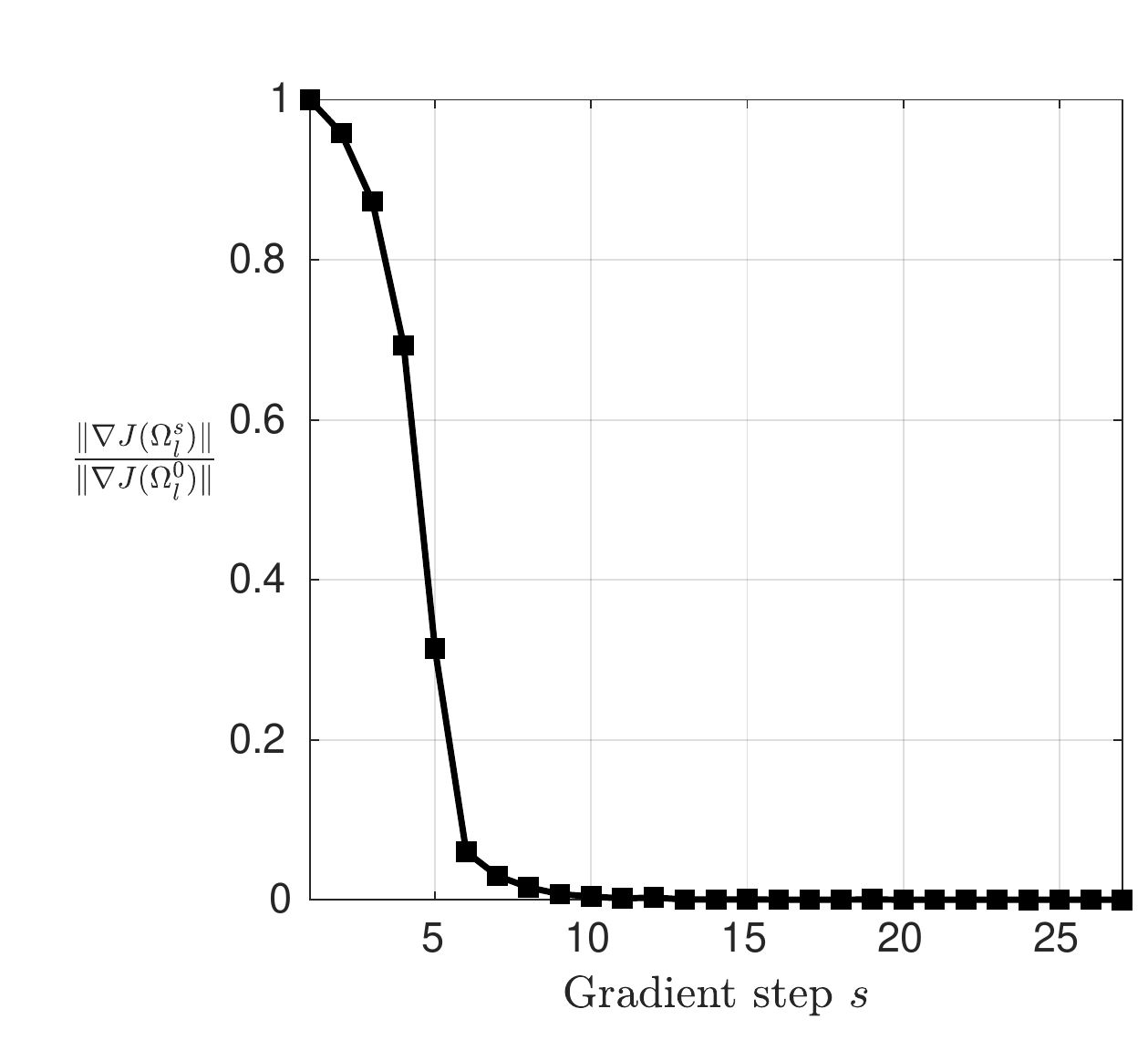}
	\caption{\textbf{left:} Relative cost change versus the number of gradient steps. \textbf{right:} Norm of the shape gradient. \label{fig:Cost_One}}
\end{figure}
\indent Figure~\ref{fig:Lens_overlay} depicts the final lens shape $\Omega_l^{\textup{final}}$ together with the initial shape $\Omega_l^0$ and the desired goal shape $\Omega^{\star}_l$. In Figure~\ref{fig:Cost_One}, the evolution of the relative change of the cost over the course of optimization is plotted as well as the relative change of the norm of the shape gradient. The cost, starting from $$J(\Omega_l^0)=2.2762\cdot10^{5}, $$ was reduced to $$J(\Omega_l^{\textup{final}})=2.7234\cdot 10^3,$$ corresponding to a decrease of around $98.8\,\%$. Similarly, the norm of the shape gradient was reduced by $99.99\,\%$ of its initial value $\|\nabla J(\Omega_l^0)\|=2.1886\cdot 10^7$. We observe that the cost function and the geometry change rapidly in the beginning and then very little closer to the optimal shape. Such behavior is common for gradient flows (see, for instance,~\cite{Dogan}). Furthermore, it can be observed that the cost values don't change significantly after the first 10 steps, hence a termination criterion that would stop the algorithm sooner can be employed. \\
\begin{figure}[h]
	\includegraphics[trim = 4cm 9cm 4.2cm 9cm, clip, scale=0.48]{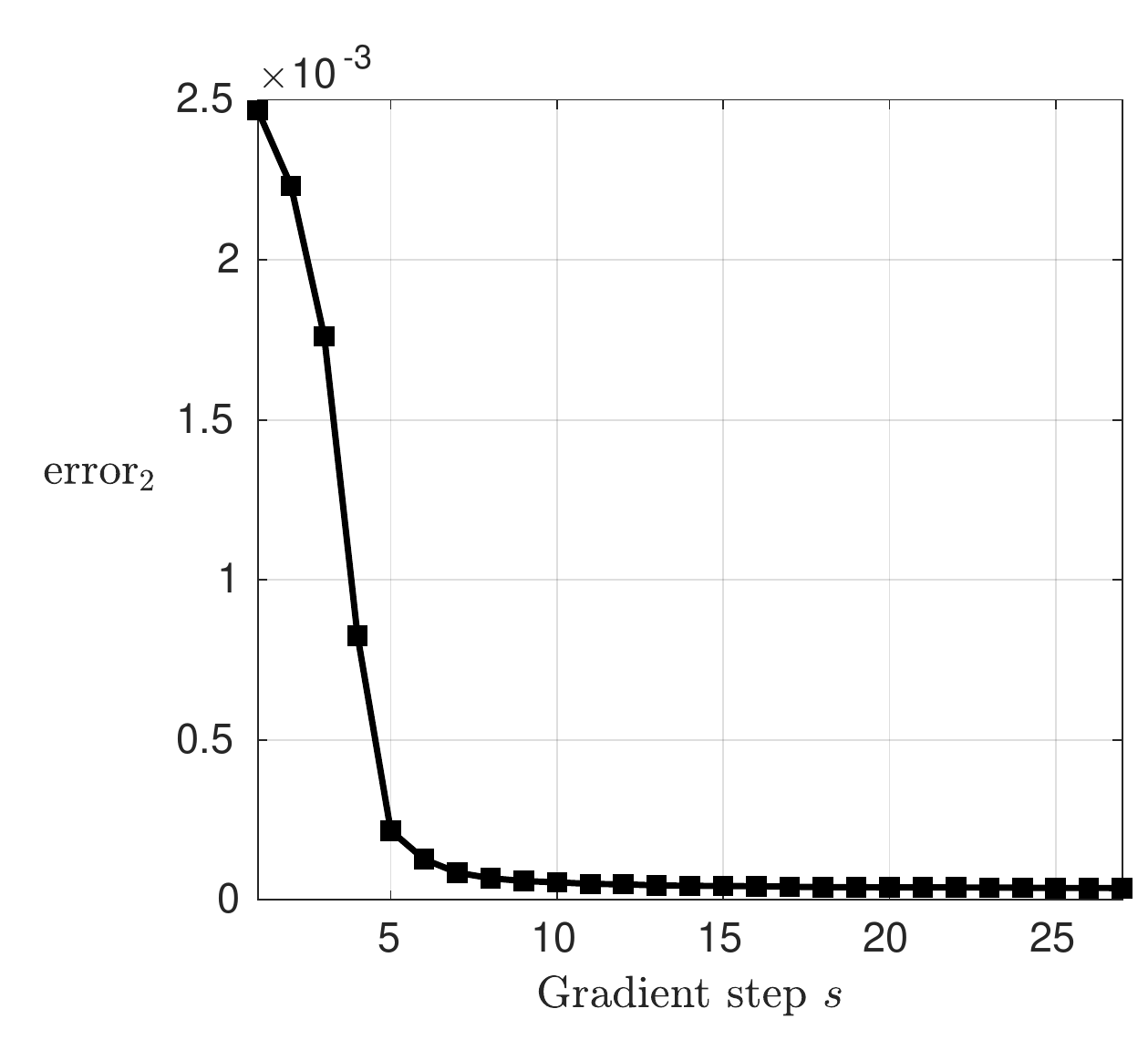}
	\caption{$L^2$-error over the course of optimization\label{fig:Lp_shape_errors}}
\end{figure}
\indent As a way of measuring errors in the lens shape, we introduce the mean $L^2$ shape error (cf. \cite{Park}):
\begin{align*}
	\textup{error}_{2} =\sqrt{\frac{1}{|\mathcal{I}_{\partial \Omega_l}|}\sum_{i\in\mathcal{I}_{\partial \Omega_l}}|y_{\textup{goal},i}-y_i|^2}.
\end{align*}
\noindent Figure~\ref{fig:Lp_shape_errors} shows the decay of this error versus the number of gradient steps. \\
\subsubsection{Both lens boundaries moving}
\begin{figure}[h]
	\vspace*{-5mm}
	\includegraphics[trim = 4cm 9cm 4.2cm 10cm, clip, scale=0.45]{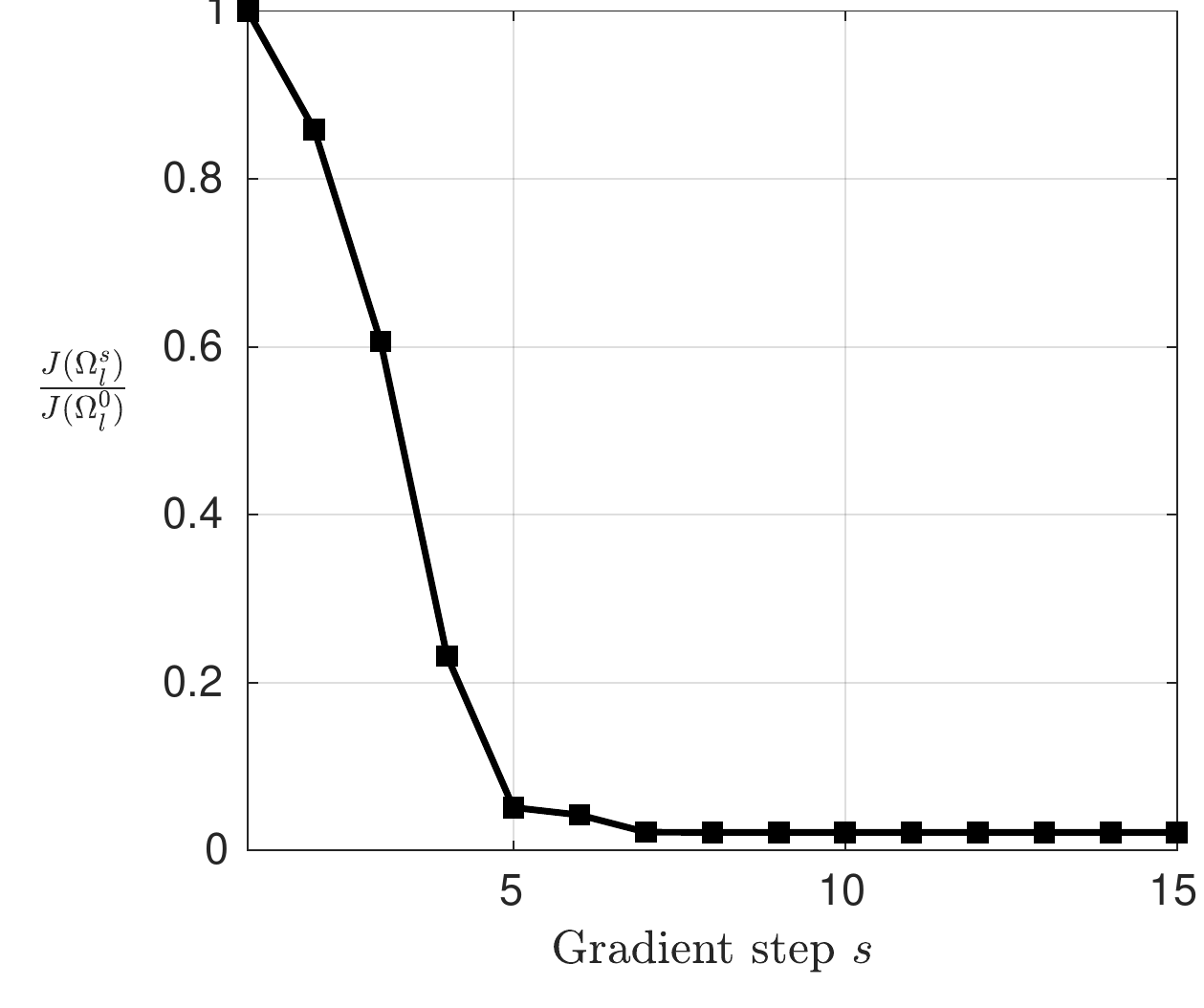}
	\includegraphics[trim = 4cm 9cm 4.2cm 10cm, clip, scale=0.45]{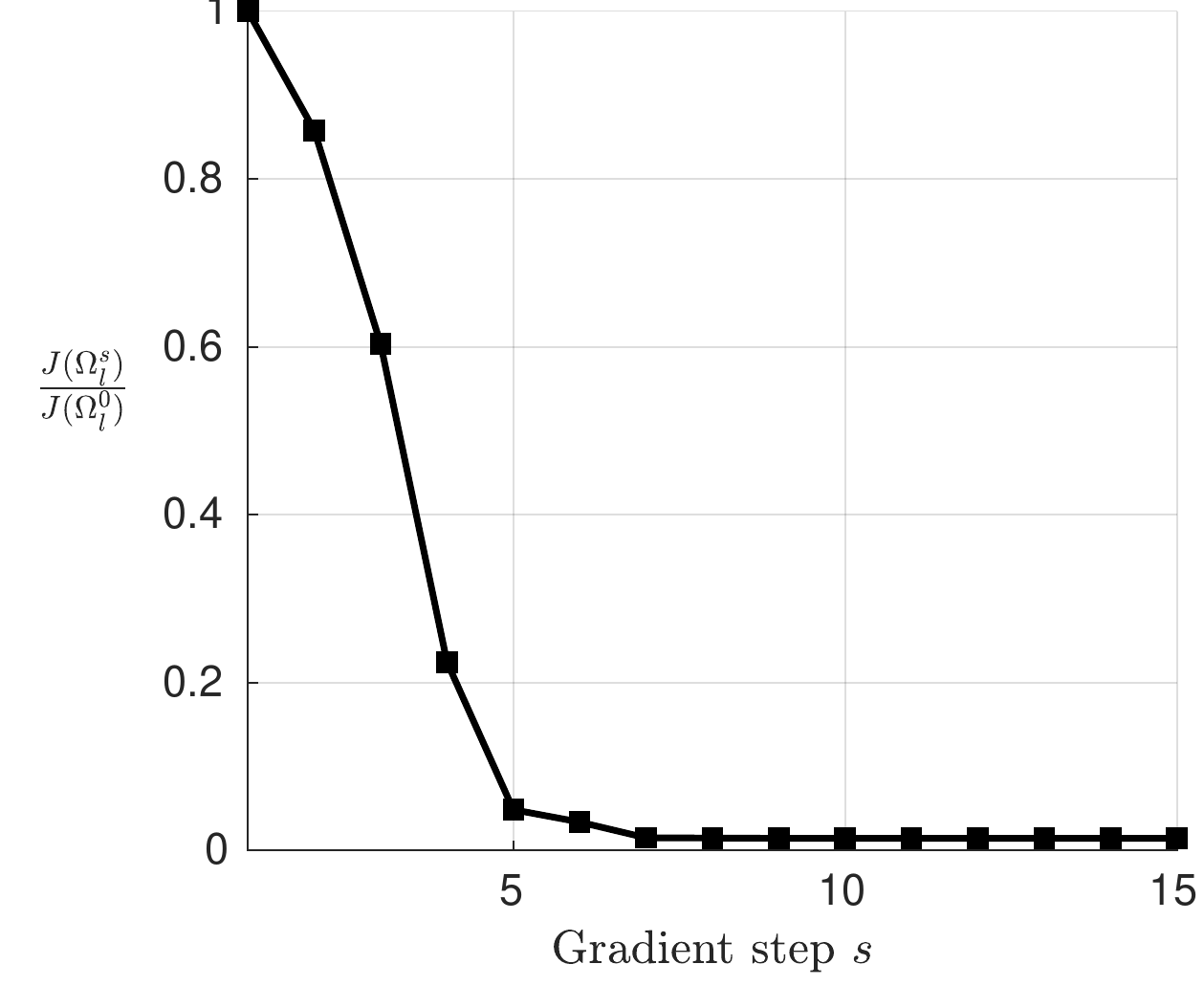}
	\caption{Relative change of the cost over the course of optimization. \textbf{left:} Linear NURBS. \textbf{right:} Quadratic NURBS. \label{fig:Cost}}
\end{figure}

\begin{figure}[h]
	\vspace*{-5mm}
	\includegraphics[trim = 4cm 9cm 4.2cm 10cm, clip, scale=0.45]{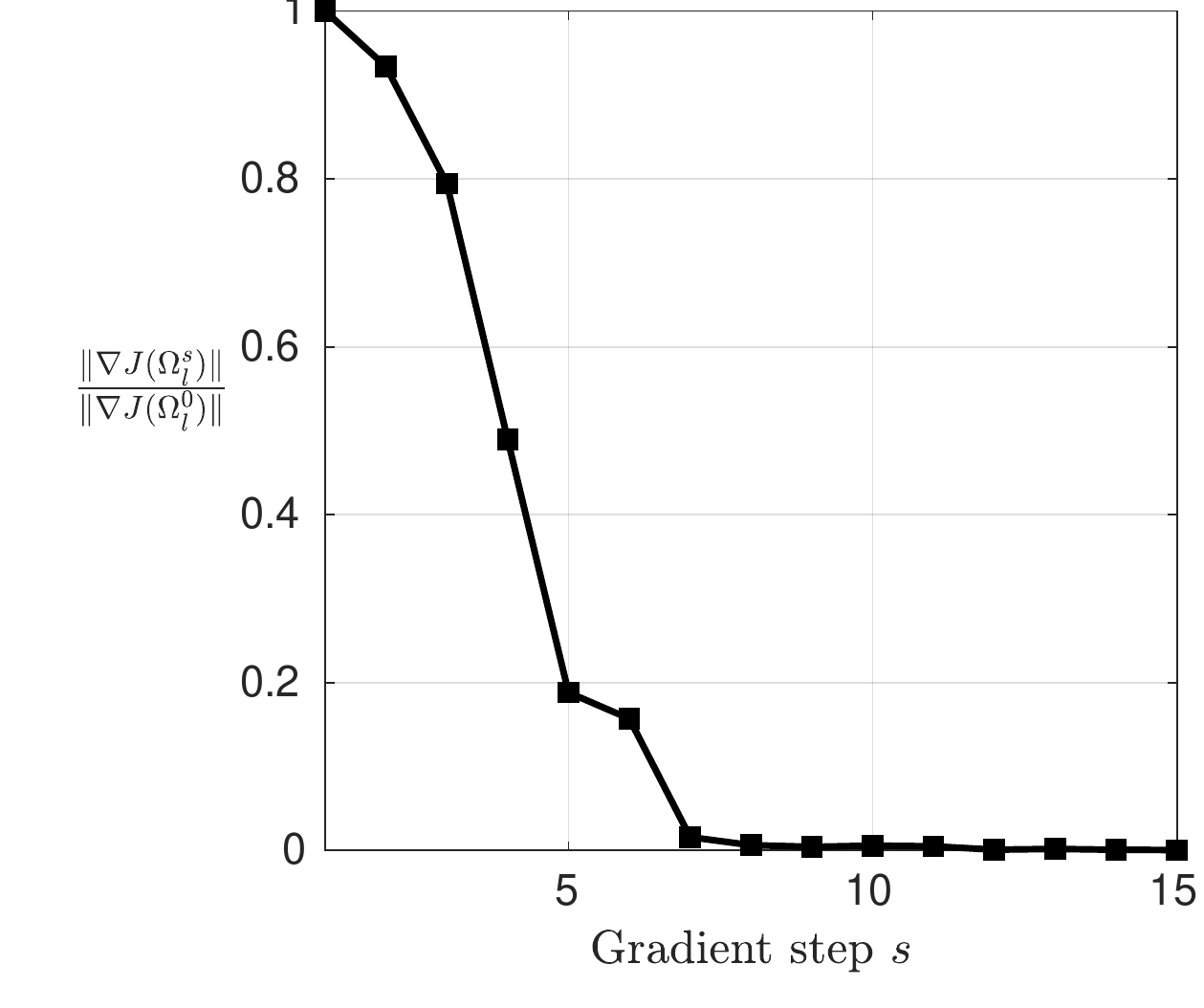}
	\includegraphics[trim = 4cm 9cm 4.2cm 10cm, clip, scale=0.45]{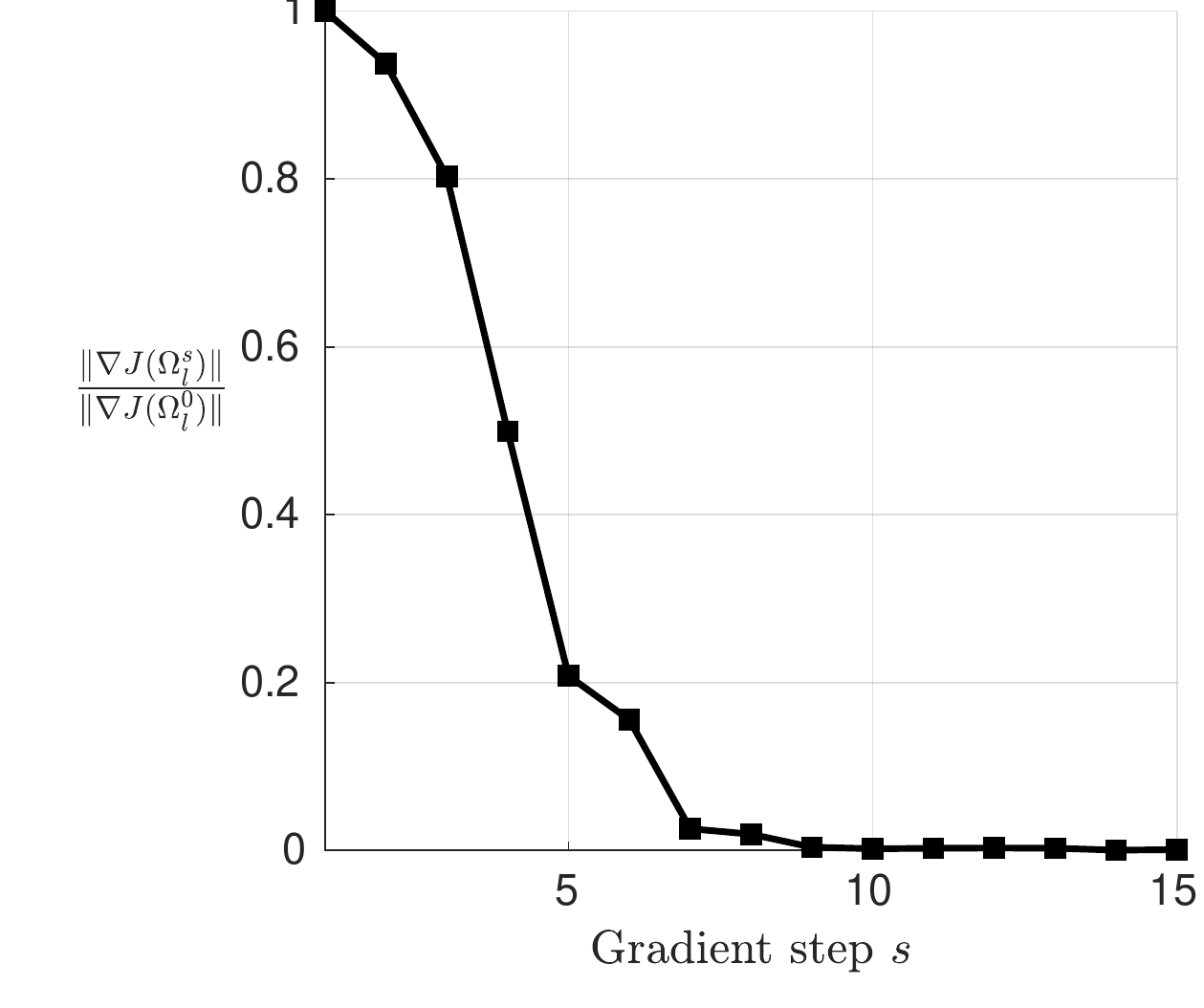}
	\caption{Norm of the shape gradient over the course of optimization. \textbf{left:} Linear NURBS. \textbf{right:} Quadratic NURBS. \label{fig:Norm_SG}}
\end{figure}

\begin{figure}[h]
	\vspace*{-5mm}
	\includegraphics[trim = 4cm 9cm 4.2cm 10cm, clip, scale=0.45]{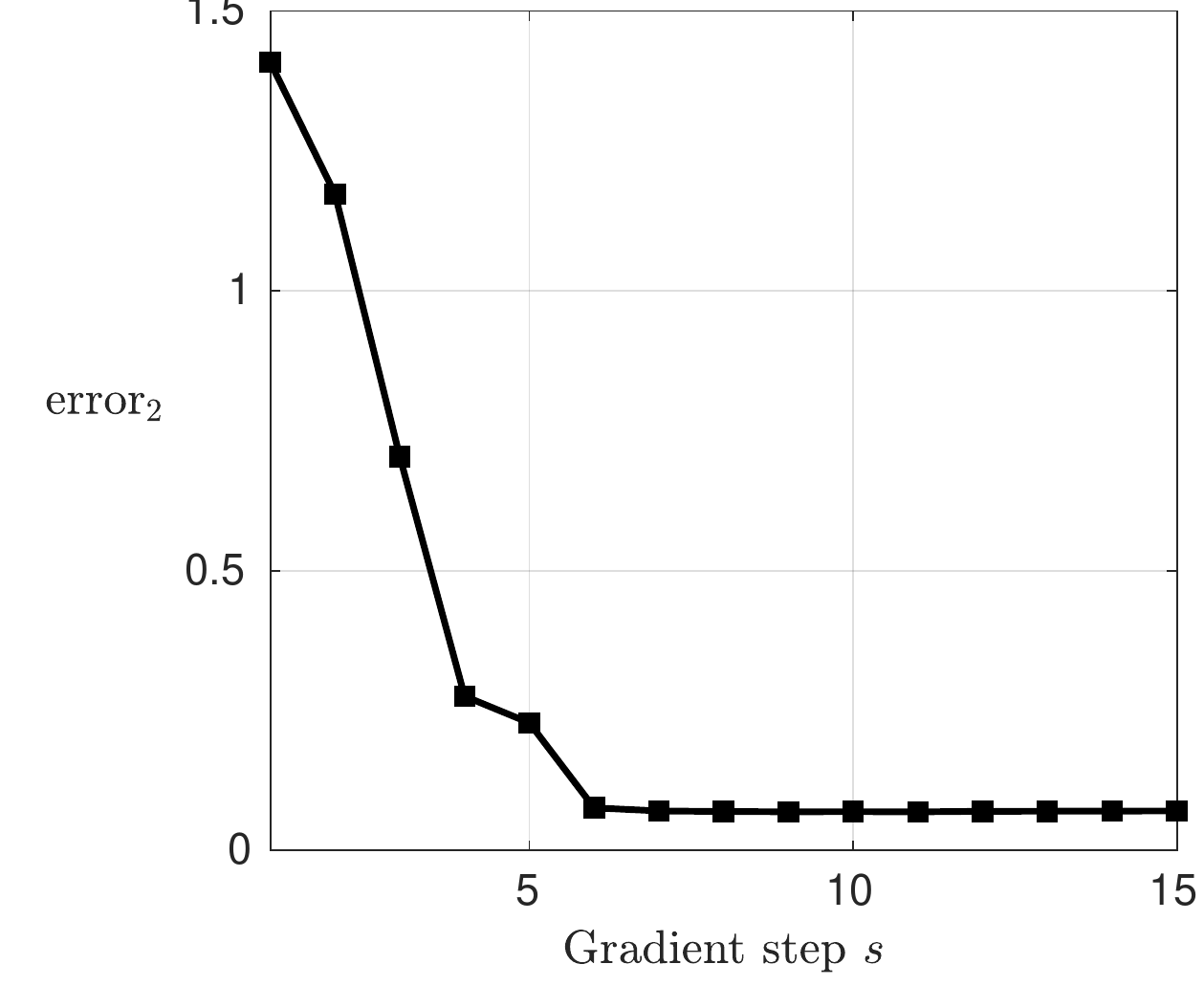}
	\includegraphics[trim = 4cm 9cm 4.2cm 10cm, clip, scale=0.45]{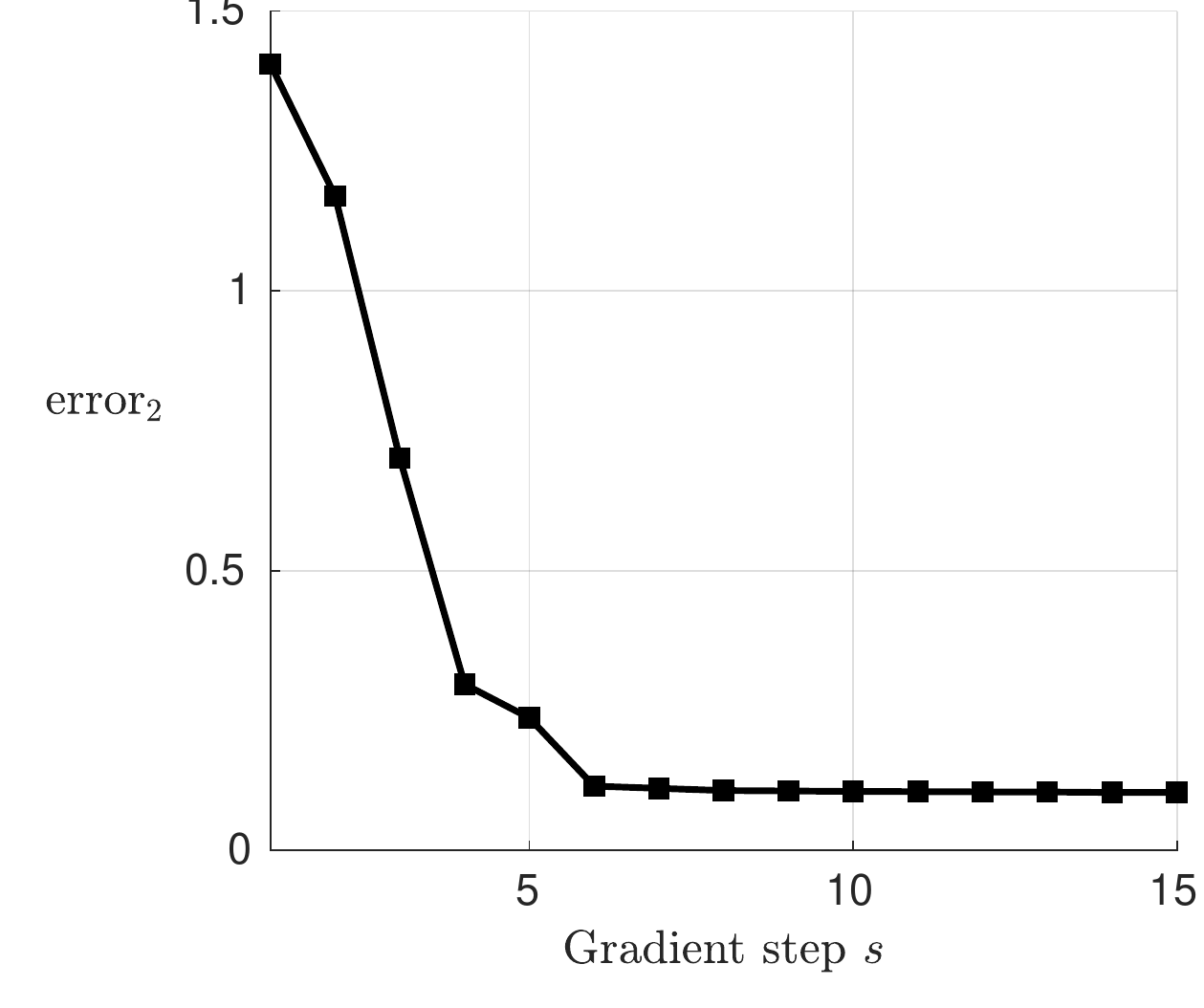}
	\caption{$L^2$ shape error over the course of optimization. \textbf{left:} Linear NURBS, \textbf{right:} Quadratic NURBS \label{fig:err_2}}
\end{figure}

For our second experiment, we allow \emph{both} the upper and lower lens boundary to move. We start with an initial lens geometry $$\Omega^0_l = \Omega_{l,\textup{both perturbed}},$$ where both lens boundaries are perturbed with respect to the goal shape. Unlike the previous case, now the initial shape has to expand in order the match the target lens. The exact domain measurements are again given in Table \ref{table_domain_size}.\\
\indent We will first, as before, perform tests where linear NURBS are used to obtain the synthetic data and then to solve the problem. Secondly, we will employ quadratic NURBS to obtain the artificial pressure data and then solve the shape optimization problem. \\
\begin{figure}[h]
	\includegraphics[trim = 0cm 7cm 0cm 7cm, clip, scale=0.3]{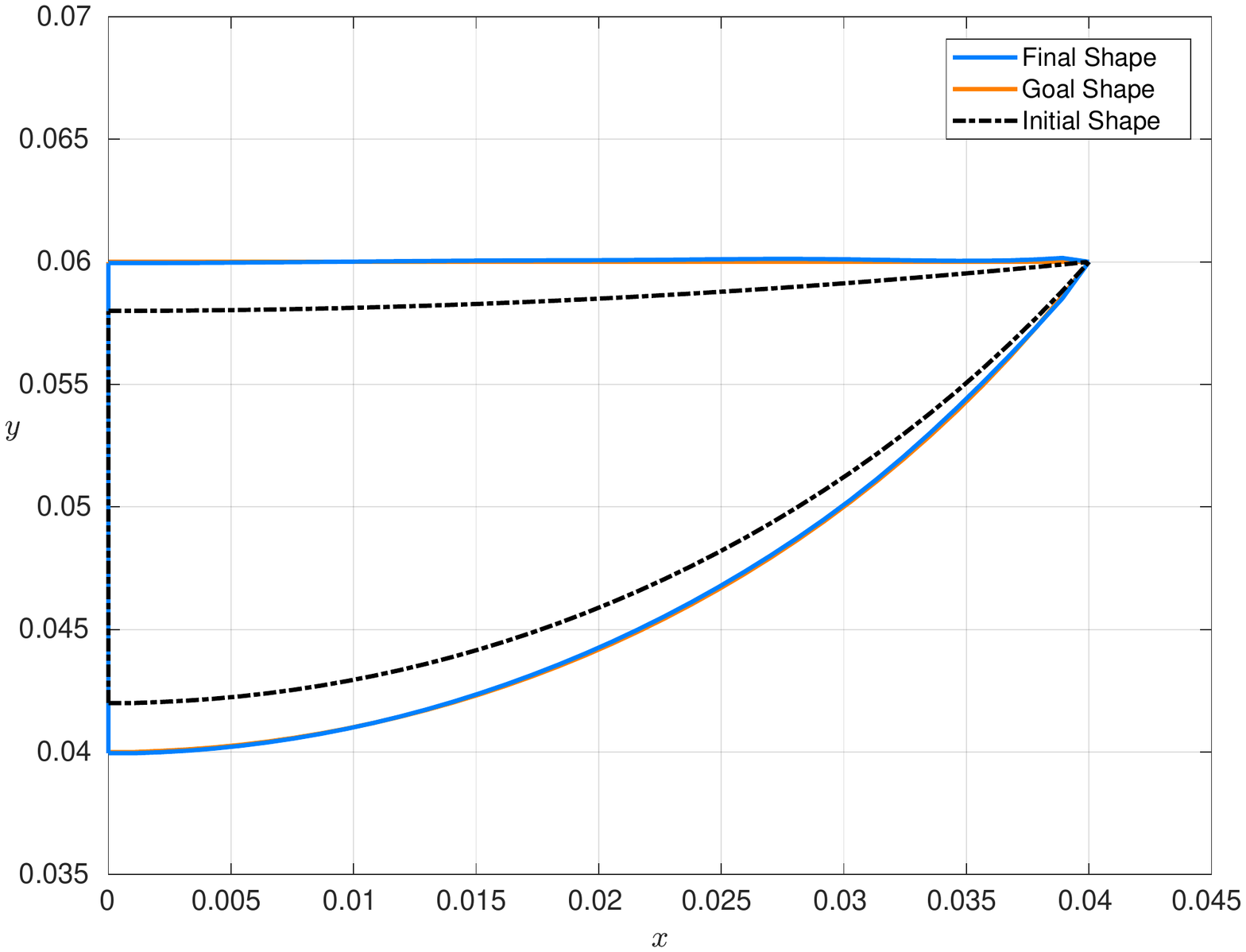}
	\includegraphics[trim = 0cm 7cm 0cm 7cm, clip, scale=0.3]{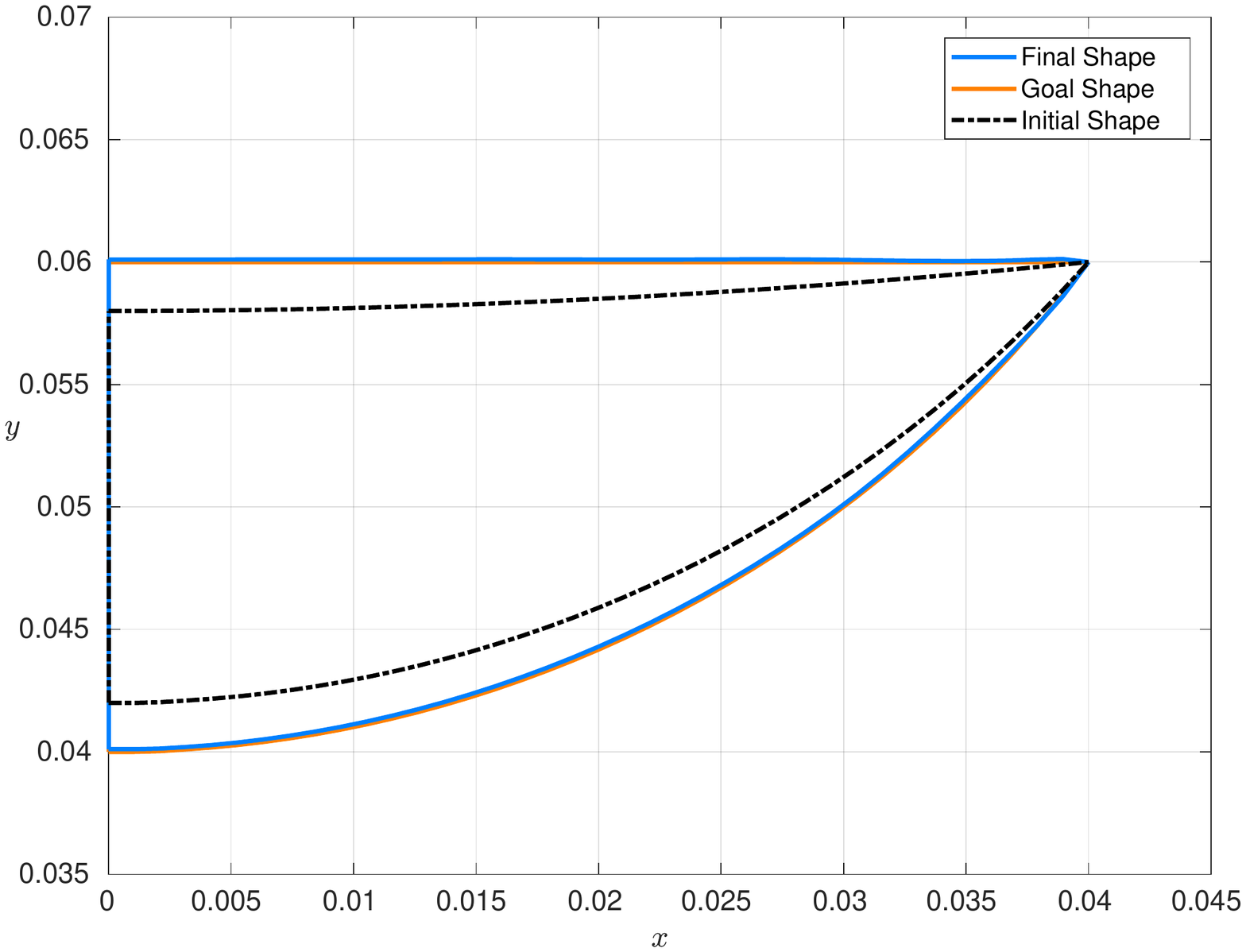}\\
	\caption{Final lens shape, together with initial and target shape. \textbf{left:} Linear NURBS. \textbf{right:} Quadratic NURBS. \label{fig:Lens_overlay2}}
\end{figure}
\indent In Figures~\ref{fig:Lens_overlay2}~-~\ref{fig:err_2}, we show the same quantities of interest as in the first example. They exhibit a similar behavior. The cost functional, its gradient as well as the shape errors decay rapidly within the first few steps, while almost no changes happen once the shape is close to its goal. When using the linear NURBS, the cost functional decreased by $97.85\,\%$ of its initial value, going from $J(\Omega_l^0)=1.532\cdot 10^5$ to $J(\Omega_l^{\textup{final}})=3.298\cdot 10^3$. 
The values are very similar in the quadratic NURBS case (see Figure~\ref{fig:Cost}). 

\subsubsection{Local, but not global minimum}
\begin{figure}[h]
	\includegraphics[trim = 0cm 6cm 1.5cm 6cm, clip, scale=0.33]{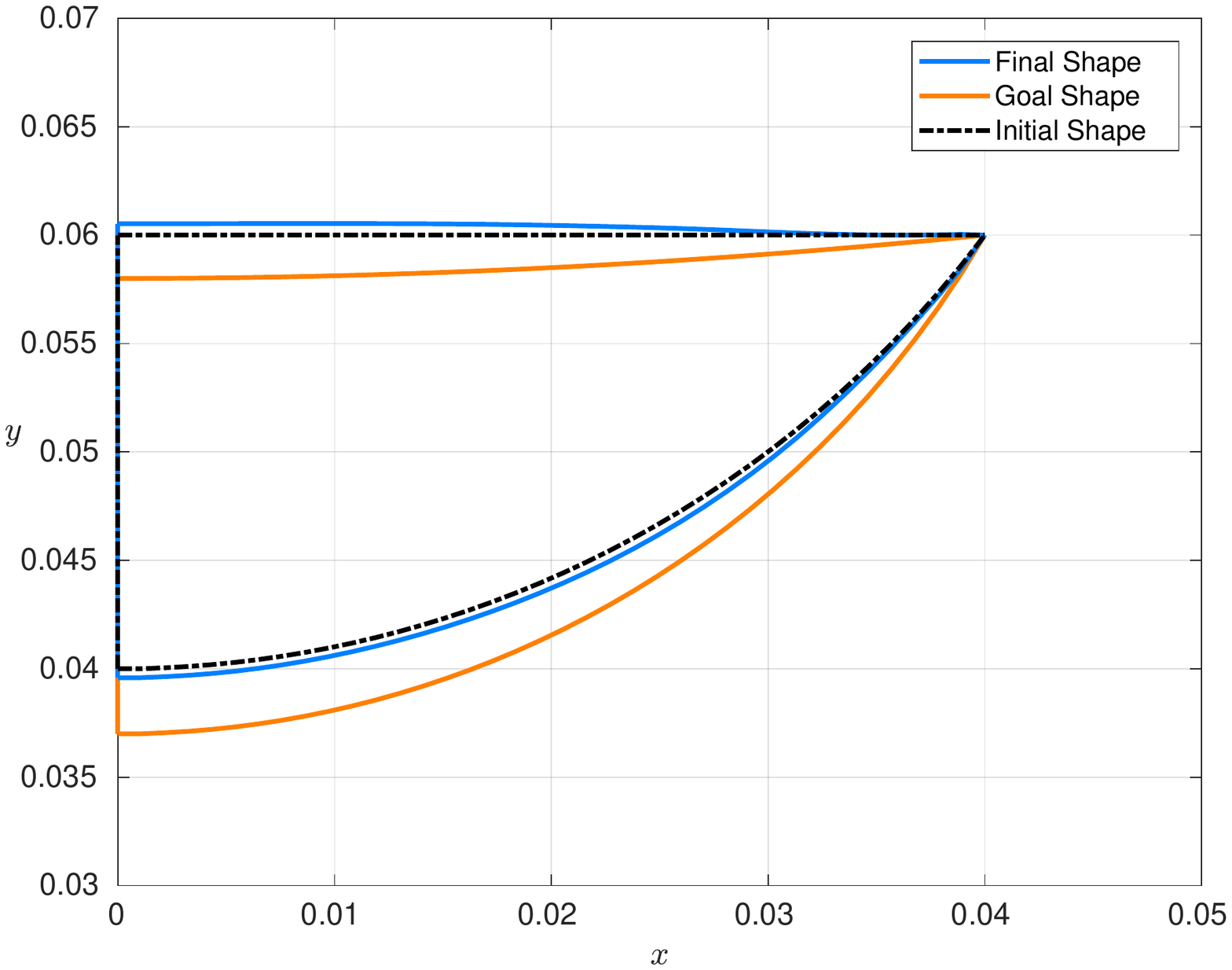}
	\includegraphics[trim = 4cm 9cm 4.2cm 9cm, clip, scale=0.48]{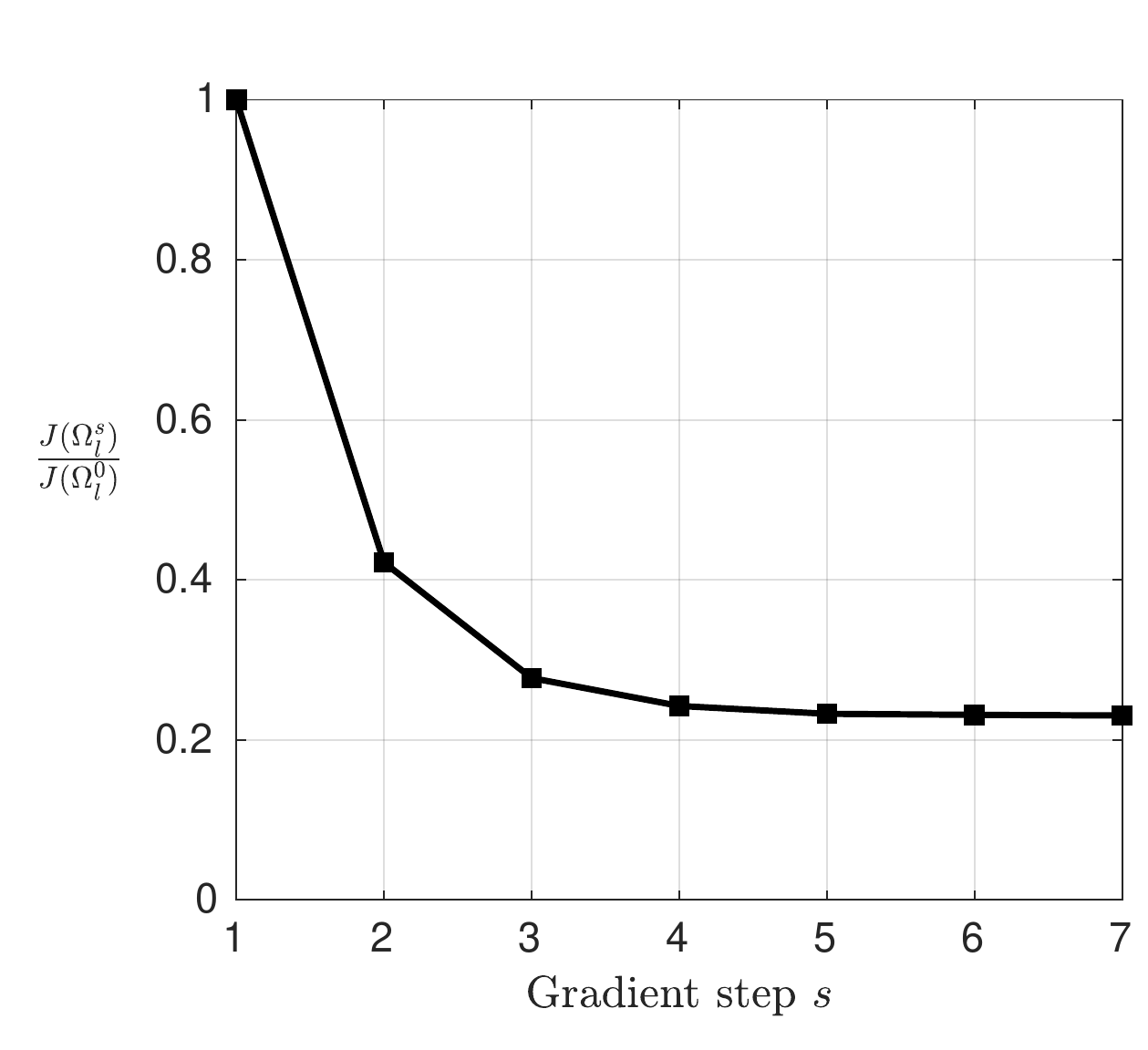}
	\caption{\textbf{left:} Initial, final and goal shape; the final shape is a local, but not a global minimum. \textbf{right:} Relative cost decay. \label{fig:local_minimum}}
\end{figure}

The third experiment with artificial data demonstrates that the final shape obtained with our algorithm does not necessarily have to be a global minimum. We start with the initial lens with a straight upper boundary 
 $$\Omega_l^0 = \Omega_{l,\textup{upper straight}},$$ and take the target lens to be $$\Omega_l^\star = \Omega_{l,\textup{both down}},$$ see Table~\ref{table_domain_size}.
Ideally, both lens boundaries should move down during optimization in order to reach the target shape. However, it seems that our initial lens is already close to a local minimum. The final shape is similar to the initial one and represents a local, but not a global minimum. However, the decrease in cost is still significant, going from $J(\Omega_l^0)=1.388630\cdot 10^4 $ to $J(\Omega_l^{\textup{final}})=3.200507 \cdot 10^3$, which amounts to a reduction of about $76.95\,$\% of the initial value. The experiment was conducted with the linear NURBS. The lens shapes and the relative change of the cost function are given in Figure \ref{fig:local_minimum}. 
\subsection{Test 2: Unknown solution}
\begin{figure}[h]
	\includegraphics[trim = 3cm 7cm 3cm 7cm, clip, scale=0.37]{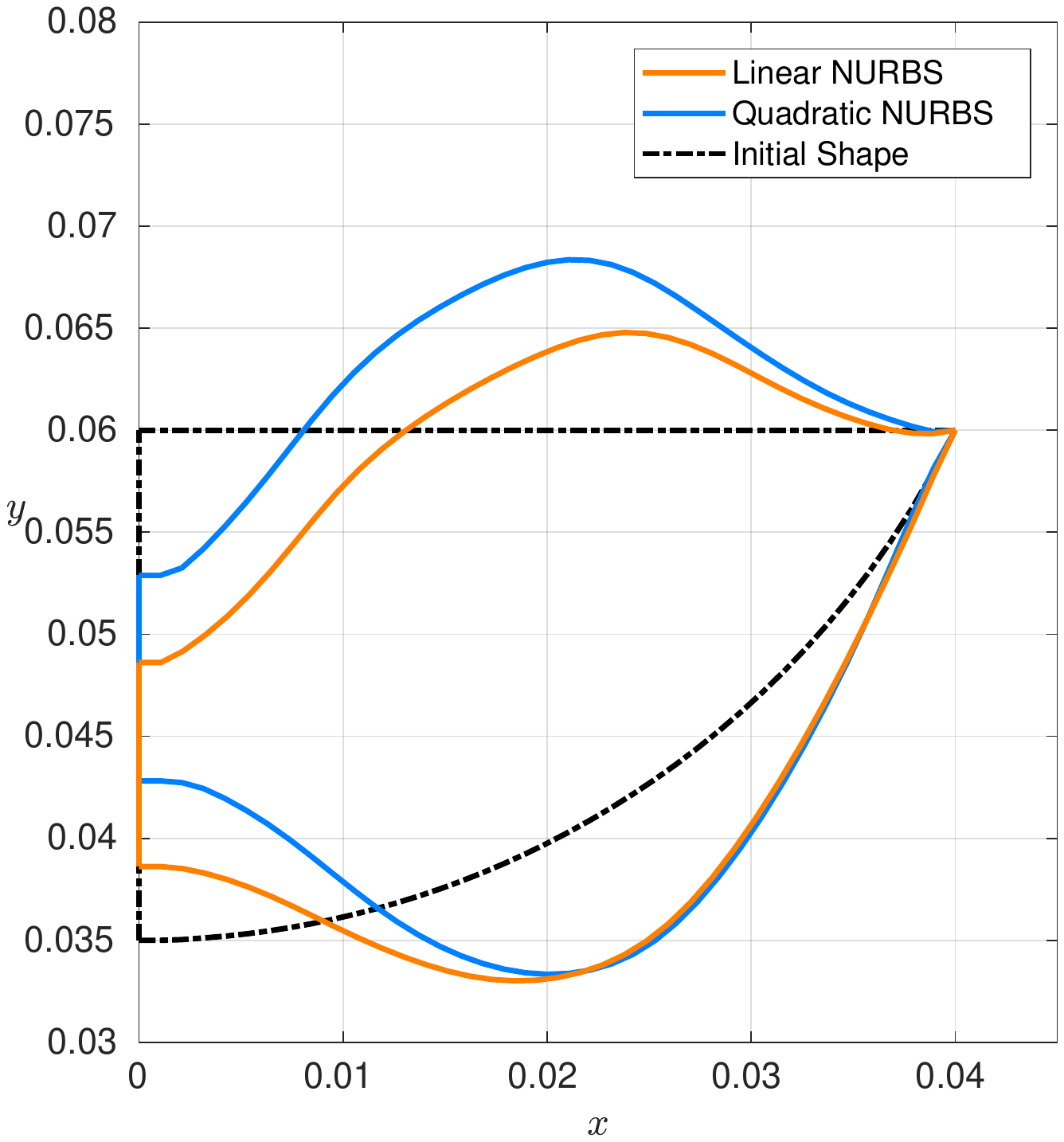}
	\includegraphics[trim = 3cm 9cm 4cm 9.9cm, clip, scale=0.5]{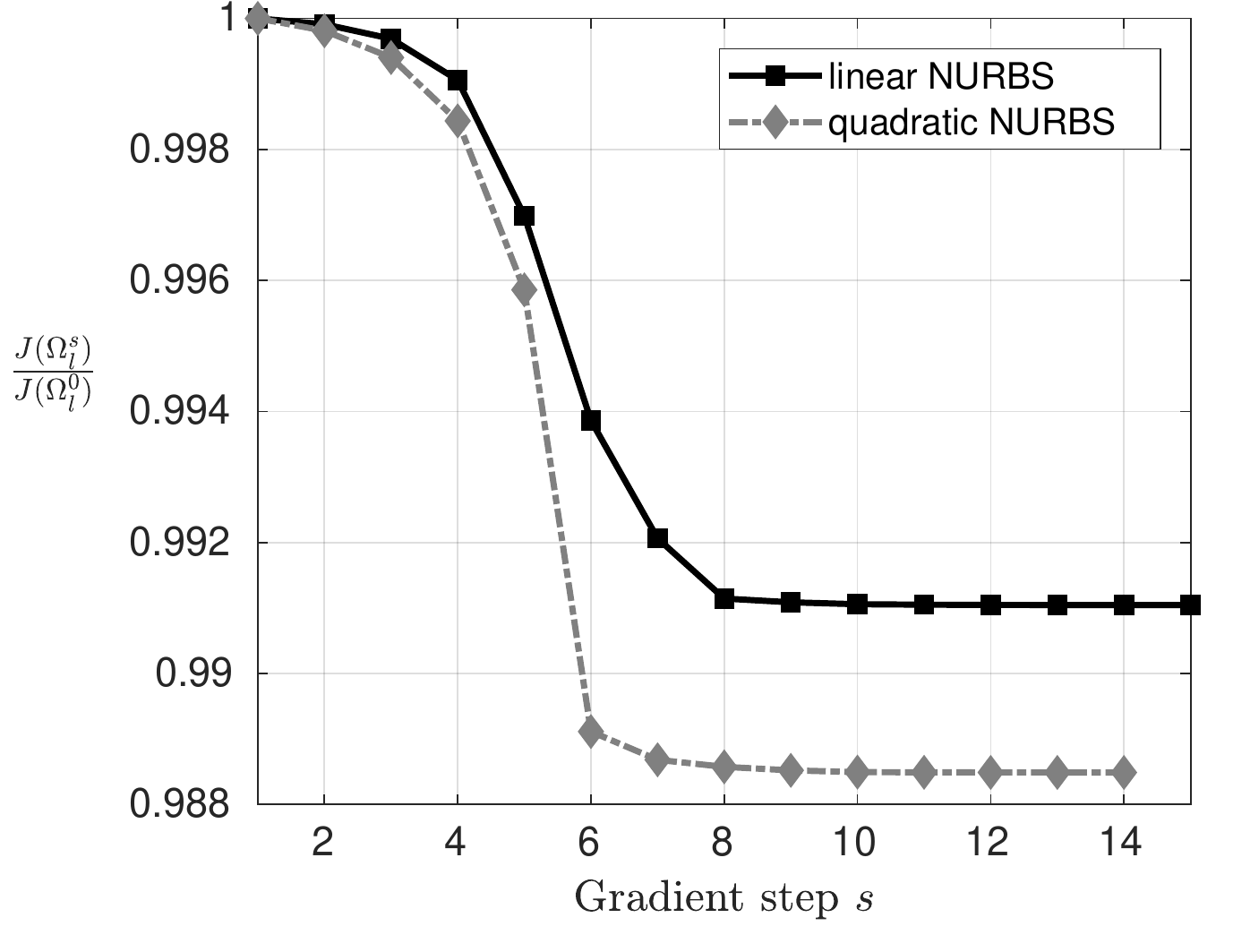}
	\caption{\textbf{left:} Initial and final lens shape. \textbf{right:} Relative cost decay versus the number of gradient steps.\label{fig:Gauss_Result}}
\end{figure}

\begin{figure}[h]
	\includegraphics[trim = 0cm 0cm 9cm 0cm, clip, scale=0.20]{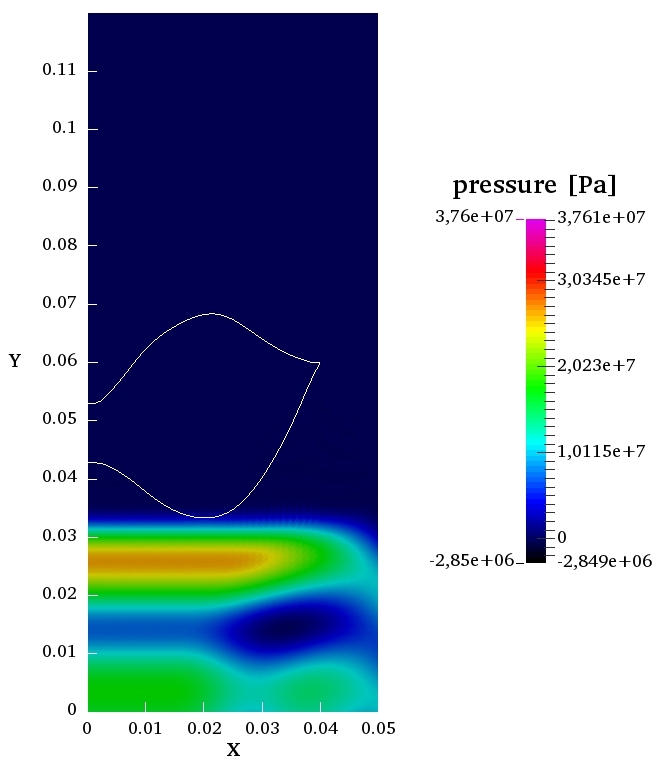}
	\includegraphics[trim = 0cm 0cm 9cm 0cm, clip, scale=0.20]{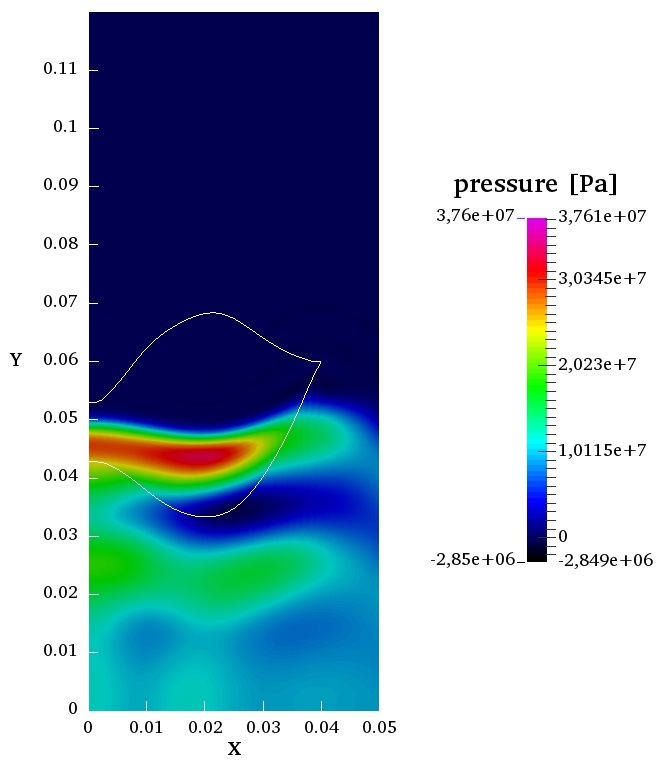}
	\includegraphics[trim = 0cm 0cm 0cm 0cm, clip, scale=0.20]{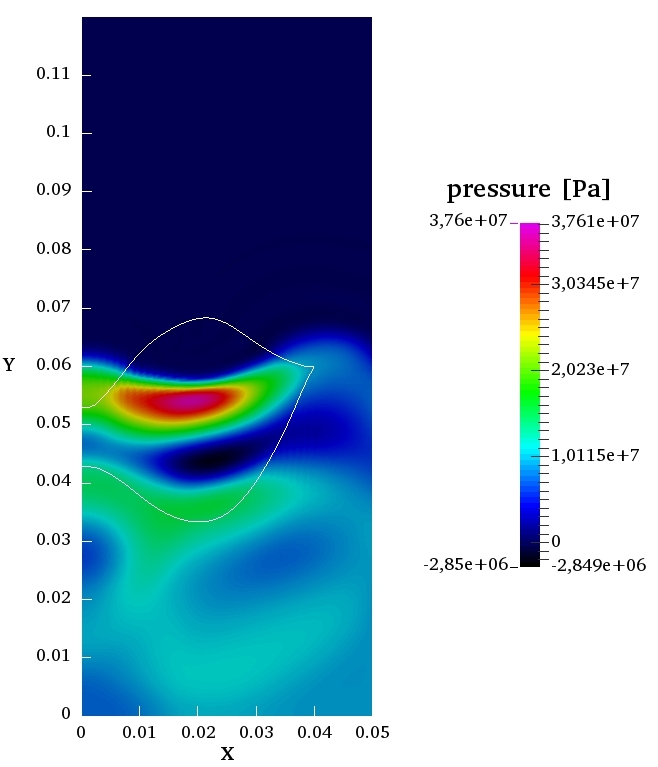}\\
	\includegraphics[trim = 0cm 0cm 9cm 0cm, clip, scale=0.20]{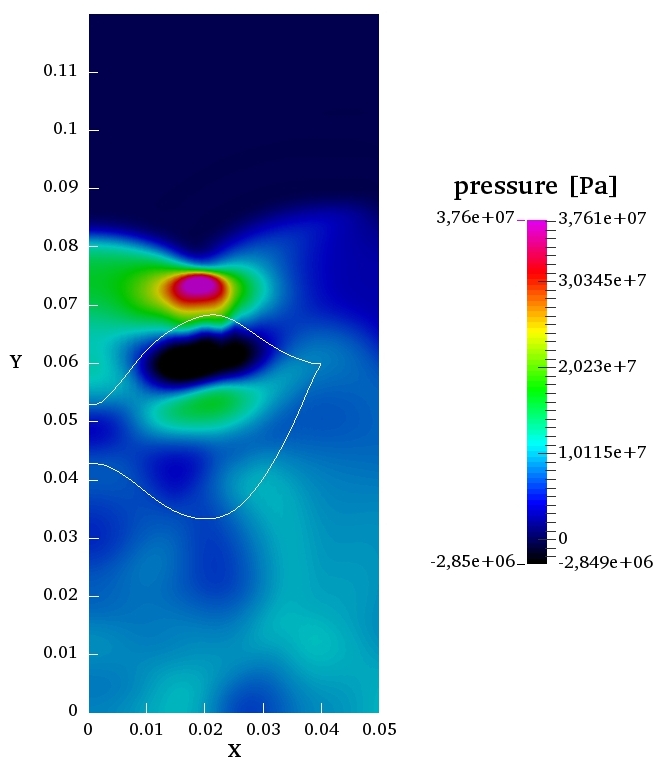}
	\includegraphics[trim = 0cm 0cm 9cm 0cm, clip, scale=0.20]{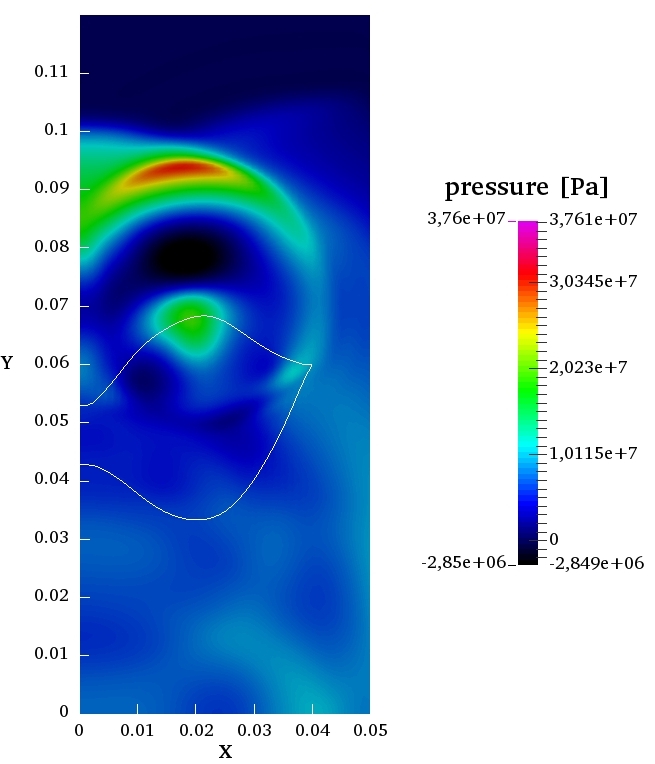}
	\includegraphics[trim = 0cm 0cm 0cm 0cm, clip, scale=0.20]{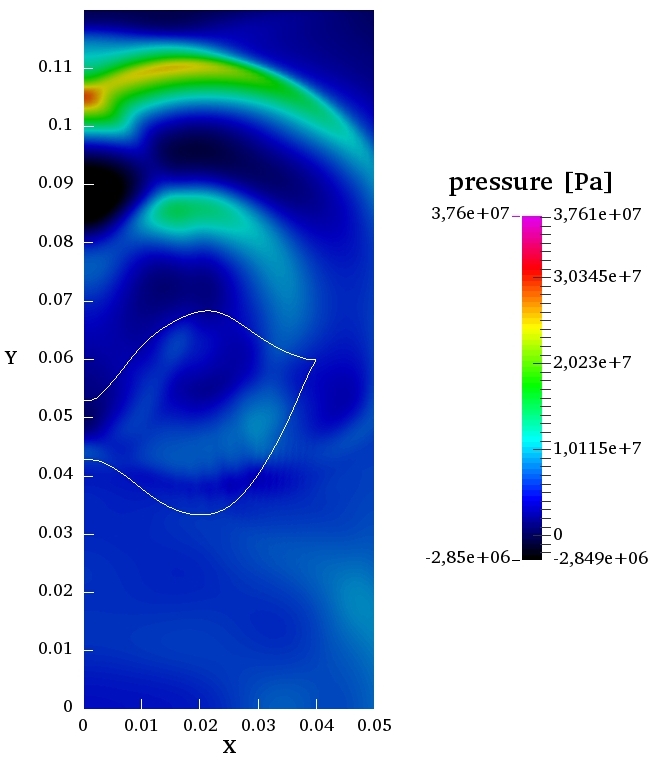}\\
	\caption{Pressure-wave propagation with the final lens shape using quadratic NURBS\label{fig:Final_Lens_Run}}
\end{figure}

\begin{figure}[h]
	\includegraphics[trim = 3cm 7cm 3cm 7cm, clip, scale=0.37]{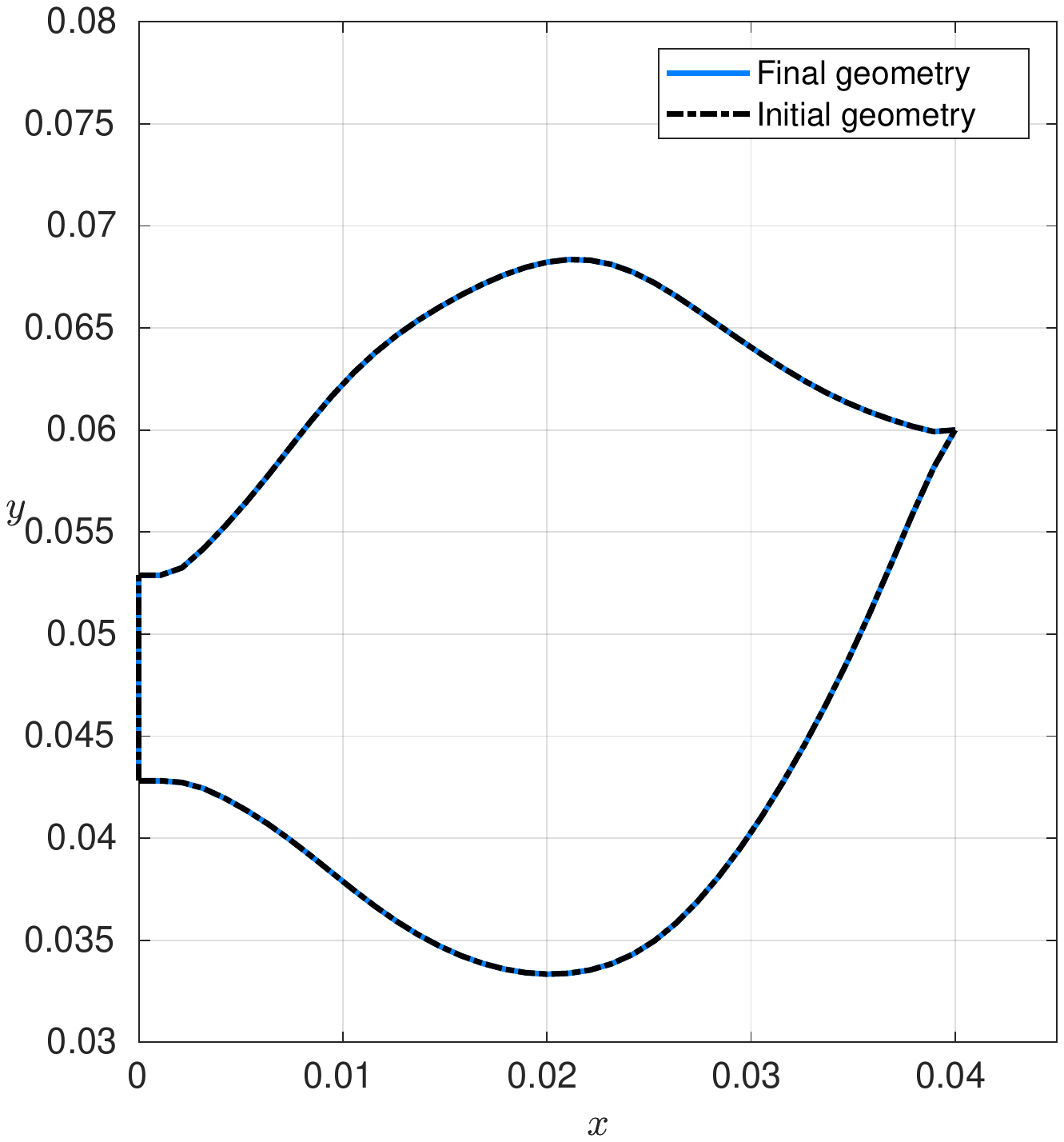}
	\includegraphics[trim = 3cm 9cm 4cm 9.9cm, clip, scale=0.5]{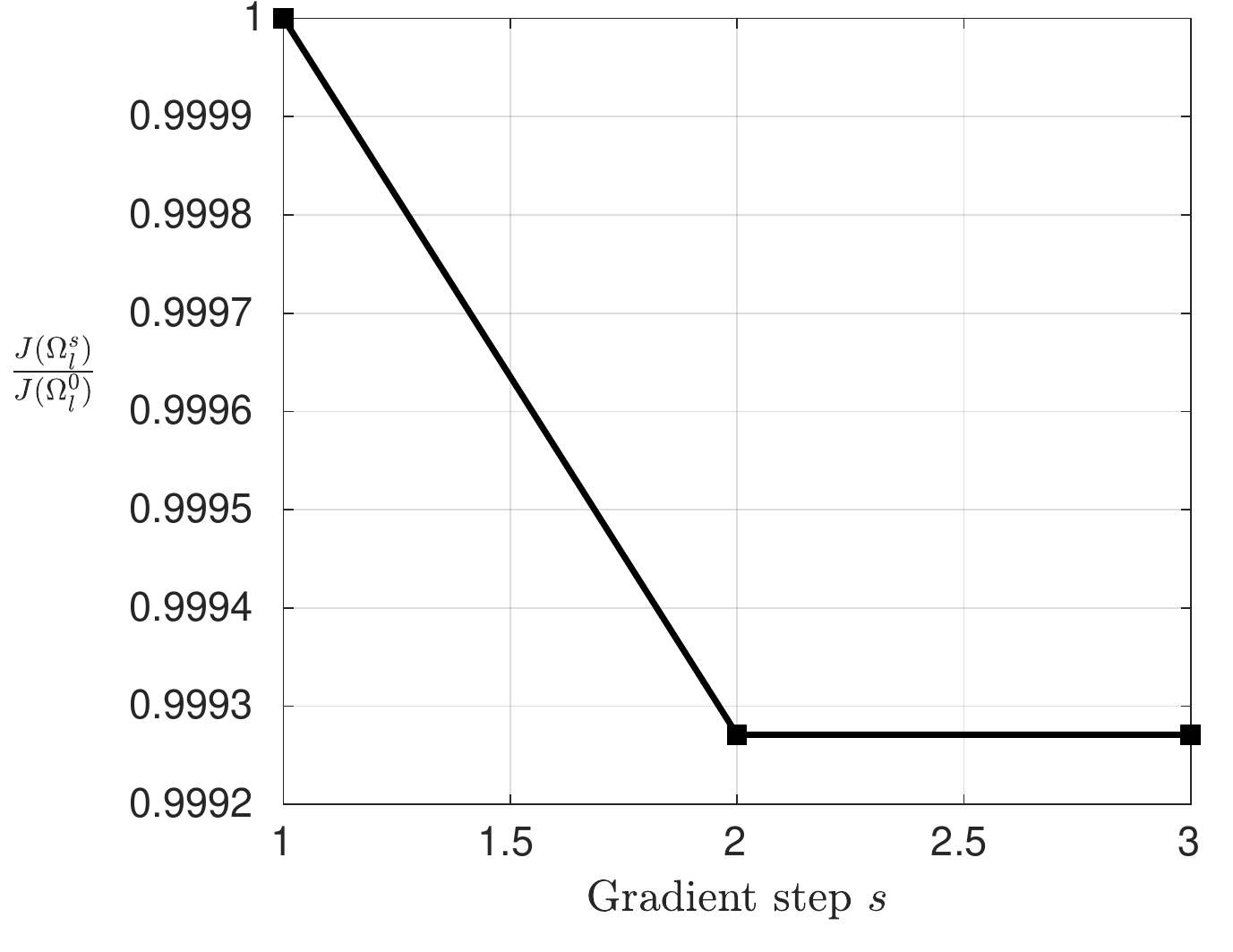}
	\caption{\textbf{left:} No visible changes between the initial and final lens in the post-processing optimization. \textbf{right:} Cost decay by $0.0729\,\%$ of its starting value.  \label{fig:PPO}}
\end{figure}

We now proceed to the case where the goal shape is \emph{not known}. We assume the desired focal point to be on the axis of symmetry with coordinates $(0,y_{\textup{fp}})^T$ and associate with it the desired pressure distribution $u_d$. \\
\indent Finding an adequate analytical expression for the nonlinear target pressure distribution is challenging. Ideally, $u_d$ should assign high pressure values to the region around the focal point and low values elsewhere. We will employ a stationary Gaussian centered around $(0,y_{\textup{fp}})^T$ with a scaling factor in the order of magnitude of the maximal pressure amplitude within the computational domain:
\begin{align*}
	u_d(x,y,t)=A\cdot \exp\left(-\frac{1}{2}(x,y-y_{\textup{fp}})\cdot\varSigma^{-1}\cdot(x,y-y_{\textup{fp}})^T\right).
\end{align*}
Here $\varSigma$ denotes the covariance matrix $$\varSigma = \textup{diag}(\sigma_x^2,\sigma_y^2),$$ with parameters $\sigma_x = 0.02, \sigma_y =  0.004 $. The $y$-coordinate of the focal point is given by $y_{\textup{fp}}=0.105$ and the amplitude of the Gaussian by $A=6\cdot 10^7$. \\
\indent Having in mind a realistic shape for the focusing lens, we will include an additional stopping criterion in the algorithm, that will stop the optimization process if an update would lead to a shape being thinner than some tolerance. \\
\begin{figure}[h]
	\includegraphics[trim = 1cm 7cm 2cm 7cm, clip, scale=0.36]{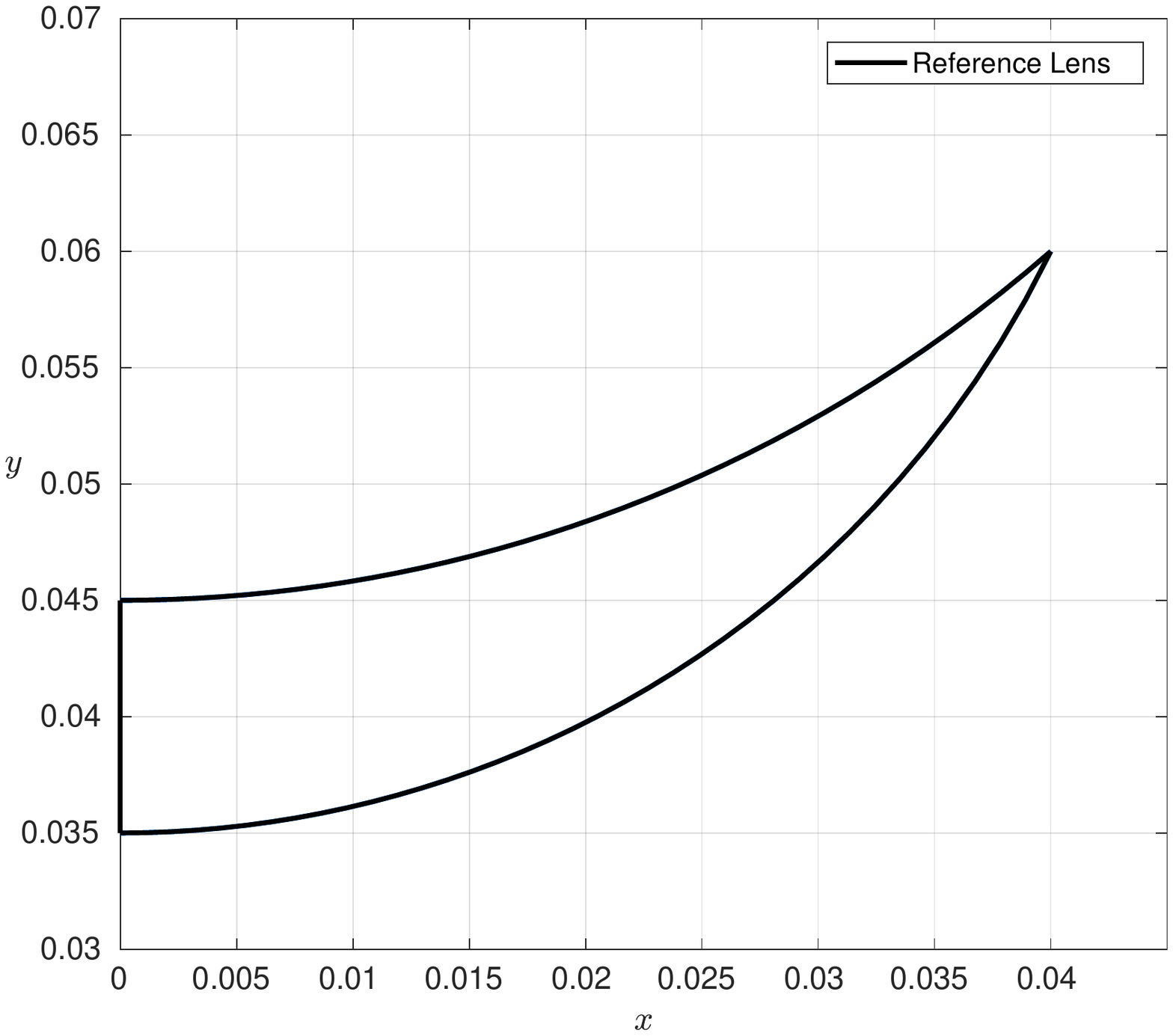}
	\caption{Reference lens that marks the minimal thickness possible to manufacture.\label{fig:Ref_Lens}}
\end{figure}\\
\indent A constant thickness tolerance over the width of the lens would already fail in the first gradient step, since the thickness of the lenses is zero at the tip. To determine a lower thickness limit that varies over the width of the lens, we set up a \emph{reference lens}, which represents the thinnest lens design that is possible to manufacture (see Figure~\ref{fig:Ref_Lens}). If the newly computed lower and upper lens boundary do not fulfill the minimal thickness condition at all control points, we reduce the descent step size until either the minimal thickness limit is matched again for all control points or until the step size undergoes a lower bound, in which case the algorithm terminates.\\
\indent We will use a reference lens which has two arcs as upper and lower boundaries. The arcs have horizontal tangents at the point on the axis of symmetry and a width on the axis of $0.01$, see~Figure~\ref{fig:Ref_Lens}. The thickness constraint is given by
\begin{equation*}
	d_{\min}:[0,0.04]\rightarrow \mathbb{R}, \quad d_{\min}(x)=0.0263-\sqrt{0.0608^2-x^2}+\sqrt{0.0445^2-x^2}.
\end{equation*}
The result of the optimization process with the given method and parameters is shown in Figure \ref{fig:Gauss_Result}. The initial shape is $\Omega_{l,\textup{Gauss}}$, specified in Table \ref{table_domain_size}, the domain $D$ over which the pressure is tracked was chosen to be $D=[0,0.015]\times [0.1,0.11]$. 
\indent Figure~\ref{fig:Gauss_Result} shows the results of the optimization process with linear and quadratic NURBS, where we took the lower bound for the descent step size to be $\textup{TOL}_{\textup{step}}=10^{-3}$. For the employed mesh sizes, we refer again to Table~\ref{spatial_discretization}, where we note that the numbers of degrees of freedom are of the same order of magnitude. In both linear and quadratic NURBS cases, we notice a decrease of the cost functional by $0.89\,\%-1.15\,\%$ of its initial value. The small decrease is likely caused by using a stationary Gaussian as the desired pressure distribution, which our moving pressure wave can hardly fit to. Future research should include other objective functions which would not require an analytically given target pressure distribution, for instance, maximizing acoustic pressure in a certain area. \\
\indent The resulting final shapes, although similar, still differ. We would like to compare which shape performs better. We observe that the quadratic NURBS give a better relative decay of the cost functional (see Figure ~\ref{fig:Gauss_Result}). In order to test if the quadratic NURBS found a better shape, we will bring the two shapes into the same setting. First, we interpolate the final quadratic NURBS shape into the linear NURBS space and compute the cost values. Secondly, we reverse the process and represent the "linear" final shape with quadratic NURBS and compute the cost in this quadratic NURBS space. The resulting values of the cost functionals are listed in Table~\ref{tab:cost_comp}, confirming that the quadratic NURBS perform better in both cases by around $0.235\,\%$. In Figure~\ref{fig:Final_Lens_Run}, snapshots of the pressure wave propagation with the final quadratic NURBS shape are shown. 
\begin{center}
	\renewcommand{\arraystretch}{1.6}
	\captionof{table}{Cost comparison with final lens shapes} \label{tab:cost_comp} 
	\begin{tabular}{| c || c | c||} 
		\hline
		& \multicolumn{2}{c||}{Interpolated into the}\\
		Optimization with & linear NURBS space & quadratic NURBS space \\
		\hline\hline
		linear NURBS & $J=2.937255\cdot 10^7$& $J=2.933371\cdot 10^7$ \\
		\hline
		quadratic NURBS & $J= 2.930353\cdot 10^7$ & $J=2.926472\cdot 10^7$ \\
		\hline
	\end{tabular}\vspace{2mm}
\end{center}
Finally, we take the optimal "quadratic" lens (interpolated into the linear NURBS space) as the initial geometry and restart the optimization algorithm in the linear NURBS setting. The results of this post-processing optimization, given in Figure \ref{fig:PPO}, confirm that this lens is already a local minimum. The algorithm performs 3 steps, resulting in a cost reduction of only $0.0729\,$\% of its starting value, with no visible changes to the shape of the lens. The steps that the algorithm needed to take are likely due to the approximation error made by reducing the degree of the NURBS curve.\\
\clearpage

\section{Conclusion and outlook}
We have developed a gradient-based algorithm to find a locally optimal shape of an acoustic lens in nonlinear ultrasound applications. The shape derivative of the pressure-tracking cost functional was first rigorously derived by employing a variational approach. The algorithm was then numerically realized within the isogeometric concept, where geometries can be exactly represented in every gradient step. Due to the nonlinear nature of the problem and high source frequencies, a very fine mesh and time discretization had to be employed, increasing the computational complexity. Numerical experiments show that the quadratic NURBS perform better with respect to the final lens shape than the linear ones, reinforcing the potential of isogeometric shape optimization. \\
\indent There are many possible directions for further research. Combining linear and higher-order NURBS in optimization in a clever way could result in a more efficient algorithm. Future research should also include considerations of different objective functionals, where an analytical expression for the target pressure distribution would not be needed. Exploring other acoustic models for nonlinear ultrasound propagation is of interest as well. Finally, from the point of view of applications, it is important to numerically tackle a 3D setting with even higher source frequencies, which would require developing a high-performance computational framework. 
\section*{Acknowledgments} The authors gratefully acknowledge the funds provided by the Deutsche Forschungsgemeinschaft under the grant numbers WO 671/11-1 and WO 671/15-1 (within the Priority Programme SPP 1748).


\begin{thebibliography}{99}
\bibitem{BuffaVazquez}
L.  Beir\~ ao  da  Veiga,  A.  Buffa,  G.  Sangalli  and  R.  V\' azquez, \emph{Mathematical  analysis  of  variational
isogeometric methods}, Acta Numerica, 23(5), 157--287, 2014.
\bibitem{Blana}
A. Blana, N. Walter, S. Rogenhofer and W. F. Wieland, \emph{High-intensity focused ultrasound for the treatment of localized prostate cancer: 5-year experience}, Urology, 63(2), 297--300, 2004.
\bibitem{Brandenburg}
C. Brandenburg, F. Lindemann, M. Ulbrich and S. Ulbrich, \emph{Advanced numerical methods for PDE constrained optimization with application to optimal design in Navier Stokes flow},  Constrained Optimization and Optimal Control for Partial Differential Equations, Springer Basel, 257--275, 2012.
\bibitem{Brereton}
G. J. Brereton and B. A. Bruno, \emph{Particle removal by focused ultrasound}, Journal of sound and vibration, 173(5), 683--698, 1994.
\bibitem{Canney}
M. S. Canney, Bailey,  L. A. Crum, V. A. Khokhlova and  O. A. Sapozhnikov, \emph{ Acoustic characterization of high intensity focused ultrasound fields: A combined measurement and modeling approach}, The Journal of the Acoustical Society of America, 124(4), 2406--2420, 2008.
\bibitem{Cho}
S. Cho and S.-H. Ha, \emph{Isogeometric shape design optimization: exact geometry and enhanced sensitivity}, Structural and Multidisciplinary Optimization, 38(1), 53--70, 2009.
\bibitem{ChungHulbert}
J. Chung and G. M. Hulbert, \emph{A time integration algorithm for structural dynamics with improved numerical dissipation: the generalized-$\alpha$ method}, Journal of applied mechanics, 60(2), 371--375, 1993.
\bibitem{ColtonKress}
D. Colton and R. Kress, \emph{Integral Equation Methods in Scattering Theory}, Wiley, New York, 1983.
\bibitem{hughes:09}
J.A. Cottrell, T.J.R. Hughes and Y. Bazilevs, \emph{Isogeometric analysis: toward integration of CAD and FEA}, John Wiley \& Sons, 2009.
\bibitem{Crighton}
D. G. Crighton, \emph{Model equations of nonlinear acoustics}, Annual Review of Fluid Mechanics, 11(1), 11--33, 1979.
\bibitem{Demengel}
F. Demengel, G. Demengel and R. Erné,  \emph{Functional spaces for the theory of elliptic partial differential equations}, London, UK: Springer, 2012.
\bibitem{Dogan}
G. Dogan, P. Morin, R. H. Nochetto and M. Verani, \emph{Discrete gradient flows for shape optimization and applications}, Computer Methods in Applied Mechanics and Engineering, 196(37), 3898--3914, 2007.
\bibitem{Engquist}
B. Engquist and A. Majda, \emph{Absorbing Boundary Conditions for the Numerical Evaluation of Waves}, Proceedings of the National Academy of Sciences, 74(5), 1765--1766, 1977.
\bibitem{Eppler1}
K. Eppler and H. Harbrecht, \emph{Coupling of FEM and BEM in shape optimization}, Numerische Mathematik, 104(1), 47--68, 2006.
\bibitem{Eppler2}
K. Eppler, H. Harbrecht and R. Schneider, \emph{On convergence in elliptic shape optimization}, SIAM Journal on Control and Optimization, 46(1), 61--83, 2007.
\bibitem{Erlicher}
S. Erlicher, L. Bonaventura and O.S. Bursi, \emph{The analysis of the Generalized-$\alpha$ method for nonlinear dynamic problems}, Computational Mechanics 28, 2002.
\bibitem{Falco}
C. de Falco, A. Reali and R. V\' azquez,  \emph{GeoPDEs: a research tool for Isogeometric Analysis of PDEs}, Advances in Engineering Software, 42(12), 1020--1034, 2011. 
\bibitem{Folds}
D. L. Folds, \emph{Speed of sound and transmission loss in silicone rubbers at ultrasonic frequencies}, The Journal of the Acoustical Society of America, 56(4), 1295--1296, 1974.
\bibitem{Simeon}
D. Fu\ss eder, A.-V. Vuong and B. Simeon, \emph{Fundamental aspects of shape optimization in the context of isogeometric analysis}, Computer Methods in Applied Mechanics and Engineering, 286(39), 313--331, 2015.
\bibitem{Simeon2}
D. Fu\ss eder and  B. Simeon. \emph{Algorithmic Aspects of Isogeometric Shape Optimization}, In B. J\"uttler and B. Simeon (editors), Isogeometric Analysis and Applications 2014, 183--207, Springer, 2015.
\bibitem{Gangl}
P. Gangl, U. Langer, A. Laurain, H. Meftahi and K. Sturm, \emph{Shape optimization of an electric motor subject to nonlinear magnetostatics}, SIAM Journal on Scientific Computing, 37(6), B1002--B1025, 2015.
\bibitem{Harbrecht1}
H. Harbrecht, \emph{Analytical and numerical methods in shape optimization}, Mathematical Methods in the Applied Sciences, 31(18), 2095--2114, 2008.
\bibitem{Haslinger}
J. Haslinger and P. Neittaanm\" aki, \emph{Finite element approximation for optimal shape, material, and topology design}, John Wiley \& Sons, 1996.
\bibitem{Henrot}
A. Henrot and M. Pierre, \emph{Variation et optimisation de formes: une analyse g\'eom\'etrique}, Springer Science \& Business Media, 48, 2006.
\bibitem{Hoellig}
K. H\"ollig, \emph{Finite Element Methods with B-Splines}, SIAM, 2003.
\bibitem{Hofmann}
S. Hofmann, M. Mitrea and M. Taylor, \emph{Geometric and transformational properties of Lipschitz domains, Semmes-Kenig-Toro Domains, and other classes of finite perimeter domains}, The Journal of Geometric Analysis, 17(4), 593--647, 2007.
\bibitem{hughes:05}
T.J.R. Hughes, J.A. Cottrell and Y. Bazilevs, \emph{Isogeometric analysis: CAD, finite elements, NURBS, exact geometry and mesh refinement}, Computer Methods in Applied Mechanics and Engineering, 194, 4135--4195, 2005. 
\bibitem{Hoffelner}
J. Hoffelner, H. Landes, M. Kaltenbacher and R. Lerch, \emph{Finite element simulation of nonlinear wave propagation in thermoviscous fluids including dissipation}, IEEE transactions on ultrasonics, ferroelectrics, and frequency control, 48(3), 779--786, 2001.
\bibitem{Maloney}
E. Maloney, and J. H. Hwang, \emph{Emerging HIFU applications in cancer therapy}, International Journal of Hyperthermia, 31(3), 302-309, 2015.
\bibitem{IKP}
K.~Ito, K.~Kunisch and G. Peichl, \emph{Variational approach to shape derivatives for a class
of Bernoulli problems}, Journal of Mathematical Analysis and Applications, 314(1), 126--149, 2006.
\bibitem{IKP08}
K.~Ito, K.~Kunisch and G. Peichl, \emph{Variational approach to shape derivatives}, ESAIM: Control, Optimisation and Calculus of Variations, 14(3), 517--539, 2008.
\bibitem{KL11}
B.\ Kaltenbacher and I.\ Lasiecka, \emph{Well-posedness of the Westervelt and the Kuznetsov equation with nonhomogeneous Neumann boundary conditions}, AIMS Proceedings, 763--773, 2011.
\bibitem{KLV10}
B.\ Kaltenbacher, I.\ Lasiecka and S.\ Veljovi\'c, \emph{Well-posedness and exponential decay for the Westervelt equation with inhomogeneous Dirichlet boundary data}, Progress in Nonlinear Differential Equations and Their Applications, 60, 357--387, 2011.
\bibitem{KaltenbacherPeichl}
B.~Kaltenbacher and G.~Peichl, \emph{Sensitivity analysis for a shape optimization problem in lithotripsy}, Evolution Equations and Control Theory(EECT), 5, 399--430, 2016.
\bibitem{KV}
B. Kaltenbacher and S. Veljovi\' c, \emph{Sensitivity analysis of linear and nonlinear lithotripter models},
European Journal of Applied Mathematics, 22(1), 21--43, 2010. 
\bibitem{manfred}
M.~Kaltenbacher, \emph{Numerical simulations of mechatronic sensors and actuators}, Springer, Berlin, 2004.
\bibitem{Kennedy1}
J. E. Kennedy, \emph{High-intensity focused ultrasound in the treatment of solid tumors}, Nature reviews cancer, 5(4), 321--327, 2005.
\bibitem{Kennedy2}
J. E. Kennedy, F. Wu, G. R. Ter Haar, F. V. Gleeson, R. R. Phillips, M. R. Middleton and D. Cranston, \emph{High-intensity focused ultrasound for the treatment of liver tumors}, Ultrasonics, 42(1), 931--935, 2004.
\bibitem{Kiendl}
J. Kiendl, R. Schmidt, R. W\"uchner and K.-U. Bletzinger, \emph{Isogeometric shape optimization of shells using semi-analytical sensitivity analysis and sensitivity weighting}, Computer Methods in Applied Mechanics and Engineering, 274(1), 148--167, 2014.
\bibitem{KuhlCrisfield}
D. Kuhl and M.A. Crisfield, \emph{Energy-conserving and decaying algorithms in nonlinear structural mechanics}, International Journal for Numerical Methods in Engineering, 45, 569--599, 1999.
\bibitem{Laporte}
E. Laporte and P. Le Tallec, \emph{Numerical methods in sensitivity analysis and shape optimization}, Springer Science \& Business Media, 2012.
\bibitem{Lee}
Y.-S. Lee, \emph{Numerical solution of the KZK equation for pulsed finite amplitude sound beams in thermoviscous fluids}, PhD Thesis, The University of Texas at Austin, 1993.
\bibitem{Lee2}
D. Lee, N. Koizumi, K. Ota, S. Yoshizawa, A. Ito, Y. Kaneko, Y. Matsumoto and M. Mitsuishi, \emph{Ultrasound-based visual serving system for lithotripsy}, Intelligent Robots and Systems, 877--882, 2007.
\bibitem{Muench}
A. M\" unch, \emph{Optimal design of the support of the control for the 2-D wave equation: a numerical method}, International Journal of Numerical Analysis and Modeling, 5(2), 331--351, 2008.
\bibitem{Maestre}
F. Maestre, A. M\" unch and P. Pedregal, \emph{A spatio-temporal design problem for a damped wave equation}, SIAM Journal on Applied Mathematics, 68(1), 109--132, 2007.
\bibitem{Mancini}
J.G. Mancini, A. Neisius, N. Smith, G. Sankin, G. M. Astroza, M. E. Lipkin, W. N. Simmons, G. M. Preminger and P. Zhong, \emph{Assessment of a modified acoustic lens for electromagnetic shock wave lithotripters in a swine model}, The Journal of urology, 190(3), 1096--1101, 2013.
\bibitem{MeyerWilke}
S. Meyer and M. Wilke, \emph{Optimal regularity and long-time behavior of solutions for the Westervelt equation}, Applied Mathematics \& Optimization, 64(2), 257--271, 2011.
\bibitem{MuratSimon}
F. Murat and S. Simon, \emph{Etudes de problems d'optimal design}, Lecture notes in computer science, 41, 54--62, 1976.
\bibitem{Neisius}
A. Neisius, N.B. Smith, G. Sankin, N.J. Kuntz, J.F. Madden, D.E. Fovargue, S. Mitran, M.E. Lipkin, W.N. Simmons, G.M. Preminger and P. Zhong, \emph{Improving the lens design and performance of a contemporary electromagnetic shock wave lithotripter}, Proceedings of the National Academy of Sciences, 111(13), E1167--E1175, 2014.
\bibitem{Newmark}
N. M. Newmark, \emph{A Method of Computation for Structural Dynamics}, Journal of Engineering Mechanics, ASCE, 85, 67--94, 1959.
\bibitem{Nguyen}
D. M. Nguyen, A. Evgrafov and J. Gravesen,\emph{Isogeometric shape optimization for scattering problems}, Progress In Electromagnetics Research B, 45, 117--146, 2012.
\bibitem{NikolicKaltenbacher}
V. Nikoli\' c and B. Kaltenbacher, \emph{Sensitivity analysis for shape optimization of a focusing acoustic lens in lithotripsy}, Applied Mathematics and Optimization, 76(2), 261--301, 2017.
\bibitem{Nortoft}
P. Nortoft and J. Gravesen, \emph{Isogeometric shape optimization in fluid mechanics}, Struct Multidisc Optim, 48, 909--925, 2013.
\bibitem{Paganini}
A. Paganini, \emph{Numerical shape optimization with finite elements}, PhD Thesis, ETH Z\" urich, 2016.
 \bibitem{Park}
S. H. Park, J. W. Yoon, D. Y. Yang and Y. H. Kim, \emph{Optimum blank design in sheet metal forming by the deformation path iteration method}, International Journal of Mechanical Sciences, 41(10), 1217--1232, 1999.
\bibitem{Paterson}
R. F. Paterson, E. Barret, T. M. Siqueira, T. A. Gardner, J. Tavakkoli, V. V. Rao, N. T. Sanghvi, L. Cheng and A. L. Shalhav, \emph{Laparoscopic partial kidney ablation with high intensity focused ultrasound}, The Journal of urology, 169(1), 347--351, 2003.
\bibitem{Piegl}
L. Piegl and W. Tiller, \emph{The NURBS Book}, Springer, 1997.
\bibitem{Qian}
X. Qian and O. Sigmund, \emph{Isogeometric shape optimization of photonic crystals via Coons patches}, Computer Methods in Applied Mechanics and Engineering, 200(25), 2237--2255, 2011.
\bibitem{Rossing}
T. D. Rossing (Ed.), \emph{Springer Handbook of Acoustics}, Springer, 2014.
\bibitem{Shumaker}
L. Schumaker, \emph{Spline Functions: Basic Theory}, 3rd Edition, Cambridge University Press, Cambridge, 2007.
\bibitem{Veljovic}
S. Veljovi\'{c}, \emph{Shape optimization and optimal boundary control for high intensity focused ultrasound (HIFU)}, PhD Thesis, University of Erlangen-Nuremberg, 2009. 
\bibitem{Wall}
W. A. Wall, M. A. Frenzel and C. Cyron, \emph{Isogeometric structural shape optimization}, Computer Methods in Applied Mechanics and Engineering, 197(33), 2976--2988, 2008.
\bibitem{Wu}
F. Wu, W. Z. Chen, J. Bai, J. Z. Zou, Z. L. Wang, H. Zhu and Z. B. Wang, \emph{Pathological changes in human malignant carcinoma treated with high-intensity focused ultrasound}, Ultrasound in medicine \& biology, 27(8), 1099--1106, 2001.
\bibitem{Yoshizawa}
S. Yoshizawa, T. Ikeda, A. Ito, R. Ota, S. Takagi and Y. Matsumoto, \emph{High intensity focused ultrasound lithotripsy with cavitating microbubbles}, Medical \& biological engineering \& computing, 47(8), 851--860, 2009.
\bibitem{Zhong}
P. Zhong, N. Smith, N. W. Simmons and G. Sankin, \emph{A new acoustic lens design for electromagnetic shock wave lithotripters}, AIP Conference Proceedings, 1359(1), 42--47, 2011.
\bibitem{Zolesio}
J.-P. Zol\' esio and M. C. Delfour, \emph{Shapes and Geometries: Metrics, Analysis, Differential Calculus, and Optimization}, SIAM, 2011. \\
\end{thebibliography}
\end{document}